\documentclass{amsbook}
\usepackage{amssymb, euscript,amsmidx}
\makeindex{a}
\makeindex{b}

\def\Cal{\mathcal}
\def\A{{\Cal A}}
\def\B{{\Cal B}}
\def\H{{\Cal H}}
\def\X{{\fr X}}

\def\M{{\Cal M}}
\def\P{{\Cal P}}
\def\D{{\EuScript{D}}}
\def\DP{{\Cal D}}
\def\S{{\Cal S}}
\def\F{{\Cal F}}

\def\I{{\Cal I}}
\def\L{{\Cal L}}
\def\T{{\Cal T}}
\def\W {{\EuScript W_m}}

\def\J{{\Cal J}}

\def\s{\EuScript{S}}

\def\d{\partial}
\def\tr{{\hbox{\rm tr}}}

\def\Ma{\frM_{n,m}}
\def\Mk{\frM_{n,m}^{(k)}}
\def\Mj{\frM_{n,m}^{(j)}}
\def\Mkc{\overline{\frM}_{n,m}^{(k)}}
\def\Mm{\frM_{n,m}^{(m)}}
\def\Mt{\frM_{n-k,m}}
\def\Mkm{\frM_{k,m}}
\def\Mlm{\frM_{\ell,m}}
\def\Mmk{\frM_{m,k}}
\def\nf{\|f\|_p}
\def\V{{\bf W}_{n,m}}
\def\G{\mathcal{G}}
\def\Z{\mathcal{Z}}
\def \A{\frA}
\def\gnk{G_{n,k}}
\def\cgnk{\Cal G_{n,k}}

\def\ft{\hat f(\tau)}
\def \vx{\check \vp (x)}
\def\f0{f_0}
\def\Fc0{\varphi_0}
\def\rn{\bbr^n}

\def\rnk{\bbr^{n-k}}
\def\irn{\intl_{\bbr^n}}

\def\I_k {I_{-}^{k/2}}
\def\I+k {I_{+}^{k/2}}

\def\icgr{\intl_{\cgnk}}
\def\vnk{V_{n,n-k}}
\def\vnm{V_{n,m}}
\def\cd{\vnk\times\frM_{n-k,m}}
\def\Gr{\frT}
\def\sigk{\sig_{n,n-k}}

\def\bbr{{\Bbb R}}

\def\bbn{{\Bbb N}}

\def\bbc{{\Bbb C}}
\def\bbz{{\Bbb Z}}
\def\bbd{{\Bbb D}}

\def\rank{{\hbox{\rm rank}}}
\def\diag{{\hbox{\rm diag}}}

\def\supp{{\hbox{\rm supp}}}

\def\tr{{\hbox{\rm tr}}}

\def\det{{\hbox{\rm det}}}

\def\min{{\hbox{\rm min}}}

\def\Pr{{\hbox{\rm Pr}}}
\def\dist{{\hbox{\rm dist}}}

\def\fc{\varphi(\xi,t)}
\def\gnk{G_{n,k}}

\def\lp{L^p(\frM_{n,m})}

\def\part{\partial}
\def\intl{\int\limits}
\def\b{\beta}

\def\Gam{\Gamma}
\def\Om{\Omega}
\def\a{\alpha}
\def\om{\omega}

\def\Del{\Delta}
\def\del{\delta}
\def\vp{\varphi}

\def\g{\gamma}
\def\gam{\gamma}
\def\Lam{\Lambda}
\def\sig{\sigma}
\def\lam{\lambda}
\def\z{\zeta}
\def\e{\varepsilon}
\def\t{\tau}
\def\eq{\xi 'x=t}
\def\eqv{\bar\xi '\bar x=\bar t}
\def\Rnm{\bbr^{nm}}

\def\Rnkm{\bbr^{(n-k)m}}
\def\rf{\hat f (\xi,t)}
\def\df{\check \varphi (x)}

\font\frak=eufm10

\def\fr#1{\hbox{\frak #1}}
\def\frA{\fr{A}}

\def\frL{\fr{L}}        
\def\frM{\fr{M}}

\def\frT{\fr{T}}

\def\pk{\P_k}

\def\vmk{V_{m,k}}
\def\gma{\Gamma_m(\a)}

\def\gk{\Gamma_k}
\def\lin{{\hbox{\rm lin}}}

\def\det{{\hbox{\rm det}}}
\def\sdx{{\hbox{\rm sgn$\,$det($x$)}}}
\def\sdy{{\hbox{\rm sgn$\,$det($y$)}}}
\def\sdv{{\hbox{\rm sgn$\,$det($v$)}}}
\def\sgn{{\hbox{\rm sgn}}}

\def\min{{\hbox{\rm min}}}

\def\p{\P_m}
\def\gm{\Gamma_m}
\def\tr{{\hbox{\rm tr}}}
\def\part{\partial}
\def\intl{\int\limits}
\def\b{\beta}

\def\Gam{\Gamma}
\def\Om{\Omega}
\def\a{\alpha}
\def\cpm{\overline\P_m}
\def\pd{\stackrel{*}{P}\!{}^\a}
\def\pdz{\stackrel{*}{P}\!{}^0}
\def\pdl{\stackrel{*}{P}\!{}^\ell}

\newtheorem{theorem}{Theorem}[chapter]
\newtheorem{lemma}[theorem]{Lemma}
\theoremstyle{corollary}
\newtheorem{corollary}[theorem]{Corollary}

\theoremstyle{definition}
\newtheorem{definition}[theorem]{Definition}
\newtheorem{example}[theorem]{Example}

\theoremstyle{remark}
\newtheorem{remark}[theorem]{Remark}

\numberwithin{section}{chapter}
\numberwithin{equation}{chapter}

\newcommand{\be}{\begin{equation}}
\newcommand{\ee}{\end{equation}}

\newcommand{\bea}{\begin{eqnarray}}
\newcommand{\eea}{\end{eqnarray}}
\newcommand{\Bea}{\begin{eqnarray*}}
\newcommand{\Eea}{\end{eqnarray*}}

\begin{document}
\frontmatter
\title{The Radon transform of functions of matrix argument}

\author{Elena Ournycheva}
\address{Institute of Mathematics,
Hebrew University, Jerusalem 91904, $\;$ISRAEL}
\email{ournyce@math.huji.ac.il}

\author{Boris Rubin}
\address{Institute of Mathematics,
Hebrew University, Jerusalem 91904, $\;$  ISRAEL}
\email{boris@math.huji.ac.il}

\thanks{The authors were supported in part by
the Edmund Landau Center for Research in Mathematical Analysis and
Related Areas, sponsored by the Minerva Foundation (Germany).}

\date{(July 2-???), 2004}






\begin{abstract}
The monograph contains a systematic treatment of a circle of
problems in analysis and integral geometry related to inversion of
the
 Radon transform on the space of real rectangular matrices. This transform  assigns to a
function $f$  on the matrix  space
 the  integrals of $f$  over the so-called matrix planes, the linear manifolds
  determined by the corresponding matrix equations.
Different inversion methods are discussed. They rely on close
connection between  the Radon transform, the Fourier transform,
the G{\aa}rding-Gindikin fractional integrals, and matrix
modifications of the Riesz potentials. A special emphasis is made
on new higher rank phenomena, in particular,  on possibly minimal
conditions under which the Radon transform is well defined and can
be explicitly inverted. Apart of the space of Schwartz functions,
we also employ $L^p $- spaces and the space of continuous
functions. Many classical results for the Radon transform on $\rn$
are generalized to the higher rank case.

\end{abstract}

\maketitle

\setcounter{page}{4}
\tableofcontents

\mainmatter
\backmatter

\chapter*{ Introduction}

\setcounter{equation}{0}

\section{Preface}

 One of the basic problems of integral geometry reads as follows
\cite{Gel}. Let $\X$ be some space of points $x$, and let $\frT $
be a certain family of manifolds
 $\t \subset \X$.
Consider the mapping $ f(x) \to  \hat f (\t)=\int_\t f$ which
assigns to $f$ its integrals over all $\t \in \frT$. This mapping
is usually called the Radon transform  of $f$. How do we recover
$f(x) $ from $ \hat f (\t)$? This rather general problem goes back
to Lorentz, Minkowski, Blaschke, Funk, Radon etc., and has
 numerous applications; see \cite{Eh}, \cite{GGiV}, \cite{GGV},
 \cite{H1}, \cite{Na}, \cite{Ru8},  \cite{RK}, \cite{SSW}, and references therein.

A typical example is the so-called  $k$-plane transform on
$\bbr^n$ which assigns to $f(x), \; x \in \bbr^n$, a collection of
integrals of $f$ over all $k$-dimensional planes  in $\bbr^n$.
This well known transformation gets  a new flavor if we formally
replace the dimension $n$ by the product $nm$, $n \ge m$,  and
regard $\bbr^{nm}$ as the space
  $\Ma$  of $n \times m$
 matrices  $x=(x_{i,j})$. Once we have accepted this point of
view, then $f(x)$ becomes a function of matrix argument and
 new ``higher rank'' phenomena  come into play.

 Let  \index{a}{Stiefel manifold, $\vnm$} $\vnk$ be the
Stiefel manifold of orthonormal $(n-k)$-frames in $\bbr^n$, i.e.,
$ \vnk =\{\xi: \, \xi \in \frM_{n, n-k}, \,  \xi' \xi=I_{n-k} \}$
 where $\xi'$ denotes
the transpose of $\xi$,  $I_{n-k}$ is the identity
 matrix, and $\xi' \xi$ is understood in the sense of matrix multiplication.
Given $\xi\in\vnk$ and $t\in\Mt$,  we define {\it the matrix
$k$-plane} \index{a}{Matrix $k$-plane, $\tau(\xi,t)$}
\index{b}{Greek!pa@$\tau(\xi,t)$}
 \be\label{pl1}\tau\equiv
\tau(\xi,t)=\{x: \, x\in\Ma, \, \, \eq\}.\ee  Let $\frT$ be the
manifold of all matrix $k$-planes in $\Ma$. Each $\t \in \frT$ is
topologically an ordinary $km$-dimensional plane in $\bbr^{nm}$,
but the set $\frT$ has measure zero in the manifold of all such
planes. The relevant  Radon transform has the form \index{a}{Radon
transform, $\hat f (\tau)$}
 \be\label{0.1} \hat f (\tau)=\int_{x \in \tau} f(x), \qquad \tau \in
 \frT,\ee
 the integration being performed against the corresponding canonical
 measure. If $m=1$, then $\frT$ is the
  manifold of $k$-dimensional affine planes  in $\bbr^n$,
and $ \hat f (\tau)$ is the usual  $k$-plane
 transform.
 The inversion problem  for (\ref{0.1}) reads as
follows: How do we reconstruct a function $f$ from the integrals
$\hat f (\tau)$,  $\t \in \frT$? Our aim is to obtain explicit
inversion formulas for the  Radon transform (\ref{0.1}) for
continuous and, more generally, locally integrable functions  $f$
subject to possibly minimal assumptions at infinity.

 A  systematic study of Radon transforms
on matrix spaces  was initiated
 by E.E. Petrov \cite{P1} in 1967 and continued in \cite{C}, \cite{Gr1},
 \cite{P2}--\cite{P4}, \cite{Sh1},
\cite{Sh2}.   Petrov \cite{P1}, \cite{P2} considered the  case
$k=n-m$ and obtained inversion formulas for functions belonging to
the Schwartz space. His method is based on the classical idea of
decomposition of the delta function in  plane waves; cf.
\cite{GSh1}. Results of  Petrov were extended in part by L.P.
Shibasov \cite{Sh2} to all $1 \le k \le n-m$.

The following four inversion methods are customarily used in the
theory of Radon transforms. These are (a) {\it the Fourier
transform method} (based on the projection-slice theorem), (b)
{\it the method of mean value operators}, (c) {\it the method of
Riesz potentials}, and (d) {\it decomposition
 in  plane waves}. Of course, this classification is vague and other  names are also
 attributed
  in the literature.
  One could also mention implementation of  wavelet transforms,
  but this can be viewed as a certain regularization or computational  procedure for
  divergent integrals arising in (b)--(d).
In Section \ref{s0.1}, we briefly recall some known facts for the
$k$-plane transform on $\bbr^n$ which corresponds to  the case
$m=1$ in (\ref{0.1}). Our task is to  extend this theory to the
higher rank case $m>1$. The main results  are
 exhibited in Section \ref{s0.2}. A considerable part of the
 monograph (Chapters 1--3) is devoted to developing a necessary
 background and analytic tools which one usually gets almost ``for free''
 in the rank-one case.
We are not concerned with such important topics as the range
characterization, the relevant Paley-Wiener theorems \cite{P3},
\cite{P4}, the support theorems,  the convolution-backprojection
method, and numerical implementation. One should also mention a
series of papers devoted to the  Radon transform on Grassmann
manifolds; see, e.g., \cite{Go}, \cite{GK1}, \cite{GK2},
\cite{Gr2}, \cite{GR}, \cite{Ka}, \cite{Ru9}, \cite{Str2}, and
references therein. This topic is conceptually  close to ours.

\section{The $k$-plane Radon transform on $\bbr^n$}\label{s0.1}

Many authors contributed to this theory; see, e.g.,  \cite{GGiV},
\cite{H1}, \cite{Ke}, \cite{Ru7} where one can find further
references. We recall some main results which  will be  extended
in the sequel to the higher rank case.

\subsection{Definitions and the  framework of the inversion problem}
Let $\cgnk$ be the
 manifold of affine $k$-dimensional planes
in $\bbr^n$,  $0< k < n$. The integral (\ref{0.1}) with $m=1$  can
be interpreted  in different ways. Namely, if we define a
$k$-plane $\tau\in\G_{n,k}$ by (\ref{pl1}) where $\xi\in\vnk$ and
$t\in\bbr^{n-k}$ is a column vector, then {\it the $k$-plane Radon
transform} is represented as \be\label{0.2} \ft\equiv\hat f(\xi,t)
=\intl_{\{y: \, y \in \bbr^n; \, \xi' y=0\} } f(y+\xi t)\,  d_\xi
y ,\ee
 where $ d_\xi y$ is the induced Lebesgue measure on the plane
$\{y: \, \xi' y=0\}$. For $ k=n-1$, this is the classical
hyperplane Radon transform \cite{H1}, \cite{GGV}. Clearly, the
parameterization $\t\equiv\t(\xi,t)$ is not one-to-one. An
alternative parameterization is as follows. Let $\eta$ belong to
the ordinary Grassmann manifold $\gnk$  of $k$-dimensional linear
subspaces of $\rn$, and  $\eta^\perp$ be the orthogonal complement
of $\eta$.  The parameterization $$\tau\equiv\tau( \eta,
\lam),\qquad \eta\in \gnk, \qquad \lam\in \eta^\perp,$$ is
one-to-one and gives  \be \ft\equiv\hat f(\eta, \lam) = \intl_\eta
f(y+\lam) dy. \label{k-pl}\ee

The corresponding {\it dual Radon transform} of a function
$\vp(\t)$ on $\cgnk$ is defined as the mean value of $\vp
  (\t)$ over all $k$-dimensional planes $\t$  through $x\in\rn$:
\be\vx=\intl_{SO(n)}\vp(\g\eta_0 +x) d\g=\intl_{\gnk }\vp(\eta,
\Pr_{\eta^\perp}x) d\eta.\label{k-pld} \ee Here, $\eta_0$ is an
arbitrary fixed $k$-plane through the origin, and
$\Pr_{\eta^\perp}x$ denotes the orthogonal projection of $x$ onto
$\eta^\perp$.  A duality relation
 \be  \irn  f (x)\check \vp (x)dx= \icgr  \hat f(\t)  \vp(\t) d\t
 ,
\ee holds provided that either side of this equality is finite for
$f$ and $\vp$ replaced by $|f|$ and $|\vp|$, respectively
\cite{H1}, \cite{So}. Here and on, $d\t$ stands for the canonical
measure on $\G_{n,k}$; see \cite{Ru7}.

 The $k$-plane transform is injective for all $0<  k <n$ on ``standard'' function spaces where it is well defined
 (see Theorem \ref{l0.1} below). If $k<n-1$, then the inversion problem for $\hat f$ is
 overdetermined because
$$
\dim \G_{n,k}=(n-k)(k+1)> n=\dim \rn,
$$
   and the dual Radon transform is non-injective.  If $k=n-1$, then
 the dimension of the affine Grassmanian $\G_{n,n-1}$
is just $n$. In this case, the Radon transform and its dual are
injective simultaneously; see also \cite{Ru9} where this fact is
extended to Radon transforms on affine Grassmann manifolds.

The following theorem  specifies standard   classes of functions
for which the $k$-plane transform is well defined.
\begin{theorem}\label{l0.1} {}\hfil

\noindent {\rm (i)} \ Let $f(x)$ be a continuous function on $\rn$
satisfying $f(x)=O(|x|^{-a})$.  If $a>k$
 then the Radon transform $\ft$  is finite for all $ \t \in \cgnk$.

\noindent {\rm (ii)}\ If $(1+|x|)^\a f(x)\in L^1$  for some
$\a\geq k-n$, in particular, if  $f \in L^p, \; 1 \le p < n/k,$
then $\ft$ is finite for almost all $\t \in \cgnk$.
 \end{theorem}

 Here, (i) is obvious,  and (ii) is due to Solmon \cite{So}.
The conditions $a>k$, $\a\geq k-n$,  and $p < n/k$ cannot be
improved. For instance, if $p \ge n/k$ and $ f(x)=(2+|x|)^{-n/p}
(\log (2+|x|))^{-1} \; (\in L^p),$
 then $\ft \equiv \infty$.

  More informative
 statements can be found in \cite{Ru7}. For instance, if $|\tau|$
 denotes the euclidean distance from the $k$-plane $\tau$ to the
 origin, then

\be \label{0.10}\intl_{\G_{n,k}} \ft \frac{|\t|^{\lam
-n}}{(1+|\t|^2)^{\lam/2}} d\t=c \intl_{\rn} f(x)\frac{
|x|^{\lam-n}}{(1+|x|^2)^{(\lam-k)/2}} dx,
 \ee
 $$
c=\frac{ \Gam ((\lam-k)/2) \, \Gam (n/2)}{\Gam (\lam/2) \, \Gam
((n-k)/2)},\qquad Re\;\lam>k
 $$
(see Theorem 2.3 in \cite{Ru7}). Choosing $\lam=n$, we obtain  \be
\label{0.11}\intl_{\G_{n,k}} \frac{\ft}{(1+|\t|^2)^{n/2}} d\t=
\intl_{\rn} \frac{ f(x)}{(1+|x|^2)^{(n-k)/2}} dx.
 \ee
 It is instructive to compare (\ref{0.11}) with    \cite[Theorem 3.8]{So} where, instead of explicit
 equality (\ref{0.11}), one has an inequality (which is of independent interest!).
  Different
 modifications of (\ref{0.10}) and (\ref{0.11}), which control the
 Radon transform and its dual at the origin and near  spherical
 surfaces, can be found  in \cite{Ru7}.

It is worth noting that ``Solmon's condition'' $(1+|x|)^{k-n}
f(x)\in
 L^1(\rn)$ is not necessary for the a.e. existence of $\ft$ because,
 the Radon transform  essentially has an exterior (or  ``right-sided'') nature;
 see (\ref{rr}) below. For any rapidly decreasing function
 $f\in C^\infty(\rn\setminus \{0\})$ satisfying $f(x)\sim |x|^{-n}$ as
 $|x|\to 0$,  the Solmon condition fails, but $\rf<\infty$ a.e. (it is finite for all
 $\t\not\ni0$). However, if for $f\geq 0$,  one also assumes $\ft\in L^1_{loc}$ (in addition to
 $\ft<\infty$ a.e.), then,
 necessarily, $(1+|x|)^{k-n} f(x)\in L^1$  \cite[Theorem 3.9]{So}.

\subsection{The $k$-plane transform and the dual $k$-plane transform of radial functions}

There is an important  connection between the transformations
$f\to\hat f$ and $\vp\to\check \vp$,  and the Riemann-Liouville
(or Abel) fractional integrals \index{a}{Riemann-Liouville
fractional integrals} \index{b}{Latin and Gothic!ia@$I^\a_{\pm}f$}
\be \label{ri-li} (I^\a_+ f)(s) \! = \! \frac{1}{\Gamma(\a)} \!
\intl^s_0 \! \! \frac{f(r)}{(s \! - \! r)^{1-\a}} dr, \qquad \!
(I^\a_- f)(s) \! = \! \frac{1}{\Gamma(\a)} \! \intl^\infty_s  \!
\! \frac{f(r)}{(r \! - \! s)^{1-\a}} dr, \ee $Re \, \a >0$.
 Let $ \; \sigma_{n-1}  = 2 \pi^{n/2}/
 \Gamma(n/2)$ \index{b}{Greek!m@$\sigma_{n-1}$}  be the area of the unit sphere $S^{n-1}$ in
 $\bbr^n$. The following statement is known; see, e.g.,
 \cite{Ru7}.
\begin{theorem}\label{t0.2} For $x \in \rn$ and $\t \equiv (\eta,\lam) \in \cgnk$,
let \be r=|x|=\dist (o,x), \qquad s=|\lam|=\dist (o, \t)\ee denote
the corresponding euclidean distances from the origin. If $f(x)$
and $\vp(\t)$ are radial functions, i.e., $f(x)\equiv f_0 (r) $
and $\vp(\t)\equiv  \vp_0 (s)$, then   $\ft$ and $\vx$ are
represented by the Abel type integrals \be\label{rr} \ft=
\sig_{k-1}\intl_s^\infty f_0(r) (r^2 -s^2)^{k/2 -1} r dr, \ee
\be\label{drr} \vx= \frac{\sig_{k-1} \, \sig_{n-k-1}}{\sig_{n-1}
\, r^{n-2}}\intl_0^r \vp_0 (s) (r^2 -s^2)^{k/2 -1} s^{n-k-1} ds,
\ee provided that these integrals exist in the Lebesgue sense.
\end{theorem}

Obvious transformations reduce  (\ref{rr}) and (\ref{drr})  to the
corresponding fractional integrals (\ref{ri-li}).

\subsection{Inversion of the $k$-plane transform}

Abel type representations (\ref{rr}) and (\ref{drr}) enable us to
invert the $k$-plane transform and its dual for radial functions.
To this end, we utilize diverse  methods   of fractional calculus
described, e.g., in  \cite{Ru1} and \cite{SKM}. In the general
case,
 the following approaches can be applied.

{\it The method of mean value operators}. The idea of the method
amounts to original papers by
 Funk and Radon and employs the fact that  the
$k$-plane transform  commutes with  isometries of $\rn$.  Thanks
to this property,  the inversion problem reduces to the case of
radial functions if one  applies to $\ft$ the so-called {\it
shifted dual $k$-plane transform} (this terminology is due to F.
Rouvi\`{e}re \cite{Rou}). The latter  is defined on functions $\vp
(\t), \; \t \in \cgnk$, by the formula \be \check \vp_r (x) =
 \intl_{SO(n)}  \vp (\gam \t_r  +  x) \,
 d\gam, \label{drsh} \ee
where $\t_r$ is  an arbitrary fixed $k$-plane at distance $r$ from
the  origin. The integral (\ref{drsh}) is the mean value of $\vp
  (\t)$ over all $\t$ at distance $r$ from $x$. For $r=0$,  it coincides with the
dual  $k$-plane transform $\vx$. Modifications of (\ref{drsh}) for
Radon transforms on the $n$-dimensional unit sphere and the real
hyperbolic space  were introduced by Helgason; see \cite{H1} and
references therein.

The core of the method is the following.

\begin{lemma} \label{l0.3}  Let $\vp=\hat f$, $ \; f \in L^p,
\;  \, 1 \le p <n/k$, and let \index{a}{Spherical mean}
\index{b}{Latin and Gothic!mr@$\M_r f$} \be (\M_r
f)(x)=\frac{1}{\sig_{n-1}}\intl_{S^{n-1}} f(x+r\theta ) \, d
\theta, \qquad r>0, \ee be  the spherical mean of $f$. If $g_x
(s)=(\M_{\sqrt{s}} f)(x)$ and $ \psi_x(s)=\pi^{-k/2} \check
\vp_{\sqrt{s}}(x)$, then \be\label{r-rl} (I_{-}^{k/2}
g_x)(s)=\psi_x(s). \ee
\end{lemma}

If $k$ is even, (\ref{r-rl}) yields \be\label{inv4}
f(x)=\pi^{-k/2} \Big (-\frac{1}{2r}\frac{d}{dr} \Big )^{k/2}\check
\vp_r (x) \Big |_{r=0}, \ee where, for non-smooth $f$, the
equality is understood in the almost everywhere sense. This
inversion formula is of local nature because it assumes that
$\check \vp_r (x)$ is known only for arbitrary small $r$. The case
of odd $k$ is more delicate. If $f$ is continuous and decays
sufficiently fast at infinity, then (\ref{r-rl}) is inverted by
the formula
 \be f(x)=\Big (-\frac{d}{ds} \Big )^m
(I_{-}^{m-k/2} \psi_x)(s) \Big |_{s=0}, \qquad \forall m \in \bbn,
\quad m>k/2, \ee which is essentially non-local.  For $f$
belonging to $L^p$, the integral $(I_{-}^{m-k/2} \psi_x)(s)$ can
be divergent if $n/2m \le p <n/k$. It means that  the inversion
method should not increase the order of the fractional integral.
This obstacle can be circumvented by making use of  Marchaud's
fractional derivatives; see \cite[Section 5]{Ru7} for details.

 {\it The method  of   Riesz potentials and the method of plane waves}.
  It is instructive to describe these methods simultaneously
  because both employ  the
Riesz potentials  \index{a}{Riesz potential!on $\rn$, $I^\a
f$}\index{b}{Latin and Gothic!ia@$I^\a f$}  \be\label{Rp-rn} (I^\a
f)(x)=\frac{1}{\gamma_n(\a)} \intl_{\bbr^n}
 f(x-y) |y|^{\a -n} \, dy,\qquad x\in\rn,
\ee  \be \index{b}{Greek!c@$\gamma_{n} (\a)$}
  \gamma_{n} (\a)=
  \frac{2^\a\pi^{n/2}\Gamma(\a/2)}{\Gamma((n-\a)/2)},
\qquad Re \, \a >0,  \qquad \a \neq n,n+2, \ldots\; .\ee The
Fourier transform of the corresponding Riesz distribution
$|x|^{\a-n}/\gamma_{n}( \a)$ is $|y|^{-\a}$ in a suitable sense.
It means   that $I^\a$ can be regarded as $(-\a/2)$th power of
$-\Del$, where $\Del=\d^2/\d x_1^2+\dots + \d^2/\d x_n^2$ is  the
Laplace operator . For $f\in L^p(\rn)$, the integral (\ref{Rp-rn})
absolutely converges if and only if $1\leq p<n/ Re\,\a$. For $Re
\, \a \le 0$ and sufficiently smooth $f$,  the Riesz potential
$I^\a f$ is defined by analytic continuation  so that $I^0 f=f$.
Numerous inversion formulas for $I^\a f$ are known in diverse
function spaces, see \cite{Ru1}, \cite{Ru7}, \cite{Sa},
\cite{SKM}, and references therein.

 The method of plane waves is based on
decomposition of the Riesz distribution in plane waves; see
\cite{GSh1} for the   case $k=n-1$. Technically, this method
realizes as follows. In order to recover
 $f$ at a point $x$ from the data
$\ft\equiv\rf$, one  first applies  the Riesz potential operator
of the negative order $-k$ (this is a differential or
pseudo-differential operator) in the $t$-variable, and then
integrates over all $k$-planes through $x$, i.e., applies the dual
Radon transform. This method requires $f$ to be smooth enough.
 Another method, referred to as the method of  Riesz potentials, is
applicable to ``rough" functions as well. Now we utilize the same
operators but
 in the reverse order: first we apply the dual  Radon transform
to $\hat f(\t)$, and then the Riesz potential operator of  order $-k$
 in the $x$-variable. This method allows us to reconstruct any
function $f \in L^p, \; 1 \le p <n/k$, for almost all $x\in\bbr^n$
\cite{Ru7}.

 Both
methods can be realized in the framework of analytic families of
intertwining fractional integrals
\be\label{pa} (P^\a f)(\t)=\frac{1}{\gamma_{n-k}(\a)}\irn f(x)
|x-\t|^{\a+k-n} \, dx, \ee
\be\label{pda} (\stackrel{*}{P}{}^\a
\vp)(x)=\frac{1}{\gamma_{n-k}(\a)}\icgr \vp (\t) |x-\t|^{\a+k-n}
\, d\t,\ee
\[  Re\, \a
> 0, \qquad \a+k-n \neq 0,2,4, \ldots \, ,\]
where $|x-\t|$ denotes the euclidian distance between the point
$x$ and the $k$-plane $\tau$.   These operators were introduced by
Semyanistyi  \footnote{Vladimir Il'ich Semyanistyi (1925--1984), a
talented mathematician whose papers have played an important role
in integral geometry and related areas of analysis.} \cite{Se} for
$k=n-1$. They were   generalized in \cite{Ru11} and subsequent
publications of B. Rubin to all  $0< k <n$ and totally geodesic
Radon transform on spaces of constant curvature; see also
\cite{Ru7}. Semyanistyi's fractional integrals have a deep
philosophical meaning.  One can write \be\label{p-rp} P^\a
f=\tilde I^\a \hat f, \qquad \stackrel{*}{P}{}^\a \vp=(\tilde I^\a
\vp)^\vee, \ee where for $\vp(\t)\equiv\vp(\xi, t), \; \tilde
I^\a\vp$ denotes the Riesz potential on $\bbr^{n-k}$ in the
$t$-variable. If $f$ and $\vp$ are smooth enough, then  the
equalities (\ref{p-rp}) extend definitions (\ref{pa}) and
(\ref{pda}) to all complex $\a \neq n-k, \, n-k+2, \ldots \;$. In
particular,  for $\a=0$, we get the $k$-plane transform and the
dual $k$-plane transform, respectively.

The following statement generalizes the well known formula of
Fuglede \be\label{fugl}(\hat f)^\vee=c_{k,n}I^{k}f, \qquad
c_{k,n}= (2\pi)^k \sig_{n-k-1}/\sig_{n-1};
 \ee
 see \cite{Fu}, \cite[Theorem 3.6] {So}, \cite[p. 29]{H1}.

\begin{theorem}\label{t0.4} Let $f \in L^p, \; \, 1 \le p <
n/(Re\,\a+k), \;  \, Re\,\a \ge 0$. Then
 \be\label{gfu1} \stackrel{*}{P}{}^\a \hat f= c_{k,n}I^{\a+k}f.
\ee
\end{theorem}
The equality (\ref{gfu1})  gives a family of inversion formulas
(at least on a formal level): \be\label{inv1} c_{k,n}f=I^{-\a-k}
\stackrel{*}{P}{}^\a \hat f. \ee In the case $\a=0$, the last
equality reads \be \label{inv2} c_{k,n}f=\bbd^k \check \vp, \qquad
\vp =\hat f.\ee Here, $\bbd^k=I^{-k}=(-\Del)^{k/2}$ denotes the
Riesz fractional derivative that can be realized  in different
ways depending on a class of functions $f$ and the evenness of
$k$. For example, if $f$ is good enough, the $k$-plane transform
$\hat f =\vp$ can be inverted by repeated application of the
operator $-\Del$. Namely, for $k$ even, \be\label{inv3}
c_{k,n}f(x)= (-\Del)^{k/2}\check \vp (x) . \ee If $k$ is odd and
$1 \le k \le n-2$, then \be\label{inv4}
c_{k,n}f(x)=(-\Del)^{(k+1)/2}(I^1 \check \vp ) (x).\ee

Inversion formulas (\ref{inv2})--(\ref{inv4}) are typical for the
method of Riesz potentials. Furthermore, if  $\a=-k$, then
(\ref{inv1}) gives \be c_{k,n}f = \stackrel{*}{P}{}^{-k}\vp ,
\qquad \vp(\xi, t)  = \hat f(\xi, t).\ee  The last formula can be
regarded as a decomposition of $f$  in plane waves. Different
realizations of  operators $\bbd^k$ and $\stackrel{*}{P}{}^{-k}$
can be found in \cite{Ru1} and \cite{Ru7}; see also \cite{Ru12}.

\section{Main results}\label{s0.2}

This section is organized so that the reader could compare our new
results with those in the rank-one case. That was the reason why
subsections below are entitled as the similar ones in Section
\ref{s0.1}. The general organization of the material is clear from
the Contents.

\subsection{Definitions and the framework of the inversion problem}

In order  to give  precise meaning to the integral (\ref{0.1}), we
use the parameterization $\tau= \tau(\xi,t)$, $\xi\in\vnk$,
$t\in\Mt$, of the matrix $k$-plane (\ref{pl1}) and arrive at the
formula \be\label{28}\rf=\intl_{\{y: \, y\in\Ma \, ; \, \xi
'y=0\}} f(y+\xi t) \, d_\xi y,\ee where $d_\xi y$ is the induced
measure on the subspace $\{y: \, y\in\Ma \; ; \; \xi 'y=0\}$.
 This
 agrees with the formula
(\ref{0.2}) for the $k$-plane transform and can be also written in
a different way, see (\ref{4.9}).     Another parameterization
$$ \tau =\tau (\eta,\lam), \qquad \eta \in G_{n,k}, \quad
\lam=[\lam_1 \dots \lam_m]\in\Ma, \quad \lam_i \in \eta^\perp, $$
which is one-to-one (unlike that in (\ref{28})), gives
\be\index{a}{Radon transform, $\hat f (\tau)$} \hat
f(\tau)=\intl_\eta dy_1\dots \intl_\eta f([y_1+\lam_1 \dots
y_m+\lam_m ]) \, dy_m, \ee cf. (\ref{k-pl}).

We first specify the set of all triples $(n,m,k)$ for which the
Radon transform is injective. It is natural to assume that the
transformed function should depend on at least as many variables
as the original one, i.e.,
$$\dim\Gr\ge \dim\Ma .$$ By taking into account
  that $\dim\vnm = m(2n-m-1)/2$ \cite[ p. 67]{Mu}, we have
\bea\nonumber
 \dim\Gr&=&\dim(\cd)/ O(n-k)\\\nonumber
 &=& (n-k)(n+k-1)/2+m(n-k)- (n-k)(n-k-1)/2\\
 &=&(n-k)(k+m) \nonumber.\eea
The inequality $(n-k)(k+m)\geq nm$ implies    the natural
framework of the inversion problem, namely,  \be\label{1.2} 1\le
k\le n-m, \qquad m\geq 1.\ee We show (Theorem \ref{inj}) that  the
Radon transform $f\rightarrow\hat f$ is injective on the Schwartz
space $\s(\Ma)$ if and only if  the triple $(n,m,k)$ satisfies
(\ref{1.2}). The proof of this statement relies  on the matrix
generalization of {\it the projection-slice theorem}
\index{a}{Projection-Slice Theorem} (Theorem \ref{CST}) which
establishes connection between the matrix Fourier transform  and
the Radon transform. For the case $m=1$, we
 refer to \cite[p. 11]{Na} ($k=n-1$) and \cite[p. 283]{Ke}
 (any $0<k<n$); see also \cite{Sh1}, \cite{Sh2}, and \cite{OR} for the
 higher rank case.

We also  specify  classes of  functions for which the Radon
transform is well defined. Let  $I_m$ be the identity $m \times m$
matrix and let $|a|$ denote the absolute value of the determinant
of a square  matrix $a$. The following statement agrees with
Theorem \ref{l0.1}.

\begin{theorem}\label{t0.5}{}\hfil

\noindent {\rm (i)} \  Let $f(x)$ be a continuous function on
$\Ma$, satisfying $f(x)=O(|I_m+x'x|^{-a/2})$.  If $a>k+m-1$, then
the Radon transform $\hat f(\tau)$ is finite for all $\tau\in\Gr$.

\noindent {\rm (ii)} \  If $|I_m+x'x|^{\a/2}f(x)\in L^1(\Ma)$ for
some $\a\geq k-n$, in particular, if $f \in \lp$, $1\leq p<p_0$, $
 p_0=(n+m-1)/(k+m-1)$, then $\hat f(\tau)$ is finite for almost
all $\tau\in\Gr$.
\end{theorem}

As in  Theorem \ref{l0.1}, the conditions $a>k+m-1$, $\a\geq k-n$,
and $p<p_0$ cannot be improved. For instance, if $p \ge p_0$ and
$f(x)=|2I_m+x'x|^{-(n+m-1)/2p}(\log|2I_m+x'x|)^{-1} \; (\in \lp),$
 then $\ft \equiv \infty$.

In fact, we have proved  the following generalization of
(\ref{0.10}) which implies (ii) above:

\be
 \begin{array}{ll}
{\displaystyle \frac{1}{\sigk} \intl_{\vnk}d\xi\intl_{\Mt}\rf\;
\frac{|t't|^{(\lam-n)/2}}{|I_m+t't|^{\lam/2}}\;dt} \\[30pt]
={\displaystyle c \intl_{\Ma} f(x)\;
\frac{|x'x|^{(\lam-n)/2}}{|I_m+x'x|^{(\lam-k)/2}} \;dx},\qquad
Re\, \lam>k+m-1,
\end{array}
\ee where the constant $c$ is explicitly evaluated; see
(\ref{5.9}), (\ref{5.21}).

We introduce  {\it the dual Radon transform} which assigns to a
function $\vp(\t)$ on $\Gr$  the integral over all matrix
$k$-planes through $x\in\Ma$: \be\label{0.4} \df=\intl_{\tau\ni
x}\vp(\tau). \ee
 In terms of the parameterization $\tau= \tau(\xi,t)$,
$\xi\in\vnk$, $t\in\Mt$, this integral reads as follows:
\be\label{0.9}\index{a}{Dual Radon transform, $\df$} \qquad
\df=\frac{1}{\sigk}\intl_{\vnk} \varphi(\xi,\xi'x)\;d\xi, \ee
where $\sigk$ is the volume of the Stiefel manifold $\vnk$; see
(\ref{2.16}). The integral (\ref{0.9}) is well defined for any
locally integrable function $\vp$. We do not
 study the inversion problem for the dual Radon transform here.
 However, it is natural to conjecture that the dual Radon transform is
injective if and only if $k \ge n- m$, and the relevant inversion
formula can be derived from that
for the Radon transform. The respective results were obtained in
\cite{Ru9} for Radon transforms on affine Grassmann manifolds. We
plan to study this topic in forthcoming publications.

\subsection{The Radon transform and the dual Radon transform of radial functions}
We call functions $f$ on $\Ma$ and $\vp$ on $\Gr$ {\it  radial} if
they are $O(n)$ left-invariant. Each radial function $f(x)$ on
$\Ma$ ($\fc$ on $\Gr$, resp.) has the form
  $f(x)=f_0(x'x)$ ($\fc=\vp_0(t't)$, resp.).
  We establish connection between the Radon transforms of radial
functions and the {\it G{\aa}rding-Gindikin fractional integrals}
\bea \label{g+}(I_{+}^\a f)(s) &=& \frac {1}{\gma} \intl_0^s
f(r)|s-r|^{\a-d} dr,\\\label{g-} (I_{-}^\a f)(s)&=& \frac
{1}{\gma} \intl_s^\infty f(r)|r-s|^{\a-d} dr\eea associated to the
cone $\p$ of positive definite symmetric $m \times m$ matrices.
Here, $d=(m+1)/2$  and  $\gma$ is the generalized  gamma
 function (\ref{2.4}) \index{a}{Gamma
 function, $\gm (\a)$} which is also known as the Siegel gamma
 function \cite{Si}.
 In (\ref{g+}),  $s$ belongs to $\p$, and integration  is performed over the ``interval" $(0,s)=\{r:
\; r\in\p, \; s-r\in\p \}.$ In (\ref{g-}),  $s$ belongs to the
closure  $\cpm$, and one integrates over the shifted cone
$s+\p=\{r: \; r\in\p ,\; r-s\in\p \}$. For $m=1$, the integrals
(\ref{g+}) and (\ref{g-}) coincide with the Riemann-Liouville
fractional integrals (\ref{ri-li}). Section \ref{s3.2} contains
historical notes related to integrals $I_{\pm}^\a f$ in the higher
rank case.

For sufficiently good $f$, the integrals $I_{\pm}^\a f$ converge
absolutely if $Re \, \a
> d-1$ and extend  to all $\a\in\bbc$ as entire functions of $\a$. If $\a$ is real and
 belongs to {\it the
 Wallach-like  set} \index{a}{Wallach-like set, $\W$} $$ \W
=\left\{0, \frac{1}{2}, 1, \frac{3}{2}, \ldots , \frac{m-1}{2}
\right\} \cup \left\{\a:\; Re\,\a> \frac{m-1}{2} \right \},$$
  then $I_{\pm}^\a f$ are convolutions with
positive measures; see
 (\ref{2.2.2}), (\ref{2.29}). The converse
  is also true;  cf. \cite[ p. 137]{FK}.
 Owing to this property, if $\a\in\W$, then one can define the
G{\aa}rding-Gindikin fractional integrals  for  all locally
integrable functions which in the case of $I_{-}^\a $ obey some
extra conditions at infinity. Different explicit formulas for
$I_{\pm}^\a f$, $\a\in\W$, are obtained in Section \ref{s2.2}.

 A natural analog of Theorem \ref{t0.2} for the higher rank case employs
 fractional integrals (\ref{g+})--(\ref{g-}) and reads as
 follows.

\begin{theorem}\label{t0.1}
For $x\in\Ma$ and $t\in\Mt$, let $r=x'x$, $s=t't$.

\noindent {\rm (i)} \
 If $f(x)$ is radial, i.e., $f(x)=\f0 (x'x)$, then   \be\label{0.7} \rf=\pi^{km/2}
(I_{-}^{k/2}f_0)(s). \ee

 \noindent {\rm (ii)} \ Let $\vp(\t)$ be radial, i.e., $ \; \fc \! = \!
 \vp_0(t't)$. We denote
 \[\Phi_0(s)=|s|^\del \vp_0 (s), \qquad \del=(n-k)/2-d, \quad
d=(m+1)/2.\] If $1 \le k \le n-m$, then for any matrix $x\in\Ma$ of rank $m$,
 \be\label{0.8}
\df=c\; |r|^{d-n/2}(I_{+}^{k/2} \Phi_0)(r), \quad
c=\pi^{km/2}\sigma_{n-k,m}/\sigma_{n,m},\ee $\sigma_{n,m}$ being
 the  volume of the Stiefel manifold $\vnm$.
 Equalities (\ref{0.7})
and (\ref{0.8}) hold for any locally integrable functions provided
that either side exists in the Lebesgue sense.
\end{theorem}

\subsection{Inversion of the Radon transform}
 Theorem \ref{t0.1} reduces  inversion
  of the Radon transform and the dual  Radon transform of radial functions to
the similar problem for the G{\aa}rding-Gindikin fractional
integrals. If $f$ belongs to the space $\D (\p)$ of infinitely
differentiable functions compactly supported away from the
boundary $\partial \p$, we proceed as follows.
 For $r = (r_{i,j}) \in  \p$, we consider the
  differential operators $D_\pm$ defined by
 $$ D_+ =
 \det \left ( \eta_{i,j} \, \frac {\partial}{\partial
 r_{i,j}}\right ), \quad
  \eta_{i,j}= \left\{
 \begin{array} {ll} 1  & \mbox{if $i=j$}\\
1/2 & \mbox{if $i \neq j$}
\end{array}
\right. ,\qquad D_{-}=(-1)^{m} D_{+}.$$ If $\a \in \bbc$ and $j
\in \bbn$, then $$ D_{\pm}^j I_{\pm}^\a f=I_{\pm}^\a D_{\pm}^j
f=I_{\pm}^{\a-j} f,\qquad f\in \D (\p).$$ Hence $I_{\pm}^{\a}f$
can be inverted by the formula \be\label{inv-g} f =D^j
_{\pm}I_{\pm}^{j-\a} g, \quad g=I_{\pm}^{\a}f, \ee where $j\geq\a$
if $\a= 1/2, 1,\dots\; , $ and  $j>Re\,\a+d-1$ otherwise. For
``rough'' functions $f$, we  interpret (\ref{inv-g}) in the sense
of distributions; see Lemmas \ref{l-inv1} and  \ref{ligm}. A
challenging open problem is to obtain ``pointwise'' inversion
formulas for the G{\aa}rding-Gindikin fractional integrals similar
to those for the classical Riemann-Liouville integrals \cite{Ru1},
\cite{SKM}. It is desirable  these formulas would not contain
neither operations in the sense of distributions nor the Laplace
(or the Fourier) transform, i.e., the problem would be resolved in
the same language as it has been stated.

In the general case, we develop the following inversion methods
for the Radon transform  which extend the classical ones to the
higher rank case.

{\it The method of mean value operators.} The Radon transform
$\ft$ commutes with  the group $M(n,m)$ of matrix motions
 sending $x\in\Ma$ to $\g x\b+b$,
where $\g\in O(n)$ and  $\b\in O(m)$ are orthogonal matrices, and
$b\in\Ma$. This property allows us to apply  a certain mean value
operator  to $\ft$ and thus reduce the inversion problem   to the
case of radial functions.

 Let us make the following definitions.
We supply the matrix space $\Ma$ and the manifold $\Gr$ of  matrix
$k$-planes with an entity that serves as a substitute for the
euclidean distance.   Specifically, given two points $x$ and $y$
in $\Ma$, {\it the matrix distance} $d(x,y)$ is defined by
$$d(x,y)=[(x-y)'(x-y)]^{1/2}.$$
 Accordingly, the ``distance''  between  $x\in\Ma$ and a matrix $k$-plane $\tau=\tau(\xi, t) \in\Gr$
is defined by $$d( x, \tau)=[(\xi ' x-t)'(\xi ' x-t)]^{1/2}.$$
Unlike the rank-one case, these quantities are not scalar-valued
and represented by  positive semi-definite matrices. Radial
functions $f(x)$ on $\Ma$ and $\vp(\tau)$ on $\Gr$ virtually
depend
 on the matrix distance to the
$0$-matrix from $x$ and $\tau$, respectively.

 For functions
$\vp (\t), \; \t \in \Gr$, we define  {\it the shifted dual Radon
transform}
 \be\label{drsh1} \check \vp_s (x)=\int_{d(x ,\tau)=s^{1/2}}
\vp(\tau),  \qquad s\in\p .\ee This is  the mean value of
$\vp(\t)$ over all matrix planes at distance $s^{1/2}$ from $x$.
For $s=0$, it coincides with the dual Radon transform $\vx$. A
matrix generalization of the classical spherical mean is
$$
(M_r f)(x)=\frac{1}{\sigma_{n,m}}\intl_{\vnm}
f(x+vr^{1/2})dv,\qquad r\in \p. $$ The following statement
generalizes  Lemma \ref{l0.3} to the case $m>1$.
 \begin{lemma}
  For  fixed  $x \in \Ma$, let $F_x(r)=(M_r f)(x)$. If $f\in\lp$, $ 1\leq
p<(n+m-1)/(k+m-1)$, then
  \be\label{0.6}
(\hat f)_s^{\vee}(x)=\pi^{km/2} (I_{-}^{k/2}F_x)(s), \qquad
s\in\p. \ee
\end{lemma}

Owing to (\ref{0.6}) and the inversion formula for the
G{\aa}rding-Gindikin fractional integrals, one can recover $f(x)$
from the Radon transform $\ft$ as follows.
\begin{theorem}
Let  $ f\in\lp$,
$$
1\le k\leq n-m, \qquad 1\leq p<\frac{n+m-1}{k+m-1}.$$
 If
$\vp(\t)=\ft$, then \be f(x)=\pi^{-km/2}\lim\limits_{r\to
0}^{(L^p)}(D_-^{k/2} \Phi_x )(r), \qquad \Phi_x (s)=\check \vp_
 s (x),\ee where  $D_-^{k/2}=I^{-k/2}_-$ is understood in the
 sense of $\D'$-distributions.
\end{theorem}

{\it The method  of   Riesz potentials and the method of plane
waves}. We extend the method of Riesz potentials to the higher
rank case and give a new account of the method of plane waves. The
second method was developed  by E.E. Petrov \cite{P1} for the case
$k=n-m, $ and outlined by L.P. Shibasov \cite{Sh2} for all $1\le k
\le n-m$. Both methods employ  the following matrix analog of the
Riesz potential (\ref{Rp-rn}): \be\label{rp0} (I^\a
f)(x)=\frac{1}{\gam_{n,m} (\a)} \intl_{\Ma} f(x-y) \, |y|_m
^{\a-n}\,  dy,\quad Re\,\a>m-1, \ee
 where $|y|_m =\det(y'y)^{1/2}$,
\be\label{gam0}\gam_{n,m} (\a)=\frac{2^{\a m} \, \pi^{nm/2}\,
\Gam_m (\a/2)}{\Gam_m ((n-\a)/2)},\qquad \a\neq n-m+1, \,  n-m+2,
\ldots . \ee \index{a}{Cayley-Laplace operator!on $\Ma$, $\Del$}
We consider {\it the Cayley-Laplace operator}
$$(\Del f)(x)=(\det(\d '\d) f)(x), \qquad x=(x_{i,j})\in \Ma, $$
where $\partial$ is the $n\times m$  matrix whose entries are
partial derivatives $\d/\d x_{i,j}$. The name ``Cayley-Laplace''
was attributed to this operator by S.P. Khekalo \cite{Kh} and
suits it admirably. On the Fourier transform side, the action of
$\Del$ represents a multiplication by $(-1)^m |y|_m^2$. Since the
Fourier transform of the Riesz distribution
$|x|_m^{\a-n}/\gam_{n,m} (\a)$ is $|y|_m^{-\a}$ in a suitable
sense, then  $I^\a$ can be regarded as $(-\a/2)$th power of the
operator $(-1)^m\Del$. For $Re\,\a\leq m-1$, $\a\neq n-m+1, \,
n-m+2, \ldots\;$ and sufficiently smooth $f$, the Riesz potential
is defined as analytic continuation of the integral (\ref{rp0}).
If $\a$ is real and belongs to the corresponding Wallach-like
 set \index{a}{Wallach-like set, $\V$}
$$
\V=\{0, 1, 2, \ldots,  k_0 \}\cup \{\a  :Re\, \a \!
> \! m \! - \! 1; \; \a\neq n-m+1, \, n-m+2, \ldots\},
$$
$$
k_0=\min(m-1, n-m),
$$
then $I^\a f$ is a convolution with a positive measure; see
 (\ref{des})--(\ref{I0}). This result allows us to extend the
 original
definition of the Riesz potential to locally integrable functions
provided  $\a\in\V$. It was shown in \cite{Ru10}, that if
\be\label{cond1} f\in \lp,  \qquad 1 \le p <\frac{n}{Re\,\a +m-1},
\qquad \a\in\V, \ee
 then $(I^\a f)(x)$ is finite for
almost all $x\in\Ma$. On the other hand,    if $$p\geq
\frac{n+m-1}{k+m-1}$$
 then there is a function $f\in\lp$ such that $(I^k
f)(x)\equiv\infty$ (see (\ref{fu}), Theorem \ref{t5.1}, and
Appendix B). It is still not clear what happens if
\be\label{cond2}\frac{n}{Re\,\a +m-1}\leq p<
\frac{n+m-1}{Re\,\a+m-1}\;.\ee This gap represents an open
problem.

As in the rank-one case,   the  method  of   Riesz potentials and
the method of plane waves employ  intertwining operators
 \be\label{ppd0} P^{\a}f=\tilde I^\a\hat f, \qquad
\pd\vp=(\tilde I^\a\vp)^\vee, \ee  $$\a\in\bbc,\qquad \a  \neq n-k
- m + 1, \, n-k  - m + 2,\dots \;.$$ Here, $1\le k \le n-m$ and
$\tilde I^\a$ is the Riesz potential on $\Mt$  acting in the
$t$-variable. If $$Re\,\a
>m-1,\qquad \a\neq n-k-m+1,\; n-k-m+2, \ldots\; ,$$ then  \bea\label{P0}
\index{a}{Generalized Semyanistyi fractional integrals,  $P^{\a}$,
$\pd$} \qquad \qquad (P^\a f)(\xi, t)\label{p0}&=&
\frac{1}{\g_{n-k,m}(\a)}\intl_{\Ma} f(x)\,|\xi ' x-t|_m^{\a+k-n}\,dx,\\
(\pd\vp)(x)\!\!&=&\!\!\frac{1}{
\g_{n-k,m}(\a)}\intl_{\vnk}\!\!d_\ast
\xi\!\!\intl_{\Mt}\!\!\!\!\!\!\vp(\xi,t)\,|\xi' x-t|_m^{\a+k-n}\,
dt\label{pd0}, \eea where $d_\ast \xi=\sigk^{-1} \, d\xi$ is the
normalized measure on $\vnk$ and $\g_{n-k,m}(\a)$ is the
normalizing constant in the definition of the Riesz potential on
$\Mt$; cf. (\ref{gam0}). \index{b}{Latin and Gothic!dst@$d_\ast
\xi$} We call $P^\a f$ and $\pd\vp$ {\it the generalized
Semyanistyi fractional integrals}. For $\a=0$,  operators
(\ref{ppd0}) include the Radon transform and its dual,
respectively.
 In the case  $m=1$, the integrals (\ref{P0})  and (\ref{pd0})  coincide with those in (\ref{pa}) and
 (\ref{pda}).

 \begin{theorem}\label{t0.9}
Let $1\le k\le n-m$, $ \; \a\in\V$. Suppose that
 $$
f\in L^p, \qquad 1\leq p<\frac{n}{Re\,\a+k+m-1}\;.$$
 Then \be \label{0.5}(\pd\hat f)(x)=
c_{n,k,m} (I^{\a+k} f)(x),\ee
 $$
 c_{n,k,m}=2^{km}\pi^{km/2}\gm\left(\frac{n}{2}\right)/\gm\left(\frac{n-k}{2}\right).
 $$
 In particular,
 \be\label{fu0}\index{a}{Fuglede
formula} (\hat f)^{\vee} (x) \! = \! c_{n,k,m} (I^k f)(x) \ee (the
generalized Fuglede formula).
\end{theorem}

Theorem \ref{t0.9} is a higher rank copy of Theorem \ref{t0.4}.
Thus, as in the case $m=1$,
   inversion of  the Radon transform reduces  to  inversion  of the
Riesz potential. Unlike the   rank-one case,  we have not
succeeded
 to obtain a pointwise inversion formula for the higher rank Riesz
potential so far,  and use the theory of distributions. The
properly defined space $\Phi(\Ma)$
 of test functions  is invariant under the action of
  Riesz potentials, and
\be\label{rd0}(I^{-\a}\phi)(x)=(\F^{-1} \, |y|_m^{\a} \,
\F\phi)(x),\qquad \phi\in\Phi,\qquad \a\in\bbc, \ee in terms of
the Fourier transform.
\begin{theorem}\label{t8.10} Let $f \in L^p (\Ma), \; 1\leq
p<n/(k+m-1)$, so that  the Radon transform $\vp(\tau)=\hat
f(\tau)$ is finite for almost all matrix $k$-planes $\tau$. The
function  $f$ can be recovered from $\vp$ in the sense of
$\Phi'$-distributions by the formula \be \index{a}{Inversion
formula!for the Radon transform} c_{n,k,m}(f,\phi)=(\check \vp,
I^{-k}\phi), \qquad \phi \in \Phi, \ee where
 $I^{-k}$ is the operator (\ref{rd0}). In particular, for  $k$ even,
\be \index{a}{Inversion formula!for the Radon
transform} c_{n,k,m}(f,\phi)=(-1)^{mk/2}(\check \vp,
\Del^{k/2}\phi), \qquad \phi \in \Phi. \ee
\end{theorem}

This theorem reflects the essence of the method of
Riesz potentials.
 The case $\a=-k$ in
(\ref{0.5}) provides the following inversion result in the framework of
 the method of plane waves.
\begin{theorem}\label{t6.13}
Let $1\le k \le n-m$. If  $f$ belongs to the Schwartz space
$\s(\Ma)$, then
 the Radon
transform $\vp=\hat f$ can be inverted by the formula
\be\label{gfu2} c_{n,k,m}
f(x)=(\stackrel{*}{P}\!{}^{-k}\vp)(x).\ee
\end{theorem}
Explicit expressions for $\stackrel{*}{P}\!{}^{-k}\vp$ are given by
(\ref{pw1})-(\ref{pw3}).

\section{Acknowledgements}
We are indebted to our numerous friends and colleagues who sent us their papers, discussed different related topics, and helped  to achieve better understanding of the matter.  Their remarks and encouragement are
 invaluable. Among these people  are  professors M. Agranovsky,
 J.J. Duistermaat,
J. Faraut, H. Furstenberg,
 S.  Gindikin,  F. Gonzalez,  S. Helgason,  S.P. Khekalo,  A. Mudrov,
G. \'Olafsson,  E.E. Petrov, E.T. Quinto,   F. Rouvi\`{e}re,   E.
M. Stein, L. Zalcman, and
 others. Main results of this
manuscript were delivered by the  second co-author at the Special
Session on Tomography and Integral Geometry during the
 2004 Spring Eastern Section Meeting of the American Mathematical
 Society.

\newpage

\chapter{Preliminaries}

 \setcounter{equation}{0}

In this chapter, we establish our notation and recall some basic
facts important for the sequel. The main references  are
\cite{FZ}, \cite{Herz}, \cite{Mu}, \cite{T}.

\section{Matrix spaces. Notation}\label{s2.1}

  Let $\frM_{n,m}$
  \index{b}{Latin and Gothic!mnm@$\frM_{n,m}$} be the
space of real matrices having $n$ rows and $m$
 columns; \index{b}{Latin and Gothic!mnmk@$\frM_{n,m}^{(k)}$} $\Mk\subset\Ma$ is the subset of all matrices of rank $k$.
   We identify $\frM_{n,m}$
 with the real Euclidean space $\bbr^{nm}$.
The latin letters $x$, $y$, $r$, $s$, etc. stand for both the
matrices and the points since it is always clear from the context
which is meant. We use a standard notation $O(n)$  \index{b}{Latin
and Gothic!on@$O(n)$} and $SO(n)$ \index{b}{Latin and
Gothic!son@$SO(n)$} for the group of real orthogonal $n\times n$
matrices and its connected component of the identity,
respectively. The corresponding invariant measures on $O(n)$ and
$SO(n)$ are normalized to be of total mass 1. We denote by
$M(n,m)$ \index{b}{Latin and Gothic!ma@$M(n,m)$} the group of
transformations sending $x\in\Ma$ to $\g x\b+b$, where $\g\in
O(n)$,  $\b\in O(m)$, and  $b\in\Ma$. We call $M(n,m)$ {\it the
group of matrix motions}.

 If $x=(x_{i,j})\in\Ma$, we write $dx=\prod^{n}_{i=1}\prod^{m}_{j=1}
 dx_{i,j}$ for the elementary volume in $\Ma$. In the following,
  $x'$ denotes the transpose of  $x$, $I_m$ \index{b}{Latin and Gothic!im@$I_m$}
   is the identity $m \times m$
  matrix, and $0$ stands for zero entries. Given a square matrix $a$,  we denote by
  $\det(a)$ \index{b}{Latin and Gothic!det@$\det(a)$}
  the determinant of $a$, and by $|a|$ \index{b}{Latin and Gothic!am@\texttt{"|}$a$\texttt{"|}}the absolute value of
  $\det(a)$;
  $\tr (a)$ \index{b}{Latin and Gothic!tr@$\tr (a)$} stands for the trace of $a$.

  Let $\S_m$ \index{b}{Latin and Gothic!sm@$\S_m$} be the space of $m \times m$ real symmetric matrices
$s=(s_{i,j}),
 \, s_{i,j}=s_{j,i}$. It is a measure space   isomorphic to $\bbr^{m(m+1)/2}$
 with the volume element $ds=\prod_{i \le j} ds_{i,j}$.
We denote by  $\p$  \index{b}{Latin and Gothic!pm@$\p$} the cone
of positive definite matrices in $\S_m$; $\cpm$ \index{b}{Latin
and Gothic!pmc@$\overline\P_m$} is the closure of $\p$, that is
the set of  all positive semi-definite $m\times m$ matrices. If
$m=1$, then $\p$ is the open half-line $ (0, \infty)$ and $\cpm=
[0, \infty)$. For $r\in\p$ ($r\in\cpm$) we write $r>0$ ($r\geq
0$). Given
 $s_1$ and  $s_2$ in  $S_m$, the inequality $s_1 > s_2$  means $s_1 - s_2 \in
 \p$. If $a\in\cpm$ and $b\in\p$, then  the
symbol $\int_a^b f(s) ds$ denotes
 the integral over the set
$$
 \{s : s \in \p, \, a<s<b \}=\{s : s-a \in \p, \, b-s \in \p \}.
$$
For $s \in  S_m$, we denote by $s^{\lam}_{+}$ \index{b}{Latin and
Gothic!sl@$s^{\lam}_{\pm}$} the function which equals $|s|^{\lam}$
for $s \in \p$ and zero otherwise; $s^{\lam}_{-}$ equals
$(-s)^{\lam}_{+}$.

  The group $G=GL(m,\bbr)$ \index{b}{Latin and Gothic!gl@$GL(m,\bbr)$} of
 real non-singular $m \times m$ matrices $g$ acts transitively on $\p$
  by the rule $r \to g'rg$.  The corresponding $G$-invariant
 measure on $\p$ is \index{b}{Latin and
Gothic!dr@$d_{*} r$} \be\label{2.1}
  d_{*} r = |r|^{-d} dr, \qquad |r|=\det (r), \qquad d= (m+1)/2 \ee \cite[p.
  18]{T}. The cone $\p$ is a $G$-orbit in $\S_m$ of the identity matrix $I_m$.
  The boundary $\partial\p$   \index{b}{Latin and
Gothic!dpm@$\partial\p$} of $\p$
  is a union   of $G$-orbits of $m\times m$ matrices
\[e_k=\left[\begin{array}{ll}  I_k&0\\
0&0
\end{array}\right], \qquad k=0,1, \ldots , m-1.\]
   More information about the boundary structure of $\p$ can be
found in \cite[p. 72]{FK}  and \cite[ p. 78]{Bar}.

 Let $T_m$ be the  \index{b}{Latin and
Gothic!tm@$T_m$} subgroup of  $GL(m,\bbr)$ consisting of upper
triangular matrices \be\label{2.17}
 t=\left[\begin{array}{ccccc} t_{1,1} & {}   & {}  & {}  & {} \\
                              {} & {.}  & {}  & t_{*}  & {} \\
                              {} & {}   & {.} & {}  & {} \\
                              {} & {0}   & {}  & {.} & {} \\
                               {} & {}   & {}  & {}  & t_{m,m}

\end{array} \right], \qquad t_{i,i} >0, \ee $$t_{*}=\{ t_{i,j}: i<j \} \in \bbr^{m(m-1)/2}.
$$ Each $r \in \p$ has a unique representation $r=t't, \; t \in
T_m$, so that
 \bea \label{2.2}
 \intl_{\p} f(r) dr &=& \intl_0^\infty  t_{1,1}^m \,
dt_{1,1}\intl_0^\infty  t_{2,2}^{m-1} \, dt_{2,2} \,
 \ldots \nonumber \\ &\times & \intl_0^\infty  t_{m,m}  \tilde f (t_{1,1}, \ldots ,t_{m,m}) \,
 dt_{m,m},
\eea
\[
 \tilde{ f} (t_{1,1}, \ldots ,t_{m,m})=2^m  \intl_{\bbr^{m(m-1)/2}} f(t't) \, dt_{*}, \quad dt_{*}=\prod_{i<j}
 dt_{i,j},
\]
 \cite[p. 22]{T} , \cite[p. 592]{Mu}. \index{b}{Latin and
Gothic!dts@$dt_{*}$}In the last integration, the diagonal
 entries of the
 matrix $t$ are given by the arguments of $\tilde f$, and the
 strictly upper
 triangular entries of $t$ are variables of integration.

We recall some useful formulas for Jacobians;  see, e.g., \cite[
pp. 57--59]
 {Mu}. \index{a}{Jacobians}
 \begin{lemma}\label{12.2} \hskip10truecm

\noindent
 {\rm (i)} \ If $ \; x=ayb$ where $y\in\Ma, \; a\in  GL(n,\bbr)$, and $ b \in  GL(m,\bbr)$, then
 $dx=|a|^m |b|^ndy$. \\
 {\rm (ii)} \ If $ \; r=q'sq$ where $s\in S_m$, and $q\in  GL(m,\bbr)$,
  then $dr=|q|^{m+1}ds$. \\
  {\rm (iii)} \ If $ \; r=s^{-1}$ where $s\in \p$,   then $r\in
  \p$,
  and $dr=|s|^{-m-1}ds$.

 \end{lemma}

Some more notation are in order:
 $\bbn$ denotes the set of all positive integers; $\bbr_+=(0,\infty)$; \index{b}{Latin and
Gothic!rp@$\bbr_+$} $\del_{ij}$ \index{b}{Greek!deli@$\del_{ij}$}
is the usual Kronecker delta.
 The  space $C(\Ma)$ \index{b}{Latin and
Gothic!cma@$C(\Ma)$} of continuous functions, the Lebesgue space
 $\lp$, \index{b}{Latin and
Gothic!lp@$\lp$}
 and the Schwartz space  $\s(\Ma)$ \index{b}{Latin and
Gothic!sma@$\s(\Ma)$}  of $C^\infty$ rapidly decreasing functions
are
 identified with the corresponding spaces on $\bbr^{nm}$.
We denote by $C_c(\Ma)$ \index{b}{Latin and
Gothic!cmac@$C_c(\Ma)$} the space of compactly supported
continuous functions on $\Ma$.
  The  spaces $C^\infty(\p)$, \index{b}{Latin and
Gothic!ci@$C^\infty(\p)$}
  $C^k(\p)$ \index{b}{Latin and
Gothic!ck@$C^k(\p)$}
  and $\s(\p)$ \index{b}{Latin and
Gothic!spm@$\s(\p)$} consist of
  restrictions onto $\p$ of the corresponding functions on $\S_m \sim \bbr^{N}, \;
  N=m(m+1)/2$. We denote by $\D(\p)$ \index{b}{Latin and
Gothic!dpm@$\D(\p)$} the space of functions  $f \in C^\infty(\p)$
with $\supp f \subset
  \p$; $ L^1_{loc}(\p)$ \index{b}{Latin and
Gothic!lo@$L^1_{loc}(\p)$} is the space of locally integrable
functions on $\p$.   The Fourier transform \index{b}{Latin and
Gothic!ff@$\F f$}\index{a}{Fourier transform, $\F f$} of a
function $f\in L^1(\Ma)$ is defined by \be\label{ft} (\F
f)(y)=\intl_{\Ma} \exp(\tr(iy'x)) f (x) dx,\qquad y\in\Ma \; .\ee
This is the usual Fourier transform on $\bbr^{nm}$ so that the
relevant Parseval formula \index{a}{Parseval formula} reads
\be\label{pars} (\F f, \F \vp)=(2\pi)^{nm} \, (f,\vp),\ee where
$$(f, \vp)=\intl_{\Ma} f(x) \vp(x) \, dx.$$
 We write $c$, $c_1$, $c_2, \dots$ for different constants
the meaning of which is clear from the context.

\section[Gamma and beta functions ]{Gamma and beta functions of the cone of positive definite matrices}\label{s2.2}
  The generalized  {\it gamma
 function} \index{a}{Gamma
 function, $\gm (\a)$} associated to the cone $\p$ of positive definite symmetric matrices
 is defined by \index{b}{Greek!c@$\gm (\a)$}
\be\label{2.4}
 \gm (\a)=\intl_{\p} \exp(-\tr (r)) |r|^{\a -d} dr, \qquad d=(m+1)/2. \ee
This is also known as the {\it Siegel integral}  \cite{Si} (1935)
arising in number theory. Actually it came up earlier in
statistics. Substantial generalizations of (\ref{2.4}) and
  historical comments can be found in \cite{Gi1}, \cite[ Chapter VII]{FK},
 \cite[ p. 41]{T}, \cite[p. 800]{GrRi}.

Using (\ref{2.2}), it is easy to check \cite[ p. 62]{Mu} that the
integral (\ref{2.4}) converges absolutely
 if and only if $Re \, \a>d-1$.  Moreover, $\gm (\a)$ is a product of
 ordinary gamma functions:
\be\label{2.5}
 \gm (\a)=\pi^{m(m-1)/4}\prod\limits_{j=0}^{m-1} \Gam (\a- j/2)\;
 .
 \ee This implies a number of
 useful formulas:
\be
\frac{\gm(d-\a+1)}{\gm(d-\a)}=(-1)^{m}\frac{\gm(\a)}{\gm(\a-1)},
\qquad d=(m+1)/2, \ee
 \be\label{Poh} (-1)^m
\frac{\gm(1-\a/2)}{\gm(-\a/2)}=2^{-m}\frac{\Gam(\a+m)}{\Gam(\a)}=2^{-m}(\a,m),
\ee where $(\a,m)=\a(\a+1)\cdots (\a+m-1)$ is the  Pochhammer
symbol. \index{a}{Pochhammer  symbol, $(\a,m)$ } Furthermore, for
$1 \le k<m, \; k \in \bbn, $ the equality (\ref{2.5}) yields
\be\label{2.5.2}
\frac{\gm(\a)}{\gm(\a+k/2)}=\frac{\gk(\a+(k-m)/2)}{\gk(\a+k/2)},\ee
\be\label{2.5.1} \gma=\pi^{k(m-k)/2}
\Gamma_{k}(\a)\Gamma_{m-k}(\a-k/2).\ee In particular, for $k=m-1$,
 \be \gma=\pi^{(m-1)/2} \Gamma_{m-1}(\a)\Gamma(\a-(m-1)/2).\ee

The {\it beta function} \index{b}{Greek!b@$B_m (\a
,\b)$}\index{a}{Beta
 function, $B_m (\a ,\b)$ }of the cone  $\p$ is defined by
\be\label{2.6}
 B_m (\a ,\b)=\intl_0^{I_m} |r|^{\a -d} |I_m-r|^{\beta -d} dr, \ee
 where, as above, $d=(m+1)/2$.
This integral converges absolutely if and only if $Re
 \, \a, Re \, \b >d-1$. The following classical relation holds
\be\label{2.6.1}
 B_m (\a ,\b)=\frac{\gm (\a) \gm (\b)}{\gm (\a+\b)}
\ee \cite[ p. 130]{FK}. In Appendix A, we present a table of some
integrals that can be easily  evaluated in terms of the
generalized gamma and beta functions.  Most of these integrals are
known but scattered all over the different sources. This table
will be repeatedly referred to in the sequel. For the sake of
completeness, we mention the paper \cite{Ner} containing many
other formulas of this type in the very general set-up.

\section{The Laplace transform}\label{s2.4}
 Let $\T_m=\p +i\S_m$ \index{b}{Greek!pa@$\T_m$}be the generalized half-plane in the space $\S_m^{\bbc}=\S_m
 +i\S_m$, \index{b}{Latin and
Gothic!smc@$\S_m^{\bbc}$}the complexification of $\S_m$. The
domain $\T_m$
 consists  of all complex symmetric
matrices $z=\sigma+i\om$ such that $\sigma= Re\;z \in \p$, and
$\om=Im \;z\in \S_m$.  Let $f$ be a locally integrable function on
$\S_m$,
 $f(r)=0$ if $r\notin \cpm$, and $\exp(-\tr(\sigma_0
r))f(r)\in L^1(\S_m)$ for some $\sigma_0\in\p$. Then the integral
\be\label{3.11}\index{b}{Latin and Gothic!lf@$Lf$}
 (Lf)(z)=\intl_{\p} \exp(-\tr(zr))f(r)dr
\ee is absolutely convergent in the (generalized) half-plane
$Re\;z>\sigma_0$ and represents a complex analytic function there.
$Lf$ is called the Laplace transform  of $f$.  If \be\label{fts}
(\F g)(\om)=\intl_{\S_m} \exp(\tr(i\om s)) g (s) ds,\qquad
\om\in\S_m, \ee
 is the Fourier transform \index{a}{Fourier transform, $\F
f$} of  a function $g$  on $\S_m$, then
$$
(Lf)(z)=(\F g_\sigma)(-\om),
$$
where $g_\sigma(r)=\exp(-\tr(\sigma  r))f(r)\in L^1(\S_m)$ for
$\sigma >\sigma_0$. Thus, all properties of the Laplace transform
are obtained from the general Fourier transform theory for
Euclidean spaces  \cite{Herz}, \cite[ p. 126]{Vl}. In particular,
we have the following statement.
\begin{theorem}\label{tCauchy}
Let $\exp(-\tr(\sigma_0 r))f(r)\in L^1(\P_m)$, $ \sigma_0 \in \p$,
and let
$$
\intl_{\S_m} |(Lf)(\sig+i\om)| d\om <\infty
$$
for some  $\sig>\sig_0$. Then the Cauchy inversion formula
 holds:
\be\label{Cauchy}\frac{1}{(2\pi i)^N}\intl_{Re\, z=\sig}
\exp(\tr(sz)(Lf)(z) dz= \left \{
\begin{array} {ll} f(s) &
  \mbox{if $s\in\p$,}\\
{} \\
 0
  & \mbox{if  $s\notin\p$},
\end{array}
\right.  \ee $ N=m(m+1)/2$. The integration in (\ref{Cauchy}) is
performed over all $z=\sig+i\om$  with fixed $\sig>\sig_0$ and
$\om$ ranging over $\S_m$.
\end{theorem}
 The following uniqueness result for the
 Laplace transform immediately follows from injectivity of
the Fourier transform of tempered distributions.
\begin{lemma}\label{lap}
  If $f_1 (r)$ and $ f_2 (r)$ satisfy \be\label{2.35} \exp(-\tr(\sigma_0
r))f_j(r)\in L^1(\P_m), \quad j=1,2, \ee
 for some $ \sigma_0 \in \p$,  and $(Lf_1)(z)=(Lf_2)(z)$ whenever
$Re\;z>\sigma_0$, then $f_1(r)= f_2(r)$  almost everywhere on $
\S_m$.
\end{lemma}

The convolution theorem for the Laplace transform
 reads as
follows: If $f_1$ and  $f_2$ obey (\ref{2.35}),
 and
$$
(f_1 \ast f_2)(r)=\intl_0^r f_1(r-s)f_2(s)ds
$$
is the Laplace convolution, then
 \be\label{2.18}
L(f_1 \ast f_2)(z)=(Lf_1)(z)(Lf_2)(z), \quad Re\;z>\sigma_0.\ee

\begin{lemma}
Let $z \in \T_m=\p +i\S_m , \; d=(m+1)/2$. Then
 \be\label{2.19}
\intl_{\p} \exp(-\tr(zr))|r|^{\a-d} dr=\gm (\a)
\det(z)^{-\a},\qquad Re\;\a>d-1. \ee
\end{lemma}
 Here and  throughout the
paper, we set $\det(z)^{-\a}=\exp(-\a \log\det(z))$ where  the
branch of $\log\det(z)$ is chosen so that  $\det(z)=|\sig|$ for
real $z=\sigma\in\p$.

 The equality (\ref{2.19}) is well known,
see, e.g., \cite[p. 479]{Herz}. For real $z=\sigma\in\p$, it can
be easily obtained by changing variable $r \to z^{-1/2}rz^{-1/2}$.
For complex $z$, it follows by analytic continuation.

\section{Differential operators $ D_\pm$}
Let $r = (r_{i,j}) \in  \p$. We define the
 following differential operators \index{b}{Latin and Gothic!db@$ D_\pm$} acting in the $r$-variable:
 \be\label{2.50} D_+ \equiv D_{+,\, r}=
 \det \left ( \eta_{i,j} \, \frac {\partial}{\partial
 r_{i,j}}\right ), \quad
  \eta_{i,j}= \left\{
 \begin{array} {ll} 1  & \mbox{if $i=j$}\\
1/2 & \mbox{if $i \neq j,$}
\end{array}
\right. \ee \be\label{2.51}   D_{-}\equiv D_{-,\, r}=(-1)^{m}
D_{+,\, r}.\ee

\begin{lemma} If $f \in \D(\p)$, and the
 derivatives $D_{\mp}g$ exist in a neighborhood of the support of $f$, then
\be\label{2.20}
 \intl_{\p} (D_{\pm}f)(r)g(r)dr=\intl_{\p}f(r)(D_{\mp}g)(r)dr.
 \ee
  \end{lemma}
  \begin{proof}
Since $\supp\, f\subset
 \p$, one can replace integrals over $\p$ by those  over the whole space
 $\S_m \sim\bbr^{m(m+1)/2}$. The function $g$ can be multiplied
by a smooth cut-off function (in necessary) which $\equiv 1$ on
the support of $f$.
   Then the ordinary integration by parts yields
 (\ref{2.20}).
\end{proof}

\begin{lemma} \hfil

\noindent {\rm (i)}   For  $s \in \p$ and $z  \in \S_m^{\bbc}$,
\be\label{3.41} D_{+,\,s} [\exp(-\tr(sz))]=(-1)^m
\det(z)\exp(-\tr(sz)). \ee

\noindent {\rm (ii)}  For  $s \in \p$ and $\a  \in \bbc$,
\be\label{3.40} D_{+}(|s|^{\a-d}/\gma)=
|s|^{\a-1-d}/\gm(\a-1),\qquad d=(m+1)/2.\ee
\end{lemma}

\begin{proof} These statements are  well known; see, e.g., \cite[p. 813]{Ga},
 \cite[ p. 481]{Herz}, \cite[p. 125]{FK}. We recall the proof for convenience of the reader.
 The  equality  (\ref{3.41}) is verified by direct
calculation. Furthermore, by (\ref{2.19}),
 \be\label{2.25}
\gm (\b) |s|^{-\b}=\intl_{\p} \exp(-\tr(rs))|r|^{\b-d} dr,\qquad
Re\;\b>d-1.\ee Applying $ D_{+}$ to (\ref{2.25}) and using
(\ref{3.41}),  we get $$ D_{+}( |s|^{-\b}\gm (\b))=(-1)^m
|s|^{-\b-1}\gm (\b+1),\qquad Re\;\b>d-1.$$ Now we set $\b=d-\a$.
Since
$$\gm (d-\a+1)/\gm (d-\a)=(-1)^m \gm (\a)/\gm (\a-1),$$
 (\ref{3.40}) follows for $Re\;\a<1$. By analyticity,  it holds generally.
\end{proof}
\begin{remark} By changing notation and using (\ref{Poh}), one can write (\ref{3.40})
as \be\label{D-det} D_{+}|s|^{\a}=b(\a) |s|^{\a-1}, \ee where
\be\label{B1} b(\a)=\a(\a+1/2)\cdots(\a+d-1), \qquad d=(m+1)/2,
\ee is the so-called {\it Bernstein polynomial} \index{b}{Latin
and Gothic!ba@$b(\a)$} \index{a}{Bernstein polynomial, $b(\a)$} of
the determinant \cite{FK}. It is worth noting that $b(\a)$ can be
written in different forms, namely, \be\label{B2} b(\a)= (-1)^m
b(1-d-\a)= 2^{-m}\Gam(2\a+m)/\Gam(2\a)=2^{-m}(2 \a,m)\ee or
\be\label{B3} b(\a)= (-1)^m
\gm(1-\a)/\gm(-\a)=\gm(\a+d)/\gm(\a+d-1).\ee
\end{remark}

\begin{lemma}
 For $f\in\D(\p)$ and  $z  \in \S_m^{\bbc}$,
\be\label{2.21}
 (LD_+f)(z)=\det(z)(Lf)(z).
 \ee

  \end{lemma}
  \begin{proof}
By (\ref{2.20}), (\ref{2.51}) and (\ref{3.41}),

$$(LD_+f)(z)=\intl_{\p}D_{-,\, r} [\exp(-\tr(rz))]f(r)dr=\det(z)(Lf)(z),$$
as desired.
\end{proof}

\section{Bessel functions of matrix argument}
 We recall some facts from
\cite{Herz} and \cite{FK}.  {\it The $\J$-Bessel function}
$\J_\nu(r)$, \index{b}{Latin and
Gothic!jn@$\J_\nu(r)$}\index{a}{Bessel function, $\J_\nu(r)$}
$r\in\p$, can be defined in terms of the Laplace transform by the
property \be \intl_{\p}
\exp(-\tr(zr))\J_\nu(r)|r|^{\nu-d}dr=\gm(\nu)\exp(-\tr(z^{-1}))\det(z)^{-\nu},
\ee
$$ d=(m+1)/2,\qquad z\in\T_m=\p +i\S_m .$$ This gives
\be\label{Bes} \frac{1}{\gm(\nu)}\;\J_\nu(r)=\frac{1}{(2\pi
i)^N}\intl_{Re\, z=\sig_0} \exp(\tr(z- rz^{-1}))\det(z)^{-\nu}
dz,\quad \sig_0\in\p, \ee $N=m(m+1)/2$. This integral   is
absolutely convergent for $Re\,\nu>m$ and extends  analytically to
all $\nu\in\bbc$. Moreover, the formula (\ref{Bes}) allows us to
extend $\J_\nu(r)$ in the $r$-variable as an entire function
$\tilde\J_\nu(z)$, $z=(z_{i,j})_{m\times m}\in \bbc^{m^2}$, so
that $\tilde\J_\nu(z)\big |_{z=r}=\J_\nu(r)$ and \be\label{sym}
\tilde\J_\nu(\g z\g')=\tilde\J_\nu(z)\quad \mbox{ for all}
\quad\g\in O(m). \ee Functions satisfying (\ref{sym}) are called
symmetric. If $\lam_1,\dots, \lam_m$ are eigenvalues of $z$, and
$$\sig_1=\lam_1+\ldots +\lam_m,\quad \sig_2=\lam_1\lam_2+\ldots
+\lam_{m-1}\lam_m,\quad\dots\; ,\quad \sig_m=\lam_1\dots \lam_m$$
are the corresponding elementary symmetric functions (which are
polynomials of $z_{i,j}$), then  each analytic symmetric function
of $z\in\bbc^{m^2}$ is  an analytic function of $m$ variables
$\sig_1, \dots, \sig_m$. Since for any $r \in\p$ and $s\in\p$, the
matrices
$$
rs,\qquad sr,\qquad s^{1/2}rs^{1/2}, \qquad r^{1/2}sr^{1/2}
$$
have the same eigenvalues, then
$$
\tilde\J_\nu(rs)=\tilde\J_\nu( sr)=\tilde\J_\nu(
s^{1/2}rs^{1/2})=\tilde\J_\nu(r^{1/2}sr^{1/2}).
$$
In the following, we keep the usual notation $\J_\nu(z)$ for the
extended function $\tilde\J_\nu(z)$.

For $m=1$, the classical Bessel function $J_\nu(r)$
\index{b}{Latin and Gothic!j@$J_\nu(r)$} expresses through
$\J_\nu(r)$ by the formula
$$
J_\nu(r)=\frac{1}{\Gam(\nu+1)}\left(\frac{r}{2}\right)^\nu
\J_{\nu+1}\left(\frac{r^2}{4}\right).
$$

There is an intimate connection between the Fourier transform and
$\J$-Bessel functions.
\begin{theorem} \cite[p. 492]{Herz}, \cite[p.
355]{FK} If $f(x)$ is an integrable  function  of the form
$f(x)=f_0(x'x)$ where $f_0$ is a function on $\p$, then  \be
\label{Boh}(\F f)(y)\equiv\intl_{\Ma} \exp(\tr(iy'x)) f_0 (x'x)
dx= \frac{\pi^{nm/2}}{\gm(n/2)}\tilde f_0
\left(\frac{y'y}{4}\right), \ee where \be\label{tf0} \tilde f_0
(s)= \intl_{\p}\J_{n/2} (rs)|r|^{n/2-d}f_0(r) dr, \qquad
d=(m+1)/2\ee (the Hankel transform of $f_0$).
\end{theorem}

This statement  represents a matrix generalization of the
classical result of Bochner, see, e.g., \cite[Chapter IV, Theorem
3.3]{SW} for the case $m=1$. \index{a}{Bochner formula} Since the
product $rs$ in (\ref{tf0})  may not be a positive definite
matrix, the expression $\J_{n/2} (rs)$ is understood in the sense
of analytic continuation explained above.

 Note that one can meet a different
notation for the $\J$-Bessel function of matrix argument in the
literature. We follow the notation from \cite{FK}. It relates to
the  notation $A_{\del} (r)$ in \cite{Herz} by the formula
$\J_\nu(r)=\gm(\nu) A_{\del} (r)$, $\del=\nu-d$. More information
about  $\J$-Bessel functions and their generalizations can be
found in \cite{Dib}, \cite{FK}, \cite{FT}, \cite{GK1}, \cite{GK2}.

\section{Stiefel manifolds}\label{s2.3}

For $n\geq m$, let $\vnm= \{v \in \frM_{n,m}: v'v=I_m \}$
\index{b}{Latin and Gothic!vnm@$\vnm$}
 be  {\it the Stiefel manifold} \index{a}{Stiefel manifold, $\vnm$} of orthonormal $m$-frames in $\bbr^n$.
 If $n=m$, then $V_{n,n}=O(n)$ is the orthogonal group in $\bbr^n$.
 It is known  that $\dim\vnm = m(2n-m-1)/2$ \cite[ p. 67]{Mu}.
 The group $O(n)$
 acts on $\vnm$ transitively by the rule $g: v\to gv, \quad g\in
 O(n)$, in the sense of matrix multiplication.  The same is true for
  the special orthogonal group  $SO(n)$  provided $n>m$.
We  fix the corresponding  invariant measure $dv$ on
 $\vnm$
  normalized by \be\label{2.16}\index{b}{Greek!ma@$\sigma_{n,m}$} \sigma_{n,m}
 \equiv \intl_{\vnm} dv = \frac {2^m \pi^{nm/2}} {\gm
 (n/2)}, \ee
\cite[p. 70]{Mu},   \cite[p.
 57]{J}, \cite[ p. 351]{FK}. This measure is also $O(m)$
 right-invariant.

\begin{lemma}

\cite[p. 354]{FK}, \cite[p. 495]{Herz}. The Fourier transform of
the invariant  measure $dv$ on $\vnm$ represents the $\J$-Bessel
function: \be\label{f-st} \intl_{\vnm} \exp(\tr(iy'v)) dv
=\sig_{n,m} \J_{n/2}\left(\frac{y'y}{4}\right). \ee
\end{lemma}

\begin{lemma}\label{l2.3} {\rm (polar decomposition).}\index{a}{Polar
decomposition}
Let $x \in \frM_{n,m}, \; n \ge m$. If  $\rank (x)= m$, then \[
x=vr^{1/2}, \qquad v \in \vnm,   \qquad r=x'x \in\p,\] and
$dx=2^{-m} |r|^{(n-m-1)/2} dr dv$.
\end{lemma}

This statement can be found in many sources, see, e.g., \cite[ p.
482]{Herz}, \cite[p. 93]{GK1}, \cite[pp. 66, 591]{Mu}, \cite[p.
130]{FT}. A modification of Lemma \ref{l2.3} in terms of upper
triangular matrices $t\in T_m$ (see (\ref{2.17})) reads as
follows.

\begin{lemma}\label{sph}\index{a}{Polar
decomposition} Let $x \in \frM_{n,m}, \; n \ge m$.
 If  $\rank (x)= m$, then
  \[
x=vt, \qquad v \in \vnm,   \qquad t \in T_m,\] so that
$$
dx=\prod\limits_{j=1}^m t_{j,j}^{n-j} dt_{j,j}\,dt_\ast\, dv,
\qquad dt_\ast=\prod\limits_{i<j} dt_{i,j}.
$$
\end{lemma}

This statement is also well known and has different proofs. For
instance, it can be easily derived from Lemma \ref{l2.3} and
(\ref{2.2}); see Lemma 2.7 in \cite{Ru10}.

\begin{lemma}\label{l2.1} {\rm (bi-Stiefel decomposition).} \index{a}{Bi-Stiefel decomposition} \ Let $k$,
$m$, and $n$  be positive integers satisfying $$1 \le k  \le n-1,
\qquad 1 \le m  \le n-1, \qquad k+ m \le n.$$

\noindent {\rm (i)} \ Almost all matrices $v \in \vnm$ can be
represented in the form
 \be\label{herz1}
  v= \left[\begin{array} {cc} a \\ u(I_m -a'a)^{1/2}
 \end{array} \right], \qquad a\in \frM_{k, m}, \quad u \in V_{n-k,m},
 \ee
  so that
   \be\label{2.10}
   \intl_{\vnm} f(v) dv
 =\intl_{0<a'a<I_m}d\mu(a)
 \intl_{ V_{n- k, m}} f\left(\left[\begin{array} {cc} a \\ u(I_m -a'a)^{1/2} \end{array}
 \right]\right) \, du,
\ee
 $$
 d\mu(a)=|I_m-a'a|^\delta da,\quad \del=(n-k)/2-d, \quad  d=(m+1)/2.
 $$

 \noindent  {\rm (ii)} \ If, moreover, $k \ge m$, then
\be\label{2.11}
 \intl_{\vnm} f(v) dv
 =\intl_0^{I_m} d\nu(r) \intl_{ V_{k, m}}dw
 \intl_{ V_{n-k, m}} f\left(\left[\begin{array} {cc} wr^{1/2} \\ u(I_m -r)^{1/2} \end{array}
 \right]\right) \ du,
\ee \be\nonumber
 d\nu(r)=2^{-m} |r|^\gam |I_m -r|^\del dr, \qquad \gam
=k/2-d.
 \ee
\end{lemma}

 \begin{proof}For m=1, this  is a well known  bispherical decomposition \cite[pp. 12,
 22]{VK2}.
 For  $k=m$, this statement is due to  \cite[p. 495]{Herz}.
 The proof of Herz was extended in \cite{GR} to all $k+m \leq n$.
  For convenience of the reader, we reproduce this proof in our notations which differ from those in
  \cite{GR}.

 Let us check (\ref{herz1}). If $v=\left[\begin{array} {ll} a \\ b
 \end{array}\right] \in \vnm, \quad a \in \frM_{k,m}, \; b \in
 \frM_{n-k,m}, $ then $I_m=v'v=a'a+b'b$, i.e., $b'b=I_m-a'a$.
  By Lemma \ref{l2.3}, for almost all $b$ (specifically, for all $b$ of rank $m$), we have  $b= u(I_m -a'a)^{1/2}$ where
  $u \in V_{n- k, m}$.  This gives (\ref{herz1}).
In order to prove (\ref{2.10}), we  show  the coincidence of the
two measures, $dv$ and $\tilde{d}v=|I_m -a'a|^\del
 da du$.  Following
 \cite{Herz},  we consider the Fourier transforms
\[
 F_1(y)=\intl_{\vnm}  \exp(\tr(iy'v))  dv \quad \mbox{and} \quad
 F_2(y)=\intl_{\vnm} \exp(\tr(iy'v))\tilde{ d}v,
\]
 $y \in \frM_{n,m}$, and
 show
that $F_1=F_2$.  Let \[
 y=\left[\begin{array} {c} y_1 \\ y_2 \end{array} \right], \qquad
 v=\left[\begin{array} {c} a \\ u(I_m -a'a)^{1/2}
 \end{array} \right], \]
where
\[ y_1
  \in \frM_{k,m}, \qquad  y_2  \in \frM_{n-k,m}; \qquad
 a \in \frM_{k,m}.
\]  Then $y'v=y'_1a + y'_2u (I_m -a'a)^{1/2}$, and we have
\bea F_2(y)&=&\intl_{a'a<I_m} \exp(\tr (iy'_1a))
 |I_m -a'a|^\del da \nonumber \\
&\times&\intl_{
 V_{n- k,m}} \exp(\tr (iy'_2 u (I_m -a'a)^{1/2}))
 du. \nonumber
\eea By (\ref{f-st}),  the inner integral is evaluated as
\[
\sig_{n-k,m} \J_{(n-k)/2}(\frac{1}{4}y'_2y_2(I_m -a'a)).
\]
 Thus
\[
F_2(y)=\intl_{a'a<I_m} \exp (\tr (iy'_1a)) \vp(a'a) da, \] \[ \vp
(r)=
 \sig_{n-k,m} |I_m -r|^{(n-k)/2-d} \J_{(n-k)/2}(\frac{1}{4}y'_2y_2(I_m -r)).
\]
The function $F_2(y)$ can be transformed by the generalized
 Bochner formula (\ref{Boh}) as follows
\bea \nonumber F_2(y)&=&\frac{\pi^{km/2}\sig_{n-k,m}}{\gm(k/2)}
\intl_0^{I_m} \J_{k/2} (\frac{1}{4} y'_1 y_1 r)|r|^{k/2-d}\\
\nonumber &\times& \J_{(n-k)/2} (\frac{1}{4} y'_2 y_2(I_m -r))
|I_m -r|^{(n-k)/2-d}
 dr.
\eea The last integral can be
 evaluated by
 the formula
\be \intl_0^s \J_\a(pr)|r|^{\a-d}\J_\b(q(s-r))|s-r|^{\b-d} dr=
B_m(\a,\b)\J_{\a+\b}((p+q)s)|s|^{\a+\b-d}, \ee

$$
Re \, \a, Re \, \b >d-1, \qquad d=(m+1)/2,\qquad p,q\in\cpm,
\qquad  s\in\p,
$$
cf.  \cite[formula (2.6)]{Herz}.  Applying this formula with
$\a=k/2$, $\b=(n-k)/2$, and $s=I_m$, we obtain
 \bea
 F_2(y)&=&\frac{\pi^{km/2}\sig_{n-k,m}}{\gm(k/2)}B_m\left(\frac{k}{2},\frac{n-k}{2}\right)\J_{n/2} (\frac{1}{4} (y'_1 y_1  +  y'_2
 y_2)) \nonumber \\
 &=& \sig_{n,m}\J_{n/2} (\frac{1}{4} y'y).\nonumber
\eea By (\ref{f-st}), this coincides with $F_1(y)$, and
(\ref{2.10}) follows. The equality (\ref{2.11}) is a consequence
of (\ref{2.10}) by Lemma \ref{l2.3}.
\end{proof}

\section{Radial functions and the Cayley-Laplace operator}\label{s2.4}
The notion of a radial function usually alludes to invariance
under rotations. In the classical analysis on $\bbr^n$ the
 radial function $f(x)$ is virtually a function of a non-negative
 number $r=|x|$. In the matrix case the situation is similar {\it but now $r$ is   a
  positive semi-definite   matrix}.

\begin{definition}\label{d5.1}
A function $f(x)$ on $\Ma$  is called radial, if it is $O(n)$
left-invariant, i.e., \be\label{4.13}f(\gamma x)=f(x), \quad
\forall \gamma\in O(n).\ee If $f$ is continuous, this equality is
understood ``for all $x$". For generic measurable functions it is
 interpreted in the almost everywhere sense.

\end{definition}

\begin{lemma} Let  $x\in\Ma$. Each function of the form $f(x)=\f0
(x'x)$ is radial. Conversely, if $n\geq m$, and  $f(x)$ is a
radial function, \index{a}{Radial functions!on $\Ma$} then there
exists $\f0 (r)$ on $\p$ such that $f(x)=\f0 (x'x)$ for all (or
almost all)  $x$.
\end{lemma}
\begin{proof} The first statement follows immediately from Definition
\ref{d5.1}. To prove the second one, we note that   for
 $n\geq m$, any matrix $x\in\Ma$  admits
representation $ x=vr^{1/2}$, $v \in \vnm$,   $r=x'x \in\cpm$; see
 Appendix C, {\bf 11}.  Fix any frame $ v_0\in\vnm$ and choose
$\gamma\in O(n)$ so that $\gamma v_0=v$. Owing to (\ref{4.13}),
$$
f(x)=f(vr^{1/2})=f(\gamma v_0r^{1/2})=f(v_0r^{1/2})\stackrel{\rm
def}{=}\f0(r),
$$
as desired. Clearly, the result is independent of the choice of $
v_0\in\vnm$.
\end{proof}

{\it The Cayley-Laplace operator} $\Del$
\index{b}{Greek!del@$\Del$} \index{a}{Cayley-Laplace operator!on
$\Ma$, $\Del$} on the space of matrices $x=(x_{i,j})\in \Ma$ is
defined by \be\label{K-L} \Del=\det(\d '\d). \ee Here, $\partial$
is an $n\times m$  matrix whose entries are partial derivatives
$\d/\d x_{i,j}$. In terms of the Fourier transform, the action of
$\Del$ represents a multiplication by the polynomial $(-1)^m
P(y)$, $y\in\Ma$, where
$$
P(y)=|y'y|=\det \left[\begin{array}{cccc} y_1 \cdot y_1 & {.}   & {.}   & y_1 \cdot y_m \\
                              {.} & {.}  & {.}   & {.} \\
                              {.} & {.}   & {.}   & {.} \\
                               y_m  \cdot y_1& {.}     & {.} & y_m \cdot y_m

\end{array} \right],
$$
$ y_1, \dots , y_m$ are column-vectors of  $y$, and $``\cdot"$
stands for the usual inner product in $\rn$. Clearly, $P(y)$ is a
homogeneous polynomial of degree $2m$ of $nm$ variables $y_{i,j}$,
and $\Del$ is a homogeneous differential operator of order $2m$.
For $m=1$, it coincides with the Laplace operator on $\rn$.

Operators (\ref{K-L}) and their generalizations were studied by
S.P. Khekalo \cite{Kh}.  For $m>1$, the operator $\Del$ is not
elliptic because $P(y)=0$  for all non-zero matrices $y$ of $\rank
<m$. Moreover, $\Del$ is not hyperbolic, \cite[p. 132]{Ho},
 although,
for some $n, m$ and $\ell$, its power $\Del^\ell$ enjoys the
strengthened Huygens' principle; see \cite{Kh} for details.

The following statement gives  explicit representation of the
radial part of the Cayley-Laplace operator.
\begin{theorem}\label{TKL-r}{\rm \cite{Ru10}.}
Let $f_0(r)\in C^{2m}(\p)$, $f(x)=f_0(x'x)$, $d=(m+1)/2$. Then for
each matrix $x\in\Ma$  of rank $m$,\be \label{KL-r}(\Del
f)(x)=(Lf_0)(x'x),\ee where \be\label{L-r} L=4^m |r|^{d-n/2}D_+
|r|^{n/2-d+1}D_+ , \ee $D_+$ being the operator (\ref{2.50}).
\end{theorem}

\begin{example}
Let $f(x)=|x|_m^\lam, \; |x|_m=\det (x'x)^{1/2}$. By Theorem
\ref{TKL-r}, we have $(\Del f)(x)=\vp(x'x)$, where $\vp(r)=4^m
|r|^{d-n/2}D_+ |r|^{n/2-d+1}D_+|r|^{\lam/2}$. This expression can
be evaluated using (\ref{D-det}): \bea\nonumber \vp(r) &=& 4^m
b(\lam/2)
|r|^{d-n/2}D_+ |r|^{(n+\lam)/2-d}\\[14pt]\nonumber&=& 4^m
b(\lam/2) b((n+\lam)/2-d) |r|^{\lam/2-1}.
 \eea
Thus, we have arrived at the following identity of the Bernstein
type  \be \label{Dxm}\Del |x|_m^\lam =\B(\lam) |x|_m^{\lam -2},\ee
(cf. \cite[p. 125]{FK}) where, owing to (\ref{B1}), the polynomial
$\B(\lam)$ has the form \be\label{B} \B(\lam)=(-1)^m
\prod\limits_{i=0}^{m-1}(\lam+i)(2-n-\lam+i). \ee

An obvious consequence of (\ref{Dxm}) in a slightly different
notation reads \be\label{vaz}\Del ^k
|x|_m^{\a+2k-n}=B_k(\a)|x|_m^{\a-n},\ee where \bea\label{bka}
B_k(\a)&=&\prod\limits_{i=0}^{m-1}\prod\limits_{j=0}^{k-1}(\a-i+2j)(\a-n+2+2j+i)
\\ &=& B_k(n-\a-2k). \nonumber \eea
\end{example}

\chapter{The G{\aa}rding-Gindikin fractional integrals}\label{s3}

\section{Definitions and comments} \label{s3.2}
In this section, we
 begin detailed investigation of  fractional
integrals associated to the cone $\p$ of positive definite $m
\times m$ matrices. These integrals were introduced by L.
G{\aa}rding \cite{Ga} and substantially generalized by S. Gindikin
\cite{Gi1}.

Let $f(r)$ be a smooth function on $\p$ which is rapidly
decreasing at infinity and bounded with all its derivatives when
$r$ approaches the boundary $\partial\p$. Consider the integral
\be\label{Ia} I_f(\a)=\intl_{\p} f(r)|r|^{\a-d} dr, \qquad
d=(m+1)/2.\ee If $f(r)=\exp(-\tr (r))$ then \[ I_f(\a)=\gm
(\a)=\pi^{m(m-1)/4}\prod\limits_{j=0}^{m-1} \Gam (\a- j/2)\; ; \]
see (\ref{2.4}), (\ref{2.5}). In the general case, by passing to
upper triangular matrices and using (\ref{2.2}), we have \bea
\nonumber I_f(\a)&=& \intl_0^\infty  t_{1,1}^{2\a-1} \,
dt_{1,1}\intl_0^\infty t_{2,2}^{2\a-2} \, dt_{2,2} \,
 \ldots  \intl_0^\infty  t_{m,m}^{2\a-m}  \tilde f (t_{1,1}, \ldots ,t_{m,m}) \,
 dt_{m,m}\\[14pt]\nonumber
&=& 2^{-m}\intl_0^\infty  y_{1}^{\a-1} \, dy_{1}\intl_0^\infty
y_{2}^{\a-3/2} \, dy_{2} \,
 \ldots  \intl_0^\infty  y_{m}^{\a-(m+1)/2}  \tilde f (y_{1}^{1/2}, \ldots ,y_{m}^{1/2}) \,
 dy_{m},\eea
\[
 \tilde{ f} (t_{1,1}, \ldots ,t_{m,m})=2^m  \intl_{\bbr^{m(m-1)/2}} f(t't) \, dt_{*}, \quad dt_{*}=\prod_{i<j}
 dt_{i,j}.
\]
The function $ \tilde f$ extends to the whole space $\bbr^m$ as an
even function in each argument, and the function  \[ f_0 (y_1,
\ldots , y_m) \equiv  2^{-m}
 \tilde f (y_{1}^{1/2}, \ldots ,y_{m}^{1/2})\] belongs to
the  Schwartz space  $\s(\bbr^m)$. Hence  $I_f(\a)$ can be
regarded as a direct product of one-dimensional distributions
\be\label{dpr} I_f(\a)=\left(\prod\limits_{j=1}^m
(y_j)_+^{\a-(j+1)/2},  f_0 (y_1, \ldots , y_m)\right ). \ee
 It follows  \cite[Chapter 1, Section 3.5]{GSh1}
that  the integral (\ref{Ia}) absolutely converges if and only if
$Re\, \a>(m-1)/2$ and  extends as a meromorphic function of $\a$
with the only poles $(m-1)/2, \; (m-2)/2,\dots\;$. These poles and
their orders coincide with those  of the gamma function $\gm(\a)$
so that  $I(\a)/\gma$ is an entire function of $\a$.

These preliminaries motivate the following definitions. We first
consider ``good'' functions $f\in \D (\p)$ which are infinitely
differentiable and compactly supported away form the boundary
$\partial \p$.  Then we proceed to the case of  arbitrary locally
integrable functions.

 \begin{definition}\label{d3.1} \index{a}{G{\aa}rding-Gindikin fractional integrals!left-sided, $I_{+}^\a f$}

 Let $f\in \D (\p)$,   $ \a \in \bbc$, $d=(m+1)/2$. The {\it left-sided  G{\aa}rding-Gindikin fractional integral}  is
 defined by

\be\label{3.1} (I_{+}^\a f)(s) = \left \{
\begin{array} {ll}\displaystyle{\frac {1}{\gma} \intl_0^s f(r)|s-r|^{\a-d}
dr}\qquad
  \mbox{ if   $Re \, \a > d-1$}, \\
{} \\
 (I_{+}^{\a+\ell}D_{+}^\ell f)(s)
  \qquad  \mbox{if $d-1-\ell<Re \, \a \le d-\ell$;} \\
 {}  \qquad \qquad    \qquad \qquad  \qquad  \mbox{$\ell=1,2, \ldots$}\;.
\end{array}
\right.  \ee Here,  $s \in \p$, $D_+$ is the differential operator
(\ref{2.50}),
 and we integrate over
the ``interval" $(0,s)=\{r: \; r\in\p, \; s-r\in\p \}.$
\index{b}{Latin and Gothic!ia@$I^\a_{\pm}f$}
\end{definition}

 \begin{definition}\label{d3.2}\index{a}{G{\aa}rding-Gindikin fractional integrals!right-sided, $I_{-}^\a f$}
 For $f\in \D (\p)$ and  $ \a \in \bbc$, the {\it right-sided G{\aa}rding-Gindikin fractional integral}  is
 defined by

 \be\label{3.2} (I_{-}^\a f)(s) = \left \{
\begin{array} {ll}\displaystyle{\frac
{1}{\gma} \intl_s^\infty f(r)|r-s|^{\a-d} dr}
 \qquad \mbox{ if   $Re \, \a > d-1$}, \\
{} \\
 (I_{-}^{\a+\ell}D_{-}^\ell f)(s)
  \qquad  \mbox{if $d-1-\ell<Re \, \a \le d-\ell$;} \\
 {}  \qquad \qquad    \qquad \qquad  \qquad  \mbox{$\ell=1,2, \ldots$}\;.
\end{array}
\right.
 \ee
Here, $s$ is a positive semi-definite matrix in $\cpm$, $D_-$ is
defined by (\ref{2.51}), and $\int_s^\infty$ denotes integration
over the shifted  cone $s+\p=\{r: \; r\in\p :r-s\in\p \}$. The
vertex $s$ of this cone may be an inner point of $\p$ or its
boundary point.
 \end{definition}
 In the rank-one case  $m=1$, when $\p$ is a
positive half-line, (\ref{3.1}) and (\ref{3.2}) are ordinary
Riemann-Liouville (or Abel) fractional integrals \cite{Ru1},
\cite{SKM}.

\begin{remark} Analogous   definitions  can be given for functions on the
ambient space $\S_m \supset \p$, for instance, for $f\in\s(\S_m)$.
 We do not do this
for two reasons. Firstly, because our purposes are exclusively
connected with analysis on the cone $\p$, and the space  $\D(\p)$
fits our needs completely. Secondly, dealing with the space
$\D(\p)$, we can still reveal basic features of our objects and
give much simpler proofs then for other classes of functions.
\end{remark}

According to Definitions \ref{d3.1} and \ref{d3.2},  it is
convenient to regard the complex plane as a union
$\bbc=\cup_{\ell=0}^\infty \Omega_\ell$, where \be\label{om}
\Omega_\ell=\left \{
\begin{array} {ll} \quad \{\a : Re \, \a >d-1 \} & \mbox{ if   $\ell =0$}, \\
{} \\ \quad \{\a : d-1-\ell <Re \, \a \le d-\ell \}  &\mbox{ if
$\ell =1,2, \ldots $.}
\end{array}
\right.
 \ee

\begin{lemma}\label{ent} For $f\in \D (\p)$, the integrals $(I_{\pm}^\a f)(s)$   are entire functions of $\a$.
\end{lemma}
\begin{proof} Let us consider the integrals $(I_{+}^\a f)(s)$.  Given an integer $j>0$, we set
$F_j(\a)=I_{+}^{\a+j}D_{+}^j f$. This function is analytic in the
half-plane $Re \, \a >d-1-j$. Owing to uniqueness property of
analytic functions, it suffices to show that $F_j(\a)=I_{+}^{\a}f$
on each strip $ \Omega_\ell$ provided that $j=j(\ell)$ is large
enough. For $\a \in  \Omega_\ell$, we have $\a+j-d>j-\ell -1$, and
therefore, by (\ref{2.20}) and (\ref{3.40}),
 \bea F_j(\a)&=&\frac {1}{\gm(\a+j)}
\intl_{\S_m} (s-r)_+^{\a+j-d} (D_{+}^{j}f)(r) dr \nonumber \\
&=&\frac {1}{\gm(\a+j)}
\intl_{\S_m} D_{-,\, r}^{j-\ell} (s-r)_+^{\a+j-d} (D_{+}^{\ell}f)(r) dr \nonumber \\
&=&\frac {1}{\gm(\a+\ell)}\intl_0^s |s-r|^{\a+\ell-d}
(D_{+}^{\ell}f)(r) dr\nonumber \\
&=&(I_{+}^{\a+\ell}D_{+}^\ell f)(s). \nonumber \eea By Definition
\ref{d3.1}, this means that  $F_j(\a)=I_{+}^{\a}f$ for  $\a \in
\Omega_\ell$, and we are done. For operators $I_{-}^{\a}$ the
proof is similar.
\end{proof}

 \begin{lemma}
Let $f\in \D (\p)$. Then \be\label{3.49} D_{\pm}^j I_{\pm}^\a
f=I_{\pm}^\a D_{\pm}^j f=I_{\pm}^{\a-j} f, \ee
 for any
$\a \in \bbc$ and $j \in \bbn$.
\end{lemma}
\begin{proof}
Fix $j\in\bbn$. For $Re\, \a>j+d-1$, the  equality $I_{\pm}^\a
D_{\pm}^j f=I_{\pm}^{\a-j} f$ can be obtained using integration by
parts (see the proof of Lemma \ref{ent}). For all $\a\in\bbc$ it
then follows  by analytic continuation. Let us show that
$D_{\pm}^j I_{\pm}^\a f=I_{\pm}^{\a-j} f$.
 Let $Re \, \a -j \in \Omega_\ell \, $; see (\ref{om}). If
$\ell=0$, then \bea I_{+}^{\a-j} f &=&\frac
{1}{\gm(\a -j)}\intl_0^s |s-r|^{\a-j-d} f(r) dr\nonumber \\
&=&\frac {1}{\gm(\a)}
\intl_{\S_m} D_{+,\, s}^{j} (s-r)_+^{\a-d}  f(r) dr \nonumber \\
&=&(D_{+}^j I_{+}^\a f)(s). \nonumber \eea If $\ell \ge 1$,
 then, by Definition \ref{d3.1},
 $I_{+}^{\a-j} f=I_{+}^{\a-j+\ell}D_{+}^\ell  f$. By the previous
 case, this coincides with $D_{+}^j I_{+}^{\a+\ell}D_{+}^\ell  f=D_{+}^j
 I_{+}^{\a}f$ (here we used Definition \ref{d3.1} again).  For the right-sided integrals the argument
  follows the same lines.
\end{proof}

Some   historical notes and comments are in order.
 Fractional integrals similar to
(\ref{3.1})--(\ref{3.2}) and associated to the light cone were
introduced by M. Riesz \cite{Ri1}, \cite{Ri2} in 1936. The
investigation
 was continued by N.E. Fremberg \cite{Fr1}, \cite{Fr2}, and by Riesz himself in
\cite{Ri3}, \cite{Ri4}; see also B.B. Baker and E.T. Copson
\cite{BC}, and E.T. Copson \cite{Co}. A key motivation of this
research was to find a simple form of the solution to the Cauchy
problem for the wave equation. M. Riesz' argument was essentially
simplified by
 J.J. Duistermaat \cite{Dui}; see also J. A.C. Kolk and V.S. Varadarajan \cite{KV}.

The results of Riesz were extended by L. G{\aa}rding \cite{Ga}
(1947), who replaced the light cone (of rank 2) by the  higher
rank cone $\p$ of positive definite symmetric matrices, and
studied the corresponding hyperbolic equations.\footnote{Prof.
Lars G{\aa}rding kindly informed the second-named author that his
construction of fractional integrals was inspired by close
developments in statistics.} The Riesz-G{\aa}rding fractional
integrals are very close to those in (\ref{3.1}) but do not
coincide with them completely.  The difference is as follows.  In
the above-mentioned papers, mainly devoted to PDE, the upper limit
$s$ in (\ref{3.1}) varies in the whole space $\bbr^N, \;
N=m(m+1)/2$ (not only in $\cpm$), and the domain of integration is
bounded by the retrograde cone $s-\p$ and a smooth surface on
which the Cauchy data are given. For $s \in \cpm$, this means that
the part of the domain of integration in (\ref{3.1}) containing
the origin is cut off smoothly.

Further progress is connected with fundamental  results of S.G.
Gindikin \cite{Gi1} (1964), who developed a general theory of
Riemann-Liouville integrals associated to  homogeneous cones. The
light cone of Riesz and the cone $\p$ of G{\aa}rding fall into the
scope of his theory. In the case of $\p$, Gindikin's fractional
integrals are  convolutions of the form \be\label{3.3} J^\a_{\pm}
f=\frac{s^{\a-d}_{\pm}}{\gma} \ast f, \qquad f \in \S(\bbr^N),
\quad N=\frac{m(m+1)}{2}, \ee where the function $s^{\a-d}_{+}$
equals $|s|^{\a-d}$ for $s \in \p$, zero otherwise, and
$s^{\a-d}_{-}$ equals $(-s)^{\a-d}_{+}$. Operators (\ref{3.3}) are
defined on the whole space $\bbr^N$, and can be treated using the
Fourier transform technique and the theory of distributions; see
\cite{Gi1}, \cite{VG}, \cite{Rab}, \cite{Wa}, \cite{FK}. See also
\cite{Rich} concerning application of such operators in
multivariate statistics.

Definitions \ref{d3.1} and \ref{d3.2} are motivated by
consideration of Radon transforms of functions of matrix argument
in the next chapters. It is important that the variable of
integration $r$ and the exterior variable $s$ do not range in the
whole space and are restricted to $\p$. In the following, the
class of functions $f$ will be essentially enlarged. We shall
define and investigate fractional integrals $I_{\pm}^\a f $ of
arbitrary locally integrable functions $f$ with possibly minimal
decay at infinity. The main emphasis is made upon $\a$ belonging
to the Wallach-like  set $ \W =\left\{0, \frac{1}{2}, 1,
\frac{3}{2}, \ldots , \frac{m-1}{2} \right\} \cup \left\{\a:\;
Re\,\a> \frac{m-1}{2} \right \}$,
 cf. \cite[ p. 137]{FK}.

 Thus, one can see
that most of the facts for fractional integrals from \cite{Ga} and
\cite{Gi1} cannot be automatically transferred to our case, and a
careful consideration is required.  The case $\a=k/2$, $k\in\bbn$,
 is of  primary importance because it arises in  our treatment of
 the Radon transform.

\section{Fractional integrals of half-integral order}\label{s2.2}
An important feature of the higher rank fractional integrals
$I_{\pm}^{\a}f$ is that for $\a=k/2$, $k=0, 1,\dots, m-1$, they
are convolutions with the corresponding  positive measures
supported on the boundary $\partial\p$. In the rank-one case
 $m=1$, this phenomenon  becomes trivial (we have only one value $\a=0$
corresponding to the delta function). For  $m>1$, when dimension
of the boundary $\partial\p$ is positive we deal with many ``delta
functions" supported by manifolds of different dimensions. This
phenomenon was drawn
 considerable attention in  analysis on symmetric domains; see
 \cite[ Chapter VII]{FK}, \cite{Gi2}, \cite{Ishi}, and references therein.

 Our nearest goal is to find precise form of $I_{\pm}^{k/2}f$ for
$k=1,\dots, m-1$. We start with the following simple observation.
 For $k\geq m$, by Lemma \ref{l2.3} and
(\ref{2.16}), one can write \bea \nonumber(I_{+}^{k/2}f)(s)&=&
\frac {1}{\gm(k/2)} \intl_0^s f(s-r)|r|^{k/2-d} dr
\\\label{2.22}&=& \pi^{-km/2}\intl_{\{\om\in\Mkm:\;\om'\om<s\}}
f(s-\om'\om)d\om,\\\nonumber (I_{-}^{k/2} f)(s) &=& \frac
{1}{\gm(k/2)} \intl_{\p} f(s+r)|r|^{k/2-d} dr\\\label{2.23} &
=&\pi^{-km/2}\intl_{\Mkm} f(s+\om'\om)d\om .\eea

The expressions (\ref{2.22}) and (\ref{2.23}) are meaningful  for
all $k\in\bbn$, and it is natural to conjecture that for
$k=1,\dots, m-1$, the integrals $I_{\pm}^{k/2}f$ are represented
by (\ref{2.22}) and (\ref{2.23}) too. Let us prove this
conjecture. We start with the left-sided integrals $I_{+}^{k/2}f$
and make use of the Laplace transform.

\begin{lemma}\label{l3.6}
For $f\in\D(\p)$ and $\a\in\bbc$,
 \be\label{3.9.1}
(LI_{+}^{\a}f)(z)=\det(z)^{-\a}(Lf)(z), \qquad Re\;z>0. \ee
\end{lemma}
\begin{proof} Let $\a \in \Omega_\ell \, $; see (\ref{om}). For $\ell =0$,
(\ref{3.9.1}) is an immediate consequence of (\ref{2.18}) and
(\ref{2.19}). The desired result for $\ell >1$  follows by
Definition \ref{d3.1} and (\ref{2.21}):
 $$
(LI_{+}^{\a}f)(z)=(LI_{+}^{\a+\ell}D_+^\ell f)(z)=\det(z)^{-\a
-\ell} (LD_+^\ell f)(z)=\det(z)^{-\a} (Lf)(z).$$
\end{proof}

\begin{corollary}\label{dif}
For $f\in\D(\p)$, \be\label{3.22}
(D_{+}^{\ell}f)(s)=(I_{+}^{-\ell} f)(s), \qquad  \ell= 0,1,2,
\ldots \, . \ee
\end{corollary}
\begin{proof}
By (\ref{2.21}) and (\ref{3.9.1}), $(LD_+^\ell
f)(z)=(LI_{+}^{-\ell} f)(z), \; Re \, z>0$. Since, by Definition
\ref{d3.1},  $I_{+}^{-\ell} f=I_{+}^d D_+^{\ell +d} f$ if
$d=(m+1)/2$ is an integer, and  $I_{+}^{-\ell} f=I_{+}^{d -1/2}
D_+^{\ell +d -1/2} f$ otherwise, the result follows owing to  the
uniqueness property of the  Laplace transform (see Lemma
\ref{lap}).

\end{proof}

\begin{theorem}
If $f\in\D(\p)$, then for all $k \in \bbn$,
 \be\label{2.2.2} \index{a}{G{\aa}rding-Gindikin fractional integrals!left-sided, $I_{+}^\a f$!of half-integral order} (I_{+}^{k/2}
 f)(s)=\pi^{-km/2}\intl_{\{\om\in\Mkm:\;\om'\om<s\}}
f(s-\om'\om)d\om .\ee Moreover,  \be\label{2.32}
 (I_{+}^{0} f)(s)=f(s). \ee
\end{theorem}
\begin{proof} Denote by  $\vp(s)$ the right-hand side of
(\ref{2.2.2}). By changing the order of integration, we have \bea
(L\vp)(z)&=&\pi^{-km/2}\intl_{\p}
\exp(-\tr(zs))ds\intl_{\{y\in\Mkm:\;y'y<s\}}f(s-y'y)dy \nonumber
\\ \label{2.24}&=&\pi^{-km/2}(Lf)(z)\intl_{\Mkm}
\exp(-\tr(zy'y))dy. \eea The integral in (\ref{2.24}) is easily
evaluated: if $z=r \in \p$, we replace  $y$ by $yr^{-1/2}$ and get
\bea\intl_{\Mkm} \exp(-\tr(ry'y))dy&=&|r|^{-k/2}\intl_{\Mkm}
\exp(-\tr(y'y))dy \nonumber \\ &=&|r|^{-k/2}\intl_{\bbr^{km}}
\exp(-y^2_{1,1}- \ldots -y^2_{k,m})dy \nonumber
\\&=&\pi^{km/2}|r|^{-k/2}. \nonumber \eea
By analytic continuation,
\[ \intl_{\Mkm}
\exp(-\tr(zy'y))dy= \pi^{km/2}\det(z)^{-k/2}, \qquad Re \, z>0.\]
Thus, $(L\vp)(z)=\det(z)^{-k/2}(Lf)(z)$,  and (\ref{3.9.1}) yields
$(L\vp)(z)= (LI_{+}^{k/2} f)(z)$, $Re \, z>0. $ Since for each
$\sig_0 \in \p$, the integrals \[ \int_{\p} \exp (-\tr (\sig_0 s)
|\vp (s)| ds, \qquad \int_{\p} \exp (-\tr (\sig_0 s)
|(I_+^{k/2}f)(s)| ds \] are finite (both are dominated by
$|\sig_0|^{-k/2} (L|f|)(\sig_0)$), by Lemma \ref{lap} we get
$\vp(s)=(I_{+}^{k/2} f)(s),$ as desired.

\end{proof}

\begin{remark} For $0<k<m$,  the integral
(\ref{2.2.2}) can be written as \be\label{3.12} (\I+k f)(s) \! =
\!
 c_{k,m}|s|^{k/2} \! \intl_{\vmk} \! dv \! \intl_{0}^{I_k} \!
  f(s \! - \! s^{1/2}vqv's^{1/2})|q|^{(m-k-1)/2}dq,\ee
  $$ c_{k,m}= \pi^{-km/2} 2^{-k}.$$
Indeed, if we transform (\ref{2.2.2}) by setting $\om=hs^{1/2}$,
where $h \in \Mkm, \; h'h<I_m$, and pass to polar coordinates, we
obtain (\ref{3.12}).
\end{remark}

 The integral
(\ref{2.2.2})  can also be  written as
 the Laplace convolution
\be\label{lc} (I_{+}^{k/2}f)(s)=\intl_0^s f(s-r) d\mu_k(r), \ee
where  $\mu_k$ is a positive measure defined by
 \be\label{3.42}
(\mu_k ,\psi)=\pi^{-km/2}\intl_{\Mkm} \psi(\om'\om)d\om, \qquad
\psi \in C_c(\S_m).
 \ee
Here, $C_c(\S_m)$ denotes the space of compactly supported
continuous functions on $\S_m$. For $k \ge m$, Lemma \ref{l2.3}
yields \[ (\mu_k ,\psi)=\frac{1}{\Gam_m (k/2)} \intl_{\p} \psi (r)
|r|^{k/2 -d} dr, \qquad d=(m+1)/2. \] In order to clarify
geometric structure of   $\mu_k$ for $k<m$, we denote
\index{b}{Greek!lamk@$\Lam_k$} \be\label{Lamj} \Lam_k =\{ s \, : s
\in \overline{\p}, \; \rank (s)=k\}, \qquad k=0,1, \ldots , m-1,
\ee
\[ G=GL(m, \bbr), \qquad e_k=\left[\begin{array}{ll}  I_k&0\\
0&0
\end{array}\right]\in\partial\p.\]
The manifold $\Lam_k$ is an orbit of $e_k$ under $G$ in $\S_m$,
i.e., $\Lam_k=Ge_k$. Indeed, $\rank (g'e_k g) =k$ for all $g \in
G$, and conversely, each matrix $s$ of rank $k$ is representable
as $s=g'e_k g$ for some $g \in G$; see, e.g., \cite[A6 (V) and
Theorem A9.4(ii)]{Mu} . The closure $\overline{\Lam_k}$ is a union
of $\Lam_1, \ldots , \Lam_k$. Since $\rank (\om'\om) \le k<m$,
then $(\mu_k ,\psi)=0$ for all $\psi$ supported away from
$\overline{\Lam_k}=\overline{Ge_k}$, and therefore $\supp \,\mu_k
=\overline{\Lam_k}$. Note also that $\mu_k
(\overline{\Lam_{k-1}})=0$ because
\[ \intl_{\overline{\Lam_{k-1}}} \psi (s) d\mu_k (s)=\intl_{\A_{k-1}}
\psi(\om'\om)d\om=0,\] $  \A_{k-1}=\{ \om \in \Mkm \, : \rank
(\om) \le k-1\}$ (the set $ \A_{k-1}$ has measure $0$ in $\Mkm$).
Thus the following statement holds.

\begin{theorem} For $0<k<m$, $\I+k f$ is
 the Laplace convolution (\ref{lc}) with a positive measure $\mu_k$ defined by
 (\ref{3.42}) and such that (a) $\supp\,
\mu_k=\overline{Ge_k}$, and (b) $\mu_k (\overline{Ge_{k-1}})=0$.
\end{theorem}

Let us consider the right-sided integrals $I_{-}^{\a}f$. They can
be expressed through the left-sided ones as follows.
\begin{lemma}\label{l3.10}
For $f\in\D(\p)$ and $\a\in\bbc$,
\be\label{3.7}(I_{-}^{\a}f)(s)=|s|^{\a-d}(I_{+}^{\a}g)(s^{-1}),\qquad
g(r)= |r|^{-\a-d}f(r^{-1}).
 \ee
In particular,
 \be\label{2.33}
 (I_{-}^{0} f)(s)=f(s). \ee
\end{lemma}
\begin{proof}
Since the integrals $(I_{+}^\a f)(s)$ and $ (I_{-}^\a f)(s)$  are
entire functions of $\a$ (see Lemma \ref{ent}), it suffices to
prove (\ref{3.7}) for $Re \, \a  >  d-1$. This can be easily done
by changing variables according to Lemma \ref{12.2} (iii).
\end{proof}
\begin{corollary}
For $f\in\D(\p)$ and $\ell\in\bbn$,
\be\label{3.7.1}(D_{-}^{\ell}f)(s)=|s|^{-\ell-d}(D_{+}^{\ell}
g)(s^{-1}),\qquad g(r)= |r|^{\ell-d}f(r^{-1}).
 \ee
\end{corollary}
\begin{proof}
The formula (\ref{3.7.1}) follows from (\ref{3.22}) and
(\ref{3.7}).

\end{proof}

Now we can  justify (\ref{2.23}) for all $k\in\bbn$.
\begin{theorem}\label{l3.11}
If $f\in\D(\p)$, $k \in \bbn$, then
 \be\label{2.29} \index{a}{G{\aa}rding-Gindikin fractional integrals!right-sided, $I_{-}^\a f$!of half-integral order} (I_{-}^{k/2}
 f)(s)=\pi^{-km/2}\intl_{\Mkm} f(s+\om'\om)d\om.
\ee
\end{theorem}
\begin{proof}
According to (\ref{3.7}),
\be\label{2.34}(I_{-}^{k/2}f)(s)=|s|^{k/2-d}(I_{+}^{k/2}g)(s^{-1}),\qquad
g(r)= |r|^{-k/2-d}f(r^{-1}).
 \ee
 By (\ref{3.12}),
$(I_{+}^{k/2}g)(s^{-1})$ is represented as \[
c_{k,m}|s|^{-k/2}\intl_{\vmk}dv\intl_{0}^{I_k}
g(s^{-1}-s^{-1/2}vqv's^{-1/2})|q|^{(m-k-1)/2} \, dq \]\[
=c_{k,m}|s|^{d}\intl_{\vmk}dv\intl_{0}^{I_k}
f((s^{-1}-s^{-1/2}vqv's^{-1/2})^{-1})\frac{|q|^{(m-k-1)/2} \, dq }
{|I_m-vqv'|^{k/2+d} }.\] Using the property \be \label{mu}
|I_m-vqv '|=|I_k-q|, \qquad q\in\pk, \qquad v\in\vmk,\ee (see
\cite[p. 575]{Mu} or Appendix C, {\bf1 }), this  can be written as
\[ c_{k,m}|s|^{d}\intl_{\vmk}dv\intl_{0}^{I_k}
  f(s^{1/2}(I_m-vqv')^{-1}s^{1/2}) \frac{|q|^{(m-k-1)/2} \, dq}{|I_k-q|^{k/2+d}} \]
  or \[c_{k,m}\sigma_{m,k}|s|^{d}\intl_{O(m)}d\g\intl_{0}^{I_k}
  f(s^{1/2}\g(I_m-v_0qv_0')^{-1}\g's^{1/2})    \frac{ |q|^{(m-k-1)/2} \,
  dq}{|I_k-q|^{k/2+d}},\]

$$ v_0 = \left[\begin{array} {c} 0 \\ I_k \end{array}
\right] \in V_{m,k}.$$ One can readily check that
\[
(I_m-v_0qv_0')^{-1}=
\left[ \begin{array} {ll} I_{m-k} &  0  \\
                             0 & I_k-q
\end{array}\right]^{-1}= \left[ \begin{array} {ll} I_{m-k} &  0  \\
                             0 & (I_k-q)^{-1}

\end{array}\right].
\]
  After changing variable
$(I_k-q)^{-1}\to q$, the integral $(I_{+}^{k/2}g)(s^{-1})$ becomes
\[ c_{k,m}\sigma_{m,k}|s|^{d}\intl_{O(m)}d\g\intl_{I_k}^{\infty}
  f \left(s^{1/2}\g\left[ \begin{array} {ll} I_{m-k} &  0  \\
                             0 & q
\end{array}\right]\g's^{1/2} \right)|I_k-q|^{(m-k-1)/2}dq,
\]
and therefore (replace  $q$ by $I_k+q$), \bea
(I_{+}^{k/2}g)(s^{-1}) \! \! &=& \! \! c_{k,m}|s|^{d}\intl_{\vmk}
\! dv\intl_{\pk}
  f(s^{1/2}(I_m+vqv')s^{1/2})|q|^{(m-k-1)/2}dq \nonumber \\
  &=&c_{k,m} 2^k |s|^{d}\intl_{\Mmk} f(s+s^{1/2}yy's^{1/2}) \,
  dy \qquad (y \! = \! s^{-1/2} z)\nonumber \\
  &=&c_{k,m} 2^k |s|^{d-k/2}\intl_{\Mmk} f(s+zz') \, dz\nonumber \\
  &=&\pi^{-km/2}|s|^{d-k/2}\intl_{\Mkm} f(s+\om'\om)d\om.\nonumber
  \eea
  Owing to  (\ref{2.34}), this gives the desired equality.
\end{proof}

\begin{remark} For $0<k<m$, the integral (\ref{2.29})
can be written in  polar coordinates as \be\label{3.13}
(I_{-}^{k/2}f)(s)=c_{k,m} \intl_{\vmk}dv\intl_{\pk}
 f(vqv'+s)|q|^{(m-k-1)/2}dq, \ee $$c_{k,m}= \pi^{-km/2}
 2^{-k}.$$
 Furthermore,
$$
(I_{-}^{k/2}f)(s)=\intl_{\p} f(s+r) d\mu_k(r),
$$
 $\mu_k$ being a positive measure defined by (\ref{3.42}).
\end{remark}

\section{The  G{\aa}rding-Gindikin distributions}\label{s2.3}

It is instructive to give   an alternative proof of the formula
(\ref{2.29}). This proof follows the argument from \cite[p.
134]{FK} and is of independent interest. It can be  useful in
different occurrences. Consider  the  G{\aa}rding-Gindikin
distribution \index{a}{G{\aa}rding-Gindikin distribution,
$\G_\a$}\index{b}{Latin and Gothic!ga@$\G_\a (f)$} \be\label{gd}
\G_\a (f)=\frac {1}{\gma} \intl_{\p} f(r)|r|^{\a-d} dr,\qquad
d=(m+1)/2, \ee where $f$ belongs to the Schwartz space $\s(\S
_m)$. The integral (\ref{gd}) converges absolutely for $Re\,\a >
d-1$ and
 admits  analytic
continuation as an entire function of $\a$ so that \be\label{r0}
\G_0(f) =f(0), \ee see \cite[pp. 132--133]{FK}. The following
statement implies (\ref{2.29}) and extends it for $f\in\s(\S _m)$.

\begin{lemma}\label{lrd} \index{a}{G{\aa}rding-Gindikin distribution, $\G_\a$!of half-integral order}
For $f\in\s(\S _m)$ and $0<k<m$, \be\label{er} \G_{k/2}
(f)=\pi^{-km/2}\intl_{\Mkm} f(\om'\om)d\om. \ee
\end{lemma}
\begin{proof}
 Let us transform
$\G_\a (f)$ to upper triangular matrices.
 By
(\ref{2.2}),
$$
\G_\a (f)=\frac{2^m}{\gm(\a)}\intl_{T_m} f(t't)
\prod\limits_{i=1}^{m} t_{i,i}^{2\a -i} \prod\limits_{i\leq j}
dt_{i,j}\;\qquad Re\,\a
>d-1  .
$$
We write $t=a+b$, where $a=(a_{i,j})$ and $b=(b_{i,j})$ are upper
triangular matrices so that the lower $n-k$ rows of $a$ and the
upper $k$ rows of $b$ consist of zeros. We denote by $A$ and $B$
the sets of all matrices of the form $a$ and $b$, respectively.
Since $t't=a'a+b'b$, then \be\label{dp} \G_\a (f)=\frac{2^k
\Gamma_{m-k}(\a-k/2)}{\gm(\a)}\intl_{A} g_\a
(a'a)\prod\limits_{i=1}^{k} a_{i,i}^{2\a-i} \prod\limits_{i\leq j}
da_{i,j} \; , \ee where
$$
g_\a (a'a)= \frac{2^{m-k}}{ \Gamma_{m-k}(\a-k/2)}\intl_{ B}
f(a'a+b'b) \prod\limits_{i=1}^{m-k} b_{k+i,k+i}^{2(\a-k/2) -i}
\prod\limits_{i\leq j} db_{k+i,k+j}\; . $$ Note that $ g_\a
(a'a)=\G_{\a-k/2}\left( f \left( \left[ \begin{array} {ll} \ast &  \ast \\
                             \ast & \bullet
\end{array}\right]\right)\right)
$ is the G{\aa}rding-Gindikin  distribution acting on the
$(\bullet)$ matrix variable belonging to $\S_{m-k}$. By
(\ref{dp}), $\G_\a (f)$ is a direct product of two distributions
which are analytic in $\a$. By taking into account (\ref{r0}), we
get $g_{k/2} (a'a)=f(a'a)$, and therefore, \be\label{rFK} \G_{k/2}
(f)= c \intl_{A} f (a'a)\prod\limits_{i=1}^{k} a_{i,i}^{k-i}
\prod\limits_{i\leq j} da_{i,j} \; .\ee Here, by (\ref{2.5}),
$$
c=2^k\lim\limits_{\a\to k/2}\frac{2^k
\Gamma_{m-k}(\a-k/2)}{\gm(\a)} =\frac{2^k
\pi^{k(k-m)/2}}{\Gamma_{k}(k/2)}\, .
$$
 This representation was established in \cite[p. 134]{FK}. Let
us show that it coincides with (\ref{er}). We replace $\om$ in
(\ref{er}) by $[\eta,\zeta]$, where $\eta\in\frM_{k,k}$,
$\zeta\in\frM_{k,m-k}$. Then
$$
\pi^{-km/2}\intl_{\Mkm}
f(\om'\om)d\om=\pi^{-km/2}\intl_{\frM_{k,k}}d\eta
\intl_{\frM_{k,m-k}} f \left ( \left [\begin{array}{ll}   \eta ' \eta &   \eta ' \zeta \\
 \zeta ' \eta &  \zeta ' \zeta
\end{array}\right]\right ) d\zeta=
$$
\begin{center}
( set  $\eta=\g q$,\; $\g\in O(k)$,\; $q\in T_k$ and use Lemma
\ref{sph})
\end{center}
\bea\nonumber &=&\pi^{-km/2}\sig_{k,k} \intl_{T_k}
\prod\limits_{i=1}^k q_{i,i}^{k-i} \prod\limits_{i\leq j} dt_{i,j}
\intl_{\frM_{k,m-k}} f \left ( \left [\begin{array}{ll}    q'q &   q ' y \\
 y 'q &  y' y
\end{array}\right]\right ) dy\\\nonumber
&\stackrel{\rm (\ref{2.16})}{=}& \frac{2^k
\pi^{k(k-m)/2}}{\Gamma_{k}(k/2)}\intl_{A} f
(a'a)\prod\limits_{i=1}^{k} a_{i,i}^{k-i} \prod\limits_{i\leq j}
da_{i,j} \nonumber \eea
where $a=\left [\begin{array}{ll}    q &    y \\
 0 &  0
\end{array}\right].$ This proves the statement.

\end{proof}

\section{Fractional integrals of locally integrable
functions}\label{s4}

In the previous sections, we studied fractional integrals
$I_{\pm}^{\a}f$ assuming  $f$ is very nice, namely, $f \in
\D(\p)$. Below we  explore these integrals for arbitrary locally
integrable functions $f$  provided that $\a$ belongs to the
Wallach-like set \index{b}{Latin and Gothic!wa@$\W$}
\index{a}{Wallach-like set, $\W$}\be\label{wal} \W =\left\{0,
\frac{1}{2}, 1, \frac{3}{2}, \ldots , \frac{m-1}{2} \right\} \cup
\left\{\a:\; Re\,\a> \frac{m-1}{2} \right \}.\ee The key question
is for which $f$ the integrals $I_{\pm}^{\a}f$ do exist. We start
with the left-sided integrals.
\begin{definition}\label{defI+} \index{b}{Latin
and Gothic!ia@$I^\a_{\pm}f$} For a locally integrable function $f$
on $\p$ and $\a \in \W$, we define \index{a}{G{\aa}rding-Gindikin
fractional integrals!left-sided, $I_{+}^\a f$}

\be\label{4.1a} (I_{+}^\a f)(s) = \left \{
\begin{array} {ll}\displaystyle{\frac {1}{\gma} \intl_0^s f(r)|s-r|^{\a-d} dr}
 & \mbox{ if   $Re \, \a > d-1$}, \\
{} \\
 \pi^{-km/2}\!\!\!\!\!\!\!\!\!\!\!\!\!\displaystyle{\intl_{\{\om\in\Mkm:\;\om'\om<s\}}\!\!\!\!\!\!\!\!\!\!\!\!
f(s-\om'\om)d\om}& \mbox { if   $\a=k/2$}. \\
\end{array}
\right.
 \ee
Here, $ s \in \p$, $d=(m+1)/2$, and $k=1,2, \ldots , m-1$. For
$\a=0$, we set $(I_{+}^0 f)(s)=f(s)$.
\end{definition}

This definition agrees with consideration in  Sections \ref{s3.2}
and \ref{s2.2}.
 We recall that the second line in (\ref{4.1a}) coincides with the
 first one if $k \ge m, \; \a=k/2$.
\begin{theorem}\label{l3.4} If $f \in L^1_{ loc}(\p)$ and $ \a \in \W$, then $(I_{+}^\a
f)(s)$ converges  absolutely  for almost all $s \in \p$.
\end{theorem}
\begin{proof}  It suffices to show that the
integral
 $ \int_0^{a}(I_{+}^\a f)(s) ds$ is finite for any  $a \in \p$ provided that $\a\in\W$  and
  $f \in L^1_{ loc}(\p)$ is  nonnegative. For $\a > d-1$,  changing the order of integration,
  and evaluating the inner integral
  according to (\ref{2.7}), we have
\bea\nonumber  \intl_0^{a}(I_{+}^\a f)(s) ds  &=&\frac {1}{\gma}
\intl_0^{a}f(r)dr\intl_r^a |s -r|^{\a-d} ds \\\nonumber
&=&\frac{\Gam_m (d)}{\Gam_m (\a +d)}\intl_0^a f(r)|a-r|^\a
dr\\\nonumber &\leq& \frac{\Gam_m (d)|a|^\a}{\Gam_m (\a
+d)}\intl_0^a f(r) dr<\infty \eea (note that $|a-r|<|a|$; see,
e.g., \cite[p. 586]{Mu} ). For the second line of (\ref{4.1a}), by
changing the order of integration, we obtain
  \bea \int_0^{a} (I_{+}^{k/2}f)(s)ds &=&\pi^{-km/2}
  \intl_{\om'\om<a}d\om\intl_{\om'\om}^a  f(s-\om'\om)ds
  \nonumber \\ &=&\pi^{-km/2}
  \intl_{\om'\om<a}d\om\intl_0^{a-\om'\om} f(r) dr  \nonumber \\
  &\le&\pi^{-km/2}c_a \intl_0^{a} f(r) dr , \nonumber \eea
where \[ c_a=\intl_{\om'\om<a}d\om\stackrel{\rm
(\ref{2.15.1})}{=}\frac{\pi^{km/2}\;\gm(d)}{\gm(k/2+d)}\;
|a|^{k/2} <\infty,
\]  as
required.
\end{proof}

The following equality can be easily obtained by changing the
order of integration and using (\ref{2.14})--(\ref{2.15}):

\be \label{2.26}\intl_{\p}\frac{(I_{+}^\a f)(s)\, ds}{|\e
I_m+s|^{\g}}=\frac{\gm(\g-\a)}{\gm(\g)} \intl_{\p}\frac{ f(r)\,
dr}{|\e I_m+r|^{\g-\a}}\; ,\ee
$$
\a\in\W,\qquad Re\,(\g-\a)>(m-1)/2,\qquad \e=0,1.
$$

Let us consider the right-sided integrals $I_{-}^{\a}f$. An idea
of the following definition is the same as above.
\begin{definition} For a locally integrable function $f$ on $\p$ and $ \a \in \W$, we
define  fractional integrals $(I_{-}^{\a}f)(s), \; s \in
\overline{\p},$ by \index{a}{G{\aa}rding-Gindikin fractional
integrals!right-sided, $I_{-}^\a f$} \index{b}{Latin and
Gothic!ia@$I^\a_{\pm}f$}

\be\label{4.1b} (I_{-}^\a f)(s) = \left \{
\begin{array} {ll}\displaystyle{\frac {1}{\gma} \intl_s^\infty f(r)|r-s|^{\a-d} dr}
 & \mbox{ if   $Re \, \a > d-1$}, \\
{} \\
 \pi^{-km/2}\displaystyle{\intl_{\Mkm}
f(s+\om'\om)d\om}& \mbox { if   $\a=k/2$.}
\end{array}
\right.
 \ee
 For $\a=0$, we set $(I_{-}^0 f)(s)=f(s)$.
As before, here we assume $d=(m+1)/2$ and $k=1,2, \ldots , m-1$.
\end{definition}

Note that unlike  Definition \ref{defI+}, here  $s$ may be a
boundary point of $\p$. The following statement extends Lemma
\ref{l3.10} (on interrelation between $I_{+}^\a$ and $I_{-}^\a$)
to locally integrable functions.

\begin{theorem}\label{l3.14}
Let $\a\in\W$, $d=(m+1)/2$. Suppose that $f$ is a measurable
function on $\p$ satisfying
 \be\label{3.34} \intl_{R}^\infty
|r|^{Re\,\a-d}|f(r)| \, dr<\infty ,\quad  \forall R \in \p,
 \ee
 or, equivalently, \be\label{3.34.1}
g(r)\equiv |r|^{-\a-d}f(r^{-1})\in L^1_{loc}(\p).\ee Then the
integral $(I_{-}^\a f)(s)$ defined by (\ref{4.1b}) converges
absolutely for almost all $s \in \p$ and
\be\label{3.7a}(I_{-}^{\a}f)(s)=|s|^{\a-d}(I_{+}^{\a}g)(s^{-1}).
 \ee

\end{theorem}
\begin{proof}
Equivalence of    (\ref{3.34}) and (\ref{3.34.1}) is obvious in
view of Lemma \ref{12.2} (iii). In fact, we have already proved
(\ref{3.7a}) provided that either side of it exists in the
Lebesgue sense. This was done in Lemma \ref{l3.10} for
$Re\,\a>d-1$, and in Theorem \ref{l3.11} for $\a=k/2$, $0<k<m$. By
Theorem \ref{l3.4}, the condition $g(r)\equiv
|r|^{-\a-d}f(r^{-1})\in L^1_{loc}(\p)$ guarantees finiteness of
the right-hand side of (\ref{3.7a}) for almost all $s$. This
completes the proof.

\end{proof}

\begin{remark}\label{r3.15} Condition (\ref{3.34}) is  best possible in the sense that there is a function
$f\geq 0$ for which (\ref{3.34}) fails and $(I_{-}^\a f)(s)\equiv
\infty$. Let \be\label{3.38} f_\lam (r)= |I_m +r|^{-\lam/2}, \ee
where $\lam \le 2\a +m-1$, $\a$ is real. Then for  $\a
>  d-1$,
\[ \intl_{R}^\infty |r|^{\a-d} f_\lam (r) \, dr=\intl_0^{R^{-1}}|r|^{\lam/2 -\a -d} |I_m +r|^{-\lam/2}\,
dr\equiv\infty,\] and  by (\ref{2.13}), \bea \label{3.37}(I_{-}^\a
f_\lam)(s) &=&\frac{|I_m +s|^{\a-\lam/2}}{\gma}
\intl_0^{I_m}|r|^{\lam/2 -\a -d} |I_m-r|^{\a-d}\, dr
\\ &=&|I_m +s|^{\a-\lam/2}B_m(\a,\lam/2 -\a)/\gma\equiv\infty. \nonumber \eea
 In the case  $\a=k/2$, $k<m$, we have \bea
\label{3.39}(I_{-}^{k/2}f_\lam )(s)&=&\pi^{-km/2}\intl_{\Mkm}
f_\lam (s+\om'\om)d\om \nonumber \\ &=&\pi^{-km/2}\intl_{\Mmk} |b+
yy' |^{-\lam/2}\, dy, \qquad b=I_m +s. \eea Owing to (\ref{2.15}),
the last integral coincides with (\ref{3.37}) (with $\a=k/2$) and
is infinite.
\end{remark}

The following formula is a consequence of (\ref{2.26}) and
(\ref{3.7a}):

\be \intl_{\p}(I_{-}^\a f)(s)\;\frac{|s|^{\g-a-d}}{| I_m+\e
s|^{\g}}\; ds=\frac{\gm(\g-\a)}{\gm(\g)} \intl_{\p}
f(r)\;\frac{|r|^{\g-d}}{| I_m+\e r|^{\g-\a}}\; dr\; ,\ee
$$
\a\in\W,\qquad d=(m+1)/2,\qquad Re\,(\g-\a)>(m-1)/2,\qquad \e=0,1.
$$

\section{The semigroup property} By the Fubini theorem, the equality
(\ref{2.7}) yields
 \be\label{3.14}
I_{\pm}^{\a}I_{\pm}^{\b}f=I_{\pm}^{\a+\b}f,\qquad Re\,\a
>d-1,\quad Re\,\b >d-1 ,\ee  provided
the corresponding integral $I_{+}^{\a+\b}f$ or $I_{-}^{\a+\b}f$ is
absolutely convergent. The equality (\ref{3.14}) is usually
referred to as the semigroup property of fractional integrals. It
extends to all complex $\a$ and $\b$ if $f$ is good enough  and
fractional integrals are interpreted as multiplier operators in
the framework of the relevant Fourier-Laplace analysis. This
approach allows us to formulate the semigroup property in the
language of the theory of distributions; see \cite{Gi1},
\cite{FK}, \cite{VG}, \cite{Rab} for details.

Below we extend (\ref{3.14})
 to the important special case when $\a$ and $\b$ lie in
 the Wallach-like set $\W$ defined by (\ref{wal}). Definitions (\ref{4.1a}) and (\ref{4.1b}) enable us
to do this under  minimal assumptions for $f$ so that all
integrals  are understood in the classical Lebesgue sense.

\begin{lemma} If  $f\in L^1_{loc}(\p)$,
  then \index{a}{G{\aa}rding-Gindikin fractional integrals!semigroup property}
  \be\label{3.14.1}
I_{\pm}^{\a}I_{\pm}^{\b}f=I_{\pm}^{\a+\b}f,\qquad \a, \b\in\W ,\ee
 provided  the  corresponding integral on the right-hand
side is absolutely convergent.
\end{lemma}

\begin{proof}
In view of (\ref{3.14}), it suffices to prove that \be\label{3.15}
I_{\pm}^{\a}I_{\pm}^{k/2}f=
I_{\pm}^{k/2}I_{\pm}^{\a}f=I_{\pm}^{\a+k/2}f,\qquad Re\,\a
>d-1, \ee and
\be\label{3.16}
I_{\pm}^{k/2}I_{\pm}^{\ell/2}f=I_{\pm}^{(k+\ell)/2}f, \ee where
$k$ and $\ell$ are positive integers less than $m$.
   By changing the order of integration, we
have \bea\label{3.18} (I_{-}^{\a}I_{-}^{k/2}f)(s)&=&
\frac{\pi^{-km/2}}{\Gamma_m(\a)}\intl_s^\infty |r-s|^{\a-d} \, dr\intl_{\Mkm} f(r+\om '\om) \, d\om \nonumber \\
&=&\frac{\pi^{-km/2}}{\Gamma_m(\a)}\intl_{\Mkm}
 d\om\intl_{s+\om'\om}^\infty f(t) |t-s-\om'\om|^{\a-d} \, dt \nonumber \\
&=&\frac{\pi^{-km/2}}{\Gamma_m(\a)}\intl_s^\infty f(t)h(t-s) \,
dt,\nonumber \eea where
$$
h(z)=\intl_{\{\om\in\Mkm:\;\om '\om<z\}} |z-\om '\om|^{\a-d}d\om.
$$
The same expression comes out if we transform
$I_{-}^{k/2}I_{-}^{\a}f$. By (\ref{2.15}), \be\label{3.19}
h(z)=\frac{\pi^{km/2}\gm(\a)}{\gm(\a+k/2)} |z|^{\a+k/2-d}, \qquad
Re\,\a>d-1,\ee and we are done. For the left-sided integrals, the
argument is  similar. Namely, \bea (I_{+}^{\a}I_{+}^{k/2}f)(s)&=&
\frac{\pi^{-km/2}}{\Gamma_m(\a)}\intl_0^s
|s-r|^{\a-d}dr\intl_{\{\om\in\Mkm\,:\, \om'\om<r\}} f(r-\om
'\om)d\om \nonumber \\ &=&
\frac{\pi^{-km/2}}{\Gamma_m(\a)}\intl_0^s f(t)h(s-t) \,
dt\nonumber \\  &=&(I_{+}^{\a+k/2}f)(s). \nonumber \eea

 Let us
prove (\ref{3.16}).  By (\ref{4.1b}),
$$
(I_{-}^{k/2}I_{-}^{\ell/2}f)(s)=\pi^{-m(k+\ell)/2}\intl_{\Mkm}d\om\
\intl_{\Mlm}f(y'y+\om'\om+s)dy.
$$
The change of variables $z=\left[\begin{array}{l} y\\ \om
\end{array}\right] \in \frM_{k+\ell,m}$ yields $z'z=y'y+\om'\om$,
and therefore
$$
(I_{-}^{k/2}I_{-}^{\ell/2}f)(s)=\pi^{-m(k+\ell)/2}\intl_{\frM_{k+\ell,m}}
f(z'z+s)dz=(I_{-}^{(k+\ell)/2}f)(s).
$$
For the left-sided integrals, owing to (\ref{4.1a}), we have
\[
(I_{+}^{(k+\ell)/2}f)(s)=\pi^{-m(k+\ell)/2}\intl_{\{z\in\frM_{k+\ell,m}\,:\,z'z<s\}}
f(s-z'z) \, dz. \] By the Fubini theorem, this reads
\[\pi^{-m(k+\ell)/2}\!\!\!\!\!\intl_{\{\om\in\Mkm\,:\,\om'\om<s\}}\!\!\!\!\!\!
d\om\
\intl_{\{y\in\Mlm:\,y'y<s-\om'\om\}}\!\!\!\!\!\!f(s-y'y-\om'\om)
\,  dy\]
 and coincides with $ (I_{+}^{k/2}I_{+}^{\ell/2}f)(s)$.
\end{proof}

\section{Inversion formulas} Below we obtain  inversion formulas for the integrals
$(I_{\pm}^{\a}f)(r), \; r \in \p$.  For our purposes, it suffices
to restrict to the case $\a\in\W$. We recall that $k \in \bbn$ and
$d=(m+1)/2$. For ``good" functions $f$, the following statement
holds.
\begin{lemma} If $\vp=I_{\pm}^{\a}f, \;  f\in \D(\p), \; \a
\in\W$, then \be\label{inv} f =D^j _{\pm}I_{\pm}^{j-\a} \vp, \ee
provided that $j\geq\a$ if $\a=k/2$ and $j>Re\,\a+d-1$ otherwise.
\end{lemma}
\begin{proof} By  the semigroup property (\ref{3.14.1}), we have
$I_{\pm}^{j-\a} \vp=I_{\pm}^j f$.  Then  application of
(\ref{3.49}) and
 (\ref{2.32}) (or (\ref{2.33})) gives (\ref{inv}).
\end{proof}

Unfortunately,   we are not able to obtain pointwise inversion
formulas for $I_{\pm}^{\a}f$ if $f$ is an arbitrary locally
integrable function. To get around this difficulty, we utilize the
theory of distributions.

\begin{lemma} \label{l-inv1} \index{a}{G{\aa}rding-Gindikin fractional
integrals!left-sided, $I_{+}^\a f$!inversion} If $\vp=I_{+}^{\a}f,
\; f\in L^1_{loc}(\p), \; \a \in\W$, then \be\label{inp1} f =D^j
_{+}I_{+}^{j-\a} \vp, \ee for any integer $j$ subject to the
condition that  $j\geq\a$ if $\a=k/2$ and $j>Re\,\a+d-1$
otherwise. The differentiation in (\ref{inp1}) is understood in
the sense of distributions, so that
 \be\label{inp2} (f,\;\phi)=(I_+^{j-\a}\vp,\;
{D_-^j\phi}),\qquad \phi\in\D (\p). \ee In particular, for $\a=j$,
 \be\label{inp3} (f,\;\phi)=
 (\vp,\; {D_{-}^j \phi}).
\ee
\end{lemma}
\begin{proof} As above, we have
$I_{+}^{j-\a} \vp=I_{+}^j f$, and therefore,
\[ (f,\;\phi)\stackrel{\rm (1)}{=}(f, \; I_-^{j} D_-^j  \phi)\stackrel{\rm (2)}{=}
 (I_{+}^{j}f,\;  D_{-}^j \phi)\stackrel{\rm (3)}{=}(I_+^{j-\a}\vp,\; D_-^j\phi).\]
The equality (1) holds thanks to (\ref{3.49}), in (2) we apply the
Fubini theorem; (3) holds owing to the semigroup
 property (\ref{3.15}).
\end{proof}

Let us invert the equation  $I_{-}^{\a} f= \vp$. Formally,
$f=D^\a_- \vp$ where $D^\a_- =I^{-\a}_-$, or  $(f,\;\phi)=(\vp, \;
I^{-\a}_+ \phi)$, $\phi\in\D (\p)$ in the sense of distributions.
The next statement justifies this equality.
\begin{lemma}\label{ligm} \index{a}{G{\aa}rding-Gindikin fractional
integrals!right-sided, $I_{-}^\a f$!inversion} Let $\a\in\W$, and
let $f$ be a locally integrable function on $\p$ subject to the
decay condition (\ref{3.34}), so that $(I_{-}^{\a}f)(s)$ exists
for almost all $s \in \p$. If $\vp=I_{-}^{\a}f$ then the integral
\[ (\vp, \; I^{-\a}_+ \phi)=\intl_{\p} \vp (r) (I^{-\a}_+ \phi)(r)
dr,\qquad  \phi\in\D (\p),
\] is finite, and $ f$ can be recovered
by the formula \be\label{inm} f= D_-^\a \vp, \ee where  $D_-^\a
\vp$ is defined as a distribution \be\label{dmin} (D_-^\a \vp,
\;\phi)= (\vp, \; I^{-\a}_+ \phi). \ee
 In particular, for  $\a=j$,  \be\label{inp3}
(f,\;\phi)=
 (\vp,\; {D_{+}^j \phi}).
\ee
\end{lemma}
\begin{proof} Let
$ \phi_1(s)=|s|^{\a-d}\phi(s^{-1})$.  This function belongs to
$\D(\p)$. By (\ref{3.7}),
$$(I_{+}^{-\a} \phi)(s)=|s|^{-\a-d}(I_-^{-\a} \phi_1)(s^{-1}).$$
If $\vp=I_{-}^{\a}f$ then \bea\nonumber (\vp, \; I^{-\a}_+
\phi)&=& ((I_{-}^{\a}f)(s), \; |s|^{-\a-d}(I_-^{-\a}
\phi_1)(s^{-1})) \nonumber \\ && \hskip -1truecm \text {\rm (use
Theorem \ref{l3.14}   with $g(r)= |r|^{-\a-d}f(r^{-1})$)}\nonumber
\\ &=&((I_{+}^{\a}g)(s^{-1}), \; |s|^{-2d}(I_-^{-\a}
\phi_1)(s^{-1})) \nonumber \\ &=&
((I_{+}^{\a}g)(s), \; (I_-^{-\a} \phi_1)(s)) \nonumber \\
&\stackrel{\rm (1)}{=}& (g, \; I_-^{\a}I_-^{-\a} \phi_1)\nonumber
\\ &\stackrel{\rm (2)}{=}&  (g, \;   \phi_1) \nonumber \\ &=& (f,
\phi).
 \nonumber \eea

Let us comment on these calculations. The equality (1) holds by
the Fubini theorem. Here we take into account that $I_-^{-\a}
\phi_1 \in \D (\p)$, because by the second equality in
(\ref{3.49}), \be\label{big} I_-^{-\a} \phi_1= I_{-}^{j-\a}D_{-}^j
\phi_1 \ee
 if $j$ is big enough. A rigorous justification of the obvious
 equality $I_-^{\a}I_-^{-\a} \phi_1=\phi_1$ used in (2) is as
 follows:
 \[ I_-^{\a}I_-^{-\a} \phi_1 \stackrel{\rm (\ref{big})}{=} I_-^{\a}I_-^{j-\a}D_{-}^j
 \phi_1 \stackrel{\rm (\ref{3.14.1})}{=}I_-^{j}D_{-}^j
 \phi_1 \stackrel{\rm (\ref{3.49})}{=}I_-^{0}
 \phi_1 \stackrel{\rm (\ref{2.33})}{=}\phi_1.\]
For $\a=j$, the result follows by Corollary  \ref{dif}, according
to which $ I_{+}^{-j} \phi=D_{+}^{j} \phi $. The proof is
complete.
\end{proof}

\chapter{Riesz potentials on matrix spaces}

\setcounter{equation}{0}

  Riesz potentials of
functions of matrix argument and their generalizations  arise in
different contexts in harmonic analysis, integral geometry, and
PDE; see \cite{Far}, \cite{Ge}, \cite{Kh},
 \cite{St2}, \cite{Sh1}--\cite{Sh3}. They are intimately
connected with representations of Jordan algebras and  zeta
functions (or zeta distributions) extensively studied in the last
three decades; see \cite{FK}, \cite{Cl1}, \cite{Ig}, \cite{SS},
\cite{Sh}  and references therein. In this chapter, we explore
basic properties  of matrix Riesz potentials which will be used
 in Chapters \ref{s5} and \ref{s6} in our study of Radon transforms.
Numerous results, bibliography, and historical notes related to
ordinary Riesz potentials on $\rn$ can be found in \cite{Ru1},
\cite{Sa}, \cite{SKM}.

\section[Zeta integrals]{Zeta integrals}

\subsection{Definition and examples}
Let $f(x)$ be a Schwartz function on the matrix space $\Ma$,
$n\geq m$. We denote \index{b}{Latin and
Gothic!xm@\texttt{"|}$x$\texttt{"|}$_m$} \be\label{xm}|x|_m =\det
(x'x)^{1/2}.\ee For $m=1$, this is the usual euclidean norm on
$\bbr^n$. If $m>1$, then  $|x|_m$ is the volume of the
parallelepiped spanned  by the column-vectors of the matrix $x$
\cite[ p. 251]{G}.
 Let us consider the following distribution \index{b}{Latin and Gothic!zf@$\Z(f,\a-n)$} \be\label{zeta}
\Z(f,\a-n)=\intl_{\Ma} f(x) |x|^{\a-n}_m dx , \qquad \a\in \bbc .
\ee For $Re\, \a >m-1$, the integral (\ref{zeta})  absolutely
converges, and for $Re\, \a \leq m-1$, it is understood in the
sense of analytic continuation (see Lemma \ref{lacz}). Integrals
like
 (\ref{zeta}) are also known as the {\it zeta integrals} or  {\it zeta
distributions}.

The case $n=m$ deserves special consideration. In this case, we
denote $\frM_m=\frM_{m,m}$ and define two zeta integrals
 \index{a}{Zeta integrals} \index{b}{Latin and
Gothic!zf@$\Z_{\pm}(f,\a-n)$}\bea\label{plus}
\Z_+(f,\a-m)&=&\intl_{\frM_m} f(x) \,|\det
(x)|^{\a-m} \,  dx \qquad (\equiv \Z(f,\a-m)), \\
\label{minus} \Z_-(f,\a-m)&=&\intl_{\frM_m} f(x) \, |\det
(x)|^{\a-m} \, \sdx \, dx. \eea We call (\ref{minus}) {\it the
conjugate zeta integral} (or distribution) by analogy with the
case $\a=0, \, m=1$, where
\[ \Z_-(f, -1)= \text{\rm p.v.} \intl_{-\infty}^{\infty}
\frac{f(x)}{x} \, dx.\] A convolution of $f$ with the distribution
p.v.$\frac{1}{x}$ represents the Hilbert transform \cite{Ne},
\cite{SW}, and is called a conjugate function.

\begin{example} It is instructive to evaluate zeta integrals for the Gaussian
functions. Let \be\label{ga} e(x)=\exp(-\tr(x'x)).\ee By Lemma
\ref{l2.3}, for $Re\, \a >m-1$ we have \be\label{epl} \Z(e,
\a-n)=2^{-m}\sig_{n,m}\intl_{\p}|r|^{\a/2 -d}
\exp(-\tr(r))dr=c_{n,m}\gm(\a/2),\ee where $\sig_{n,m}$ is the
constant (\ref{2.16}), $d=(m+1)/2$,
\[
c_{n,m}=\frac{\pi^{nm/2}}{\gm(n/2)}\;,  \quad n\ge m.\] In the
case $n=m$, for the function
 \be\label{ga1} e_1(x)= e(x)\,\det (x)\ee
we have \be\label{min1}\Z_-(e_1,\a-m)=\Z_+(e, \a+1-m)=c_m \Gam_m
((\a +1)/2),\ee
\[c_m \equiv c_{m,m}=\frac{\pi^{m^2/2}}{\gm(m/2)}, \qquad Re \, \a
>m-2.\]
Formulas (\ref{epl}) and (\ref{min1}) extend meromorphically to
all complex $\a$, excluding  $\a=m-1, m-2, \dots$ for (\ref{epl}),
and $\a=m-2, m -3, \dots$ for (\ref{min1}). Excluded values are
poles of the corresponding gamma functions written  for $m\geq 2$.
If $m=1$, these poles proceed with step $2$, namely, $\a= 0,-2,-4,
\dots\; $ and $\a=-1, -3, \dots$, respectively.
\end{example}

\subsection{Analytic continuation}

In the following, throughout the paper, we assume $m\geq 2$.

\begin{lemma}\label{lacz} \index{a}{Zeta integrals!analytic continuation}
Let $f\in\s(\Ma)$. For $Re\, \a > m-1$,  the integrals
(\ref{zeta})-(\ref{minus}) are  absolutely convergent. For $Re\,
\a \leq m-1$, they extend as  meromorphic functions of $ \; \a$
with the only poles  \bea &{}&m-1, \; m-2,\dots , \qquad \text{\rm
for (\ref{zeta}) and (\ref{plus})}, \nonumber \\
&{}&  m-2, \; m-3, \dots , \qquad \text{\rm for
 (\ref{minus})}.\nonumber \eea  These poles and their orders
are exactly the same as those of the gamma functions $\gm(\a/2)$
and $\Gam_m ((\a +1)/2)$, respectively. The normalized  zeta
integrals \be\label{nor} \frac{\Z(f,\a-n)}{\gm(\a/2)}, \qquad
\frac{\Z_+(f,\a-m)}{\gm(\a/2)}, \qquad \frac{\Z_-(f,\a-m)}{\Gam_m
((\a +1)/2)}\ee are
 entire
functions of $\a$.
\end{lemma}
\begin{proof}
This statement is known; see, e.g., \cite{P5}, \cite{Kh},
\cite{Sh1}. We present the proof for the sake of completeness. The
equalities (\ref{epl})  and (\ref{min1}) say
 that the functions $\a \to \Z(f,\a-n)$ and  $\a \to \Z_-(f,\a-m)$
 have
poles at least at the same points and   of the same order as the
gamma functions $\gm(\a/2)$ and $\Gam_m ((\a +1)/2)$,
respectively. Our aim is to show that no other poles occur, and
the orders cannot exceed those of
 $\gm(\a/2)$ and $\Gam_m ((\a +1)/2)$.
Let us transform (\ref{zeta}) by passing to upper triangular
matrices $t\in T_m$ according to Lemma \ref{sph}. We have
\be\label{z-tr} \Z(f,\a-n)=\intl_{\bbr^m_+} F(t_{1,1},\dots ,
t_{m,m}) \prod\limits_{i=1}^m t_{i,i}^{\a-i} dt_{i,i}\; , \ee
$$
F(t_{1,1},\dots , t_{m,m}) =\!\!\!\!\!\! \intl_{\bbr^{m(m-1)/2}}
\!\!\!\! dt_{*} \intl_{\vnm} f(vt) dv, \quad
dt_{*}=\prod\limits_{i<j} d t_{i,j}.
$$
Since $F$ extends as an even Schwartz function in each argument,
it can be written as
$$
F(t_{1,1},\dots , t_{m,m}) =F_0(t^2_{1,1},\dots , t^2_{m,m}),
$$
where $F_0\in\s(\bbr^m)$ (use, e.g., Lemma 5.4 from \cite[p.
56]{Tr}). Replacing $t_{i,i}^2$ by $s_{i,i}$, we represent
(\ref{z-tr}) as a direct product of one-dimensional distributions
\be\label{z-h} \Z(f,\a-n)=2^{-m}(\prod\limits_{i=1}^m
(s_{i,i})_+^{(\a-i-1)/2},\;  F_0(s_{1,1},\dots , s_{m,m})),\ee
which is a meromorphic function of $\a$ with the poles $m-1, m-2,
\dots \;$, see \cite{GSh1}. These poles and their orders coincide
with those  of the gamma function $\gm(\a/2)$. To normalize the
function (\ref{z-h}), following \cite{GSh1}, we divide it by the
product
$$\prod\limits_{i=1}^m\Gam((\a-i+1)/2)=\prod\limits_{i=0}^{m-1}\Gam((\a-i)/2)=
\gm(\a/2)/\pi^{m(m-1)/4}.$$ As a result, we obtain an entire
function.

For the distribution $\Z_-(f,\a-m)$, the argument is almost the
same. Namely, by Lemma \ref{sph}, the integral (\ref{minus}) can
be written  as
$$
 \intl_{\bbr^m_+}
 \Phi (t_{1,1},\dots , t_{m,m})  \,\prod\limits_{i=1}^m t_{i,i}^{\a -i} \, dt_{i,i}
$$ where \be\label{ffi} \Phi (t_{1,1},\dots ,
t_{m,m})=\sig_{m,m}\!\!\!\!\intl_{\bbr^{m(m-1)/2}} \, dt_{*}
\intl_{O(m)} \, f(vt) \,  \sdv
 \, dv.\ee Since $\Phi$ extends to $\bbr^m$ as an odd Schwartz
function in each argument, then \be\label{ffii} \Z_-(f,\a-m)
=\frac{1}{2^m}\intl_{\bbr^m} \Phi (t_{1,1},\dots , t_{m,m})
\prod\limits_{i=1}^{m} |t_{i,i}|^{\a -i} \, \sgn (t_{i,i}) \,
dt_{i,i}.\ee By the well known theory \cite[Chapter 1, Section
3.5]{GSh1}, this integral extends as a meromorphic function of
$\a$ with the only poles $m-2, m-3,\dots\;$. To normalize this
function, we divide it by the product
$$\prod\limits_{i=1}^m\Gam \left(\frac{\a-i+2}{2}\right)=\prod\limits_{i=0}^{m-1}
\Gam \left(\frac{\a-i+1}{2}\right)= \Gam_m \left(\frac{\a
+1}{2}\right)/\pi^{m(m-1)/4}.$$ As a result, we get an entire
function.
\end{proof}

Analytic continuation of integrals (\ref{zeta})-(\ref{minus}) can
be performed with the aid of the corresponding differential
operators. In particular, for $\Z_{\pm}(f, \a-m)$, one can utilize
{\it the Cayley differential operator} \index{b}{Latin and
Gothic!d@$\DP$}\index{a}{Cayley operator!on $\Ma$, $\DP$}
\be\label{cay} \DP=\det(\d/\d x_{i,j}), \quad x=(x_{i,j}) \in
\frM_m , \ee that enjoys the following  obvious relations in the
Fourier terms: \be\label{fu1}
 (\F[\DP f]) (y)= (-i)^m \, \det (y) \, (\F f)(y), \ee
\be\label{fu2}
 \DP (\F f) (y)= i^m \, (\F [ f(x) \, \det (x)])(y). \ee
\begin{lemma} Let $x\in \frM_m, \; \rank (x)=m$. For any $\lam
\in \bbc$, \be\label{d1} \DP \, [|\det (x)|^\lam ]=(\lam,m) \,
|\det (x)|^{\lam -1} \, \sdx, \ee \be\label{ds} \DP \, [|\det
(x)|^\lam \, \sdx] =(\lam,m) \, |\det (x)|^{\lam -1}. \ee
\end{lemma}
\begin{proof}
Different proofs of these important formulas can be found in
\cite{Ra} and \cite{P2}; see also \cite[p. 114]{Tu}. All these
proofs are very involved. Below we give an alternative proof which
is elementary. Note that (\ref{d1}) and (\ref{ds}) follow one from
another. Furthermore, it suffices to assume that $\lam$ is not an
integer ($\lam \notin \bbz$). Once (\ref{d1}) and (\ref{ds}) are
proved for such $\lam$, the result for $\lam \in \bbz$ then
follows by continuity.

We start with the formula \be \label{svy} \DP_x [f(ax)]=\det (a)
\, (\DP f)(ax), \qquad a \in GL(m, \bbr), \ee which can be easily
checked by applying the Fourier transform to both sides. Indeed,
if $f$ is good enough at infinity (otherwise we can multiply $f$
by a smooth cut-off function) then  by (\ref{fu1}), the Fourier
transform of the left-hand side of (\ref{svy}) reads
\[ (-i)^m \,\det (y) \, \F[f(ax)](y)=
\frac{(-i)^m \,\det (y)}{|\det (a)|^m}(\F f)((a^{-1})' y)\] which
coincides with  the Fourier transform of the right-hand side. If
$f(x)=|\det (x)|^\lam $ then (\ref{svy}) yields \[ |\det (a)|^\lam
\, \DP \, |\det (x)|^\lam = \det (a) \, [\DP \, |\det
(\cdot)|^\lam](ax).\] By setting $a=x^{-1}$ (recall that $\rank
(x)=m$ so that $x$ is non-singular) we obtain \[  \DP \,  |\det
(x)|^\lam = A |\det (x)|^{\lam -1} \, \sdx, \quad A=[\DP \, |\det
(x)|^\lam](I_m),\] and therefore \be\label{dss} \DP \, [|\det
(x)|^\lam \, \sdx] =A|\det (x)|^{\lam -1}. \ee In order to
evaluate $A$, we make use of the Gaussian functions
$$
 e(x)=\exp(-\tr(x'x)) \qquad \text{\rm and} \quad  e_1(x)= e(x)\,\det
 (x)=(-2)^{-m}(\DP e)(x),$$
see (\ref{ga}) and (\ref{ga1}). By (\ref{min1}),
\[\Z_- (e_1, \lam)=\Z_+ (e, \lam +1)=c_m \Gam_m ((\lam +1+m)/2)\]
(here and on we do not care about the poles because we assumed in
advance that $\lam \notin \bbz$). On the other hand, by
(\ref{dss}) and (\ref{epl}), \bea \Z_- (e_1,
\lam)&=&(-2)^{-m}\Z_-(\DP
e,\lam)\nonumber \\
&=&2^{-m}(e(x), \DP \, [|\det (x)|^\lam  \, \sdx]) \nonumber \\
&=&2^{-m} A \Z_+ (e, \lam -1)\nonumber \\
&=&c_m 2^{-m} A \,\Gam_m ((\lam -1+m)/2). \nonumber \eea Hence,
owing to (\ref{Poh}), we obtain \[ A=\frac{2^m \,\Gam_m ((\lam
+1+m)/2) }{\Gam_m (\lam -1+m)/2)}= (-1)^m \frac{\Gam
(1-\lam)}{\Gam(1-\lam -m)}=(\lam, m).\]
\end{proof}
Formulas (\ref{d1}) and  (\ref{ds}) imply the following connection
between zeta integrals $\Z_+ (\cdot)$ and $\Z_- (\cdot)$.
\begin{corollary}
\be\label{conn} \index{a}{Zeta integrals!analytic continuation}
\Z_{\mp}(f, \a -m)=c_\a  \, \Z_{\pm}(\DP f, \a +1-m), \quad
c_\a=(-1)^m \frac{\Gam (\a+1-m)}{\Gam (\a+1)}. \ee
\end{corollary}

The formula (\ref{conn}) can be used for analytic continuation of
the zeta integrals  $\Z_{\mp}(f, \a -m)$.

\subsection{Functional equations for the zeta integrals}

The following statement is a core of the theory of zeta integrals
(\ref{zeta})-(\ref{minus}).
\begin{theorem}\label{tZ-F} \index{a}{Zeta integrals! functional
equations}
 Let $f\in\s(\Ma)$, $n\geq m$. Then
 \be\label{az-f}
\frac{\Z(f, \a-n)}{\gm(\a/2)}=\pi^{-nm/2} \,  2^{m(\a - n)} \,
\frac{\Z(\F f, -\a)}{\gm((n-\a)/2)}, \ee

\be\label{az-fp} \frac{\Z_+(f, \a-m)}{\gm(\a/2)}=\pi^{-m^2/2} \,
2^{m(\a - m)} \, \frac{\Z_+(\F f, -\a)}{\gm((m-\a)/2)}, \ee

\be\label{az-fm} \frac{\Z_-(f, \a-m)}{\Gam_m ((\a +1)/2)}=(-i)^m
\, \pi^{-m^2/2} \,  2^{m(\a - m)} \, \frac{\Z_-(\F f, -\a)}{\Gam_m
((m-\a +1)/2)}. \ee
\end{theorem}
\begin{proof} We recall  that both sides of each equality  are understood in the
sense of analytic continuation and represent entire functions of
$\a$. The equalities (\ref{az-f}) and (\ref{az-fp}) were obtained
in \cite{Far}, \cite{FK}, \cite{Ge}, \cite{Ra}, \cite{P2} in the
framework of more general considerations. A self-contained proof
of them and detailed discussion can be found in \cite{Ru10}. The
equality (\ref{az-fm}) was implicitly presented in
 \cite[p. 289]{P2}. In fact, it follows
from (\ref{az-fp}) owing to the formulas (\ref{conn}) and
(\ref{fu1}). We have \bea \frac{\Z_-(f, \a-m)}{\Gam_m ((\a
+1)/2)}&=&\frac{c_\a \,\Z_+(\DP f, \a+1-m)}{\Gam_m ((\a +1)/2)}
\nonumber \\&=& c_\a \,\pi^{-m^2/2} \, 2^{m(\a+1 - m)} \, \frac{
\Z_+(\F [\DP f], - \a-1)}{\Gam_m ((m-\a -1)/2)} \nonumber
\\&=&(-i)^m \, c_\a \,\pi^{-m^2/2} \,  2^{m(\a+1 - m)} \, \frac{ ((\F
f)(y), \det (y)\, |\det (y)|^{-\a -1}) }{\Gam_m ((m-\a -1)/2)}
\nonumber
\\&=&c \;\frac{\Z_-(\F f, -\a)}{\Gam_m ((m-\a
+1)/2)}\;, \nonumber \eea where (use (\ref{conn}) and (\ref{Poh}))
\[c=(-i)^m \,
c_\a \,\pi^{-m^2/2} \,  2^{m(\a+1 - m)} \, \frac{\Gam_m ((m-\a
+1)/2)}{\Gam_m ((m-\a -1)/2)}=(-i)^m  \, \pi^{-m^2/2} \,  2^{m(\a
- m)}.\]
\end{proof}

\section{Normalized zeta distributions }

It is convenient to introduce  a special notation for the
normalized zeta integrals  (\ref{nor}) which are  entire functions
of $\a$. Let \index{b}{Greek!f@$\z_\a (x)$}
\index{b}{Greek!fa@$\z^{\pm}_\a (x)$}
 \be\label{nzd} \z_\a (x)
=\frac{|x|_m^{\a -n}}{\gm(\a/2)}, \qquad x \in \Ma, \quad n \ge m,
\ee \be \z_{\a}^+ (x) \! =\!\frac{|\det (x)|^{\a -m}}{\gm(\a/2)},
\quad \z_{\a}^- (x) \! =\!\frac{|\det (x)|^{\a -m} \, \sdx}{\Gam_m
((\a +1)/2)}, \quad x \in \frM_{m}.\ee Given a Schwartz function
$f$, we denote \be\label{nins} \index{a}{Zeta
integrals!normalized} (\z_\a, f)=a.c.\intl_{\frM_{n,m}} \!
f(x)\z_\a (x) \, dx, \quad (\z_{\a}^{\pm}, f)=a.c.\intl_{\frM_{m}}
\! f(x)\z_{\a}^{\pm} (x) \, dx \ee where``$a.c.$'' abbreviates
analytic continuation. We call $\z_\a$ and $\z_{\a}^{ \pm}$ {\it
normalized zeta distributions of order $\a$} on $\Ma$ and
$\frM_{m}$, respectively. By (\ref{conn}), \be
\label{kap1}(\z_\a^-,f)=c_\a \, (\z^+_{\a+1}, \DP f), \qquad
(\z_\a^+,f)=d_\a \, (\z^-_{\a+1}, \DP f),\ee where \bea
c_\a&=&(-1)^m \frac{\Gam (\a+1-m)}{\Gam (\a+1)} \, , \nonumber
\\d_\a &=& \frac{ c_\a \,
\gm(\a/2+1)}{\gm(\a/2)} =\frac{\Gam (\a+1-m) \, \Gam (m-\a)}{2^m
\, \Gam (\a+1) \, \Gam (-\a)}  \qquad\text{\rm(use
(\ref{Poh}))}.\nonumber \eea

 Our next goal is to obtain explicit representations of the
 distributions (\ref{nins}) when the corresponding integrals are not absolutely convergent.
Evaluation of $(\z_\a, f)$ when  $\a$ is a non-negative integer
 is  of primary importance in view of subsequent
application to the Radon transform in Chapters 5 and 6.
\begin{theorem}\label{tzk}  Let $f  \in  \s(\Ma)$. For
 $\a=k$, $k=1, 2,
\ldots,  n$, \be\label{znk} \index{a}{Zeta integrals!normalized!of
 integral order}(\z_k, f)=\frac{\pi^{(n-k)m/2}}{\gm(n/2)} \intl_{O(n)} d\gam
\intl_{\frM_{k,m}}f \left (\gam \left[\begin{array} {c} \om \\ 0
\end{array} \right]  \right ) \, d\om. \ee
Furthermore, in the case
 $\a=0$ we have \be\label{zn0} (\z_0, f)=\frac{\pi^{nm/2}}{\gm(n/2)}f(0). \ee
\end{theorem}
\begin{proof} STEP 1.
Let first $k>m-1$ . In polar coordinates we have
 \bea\nonumber
(\z_k, f)&=&\frac{1}{\gm(k/2)}\intl_{\Ma} f(x) |x|^{k-n}_m
dx\\\nonumber &=&
\frac{2^{-m}\sig_{n,m}}{\gm(k/2)}\intl_{\P_{m}}|r|^{k/2 -d} dr
\intl_{O(n)} f\left ( \gam \left[\begin{array} {c} r^{1/2} \\
0\end{array} \right]\right ) d\gam .\eea Now we replace $\gam$ by
$\gam \left[\begin{array} {cc} \b& 0 \\
0& I_{n-k}\end{array} \right]$, $\b\in O(k)$, then integrate in
$\b\in O(k)$, and replace the integration over $O(k)$ by that over
$V_{k,m}$. We get

\bea\nonumber (\z_k, f)&=&
\frac{2^{-m}\sig_{n,m}}{\sig_{k,m}\gm(k/2)}\intl_{O(n)} d\gam
\intl_{\P_{m}}|r|^{k/2 -d} dr \intl_{V_{k,m}}
 f\left ( \gam \left[\begin{array} {c} v r^{1/2} \\
0\end{array} \right]\right ) dv \\\nonumber
&{}&({\mbox set \quad  \om=v r^{1/2}\in \Mkm})\\
\nonumber &=&
\frac{\sig_{n,m}}{\sig_{k,m}\gm(k/2)}\intl_{\frM_{k,m}}d\om
\intl_{O(n)}f \left (\gam \left[\begin{array} {c} \om \\ 0
\end{array} \right]  \right ) \, d\gam.\eea
This coincides with (\ref{znk}).

 STEP 2.
Our next task is to prove that analytic continuation of $(\z_\a,
f)$ at the point $\a=k$ ($\leq m-1$) has the form (\ref{znk}). To
this end, we express $\z_\a$ through the   G{\aa}rding-Gindikin
distribution; see Section \ref{s2.3}.  For $Re\;\a>m-1$, by
passing to polar coordinates, we have $ (\z_\a,f)=
2^{-m}\sig_{n,m} \G_{\a/2}(F),$ where
$$
\G_{\a/2} (F)= \frac {1}{\gm(\a/2)} \intl_{\p} F(r)|r|^{\a/2-d}
dr, \qquad F(r)=\frac{1}{\sigma_{n,m}}\intl_{\vnm}
f(vr^{1/2})dv.$$ To continue the proof, we need the following
\begin{lemma}\label{lFar}
Let  $\s(\cpm)$ be the space of restrictions onto $\cpm$ of the
Schwartz functions on $S_m\supset \cpm$ with the induced topology.
The map
$$
\s(\cpm)\ni F\ \rightarrow f(x)=F(x'x)
$$
is an isomorphism of $\s(\cpm)$ onto the space $\s(\Ma)^\natural$
of $O(n)$ left-invariant functions on $\Ma$.
\end{lemma}
This important statement, which is well known for $m=1$ (see,
e.g., Lemma 5.4 in \cite[p. 56]{Tr}), was presented in a slightly
different form by J. Farau \cite[Prop. 3]{Far} and derived from
the more general result of G. W. Schwarz \cite[Theorem 1]{Sc}.
 According to (\ref{er}),
analytic continuation of $\G_{\a/2}(F)$ at the point $\a=k$, $k=1,
2, \ldots, m-1$, is evaluated as follows:   \bea\nonumber
\G_{k/2}(F)&=& \pi^{-km/2} \intl_{\Mkm}
F(\om'\om)d\om\\&=&\pi^{-km/2} \intl_{\Mkm} d\om
\intl_{O(n)} f\left(\g \left[\begin{array} {c} (\om'\om)^{1/2} \\
0
\end{array} \right]\right)d\g. \nonumber \eea
By making use of  the polar coordinates, one can write $\om '\in
\Mmk$ as
$$
\om'= \b u_0(\om\om')^{1/2}, \quad \b\in O(m),\quad u_0=\left[\begin{array} {c} I_k \\
0
\end{array} \right]\in V_{m,k}.
$$
Hence, $\om=(\om\om')^{1/2}u'_0\b'$, and

$$
 (\om'\om)^{1/2}=(\b u_0\om\om'u'_0\b')^{1/2}=\b u_0(\om\om')^{1/2}
 u'_0\b'=\b u_0\om.
$$
By changing variable  $\g\to\g\left[ \begin{array}{ll}  \b' & 0 \\
0& I_{n-m}
\end{array}\right]$, we obtain
\bea\nonumber \G_{k/2}(F)&=& \pi^{-km/2}\intl_{\Mkm} d\om
\intl_{O(n)} f\left(\g \left[\begin{array} {c}\b u_0\om\\
0
\end{array} \right]\right)d\g \\[14pt]
&=&\pi^{-km/2}\intl_{\Mkm} d\om
\intl_{O(n)} f\left(\g \left[\begin{array} {c} \om\\
0
\end{array} \right]\right)d\g,\nonumber
\eea and  (\ref{znk}) follows. For $\a=0$, owing to (\ref{r0}), we
have
$$
(\z_0, f)= 2^{-m}\sig_{n,m}F(0)=\frac{\pi^{nm/2}}{\gm(n/2)}f(0).
$$
\end{proof}

\begin{remark}
For $n>m$, the integration in (\ref{znk}) over $O(n)$ can be
replaced by that over $SO(n)$.
\end{remark}
The following formulas for $(\z_k,f)$ are consequences of
(\ref{znk}). We denote \be\label{con1} c_1=
\frac{\pi^{(nm-km-nk)/2} \, \Gam_k(n/2)}{ 2^{k} \, \Gam_m (n/2)},
\qquad c_2= \frac{\pi^{(m-k)(n/2-k)}}{\Gam_k (k/2) \,
\Gam_{m-k}((n-k)/2)}. \ee
\begin{corollary}\label{nf}
For all  $k=1, 2, \ldots,  n$, \be\label{r1}  \index{a}{Zeta
integrals!normalized!of
 integral order}
(\z_k,f)=c_1\intl_{V_{n,k}} dv\intl_{\frM_{k,m}} f(v \om)d\om. \ee
 Moreover, if $k=1, 2, \ldots, m-1$,
then \bea\label{r3}(\z_k,f)&=&c_1\intl_{V_{m,k}} du
\intl_{\frM_{n,k}}f(yu') |y|_k^{m-n} dy \\ \label{r2}
&=&c_2\intl_{\frM_{n,k}}
\frac{dy}{|y|_k^{n-m}}\intl_{\frM_{k,m-k}} f([y; yz]) dz. \eea
\end{corollary}
\begin{proof} From (\ref{znk}) we have
\bea  (\z_k, f)&=&\frac{\pi^{(n-k)m/2}}{\gm(n/2)}
\intl_{\frM_{k,m}}d\om \intl_{O(n)}f(\gam\lam_0 \om)d\gam \quad
\left ( \lam_0= \left[\begin{array} {c} I_k \\ 0
\end{array} \right]  \in V_{n,k} \right )\nonumber \\
&=&\frac{\pi^{(n-k)m/2}}{\sig_{n,k} \,
\gm(n/2)}\intl_{\frM_{k,m}}d\om \intl_{V_{n,k}}f(v \om)dv\nonumber
\eea which coincides with (\ref{r1}).  To prove (\ref{r3}), we
pass to polar coordinates  in (\ref{r1}) by setting
$\om'=ur^{1/2}, \; u \in V_{m,k}, \; r \in \P_k$. This gives \bea
(\z_k, f)&=& 2^{-k} \, c_1
\intl_{V_{n,k}}dv\intl_{\P_k}|r|^{(m-k-1)/2}
dr\intl_{V_{m,k}} f(vr^{1/2}u') du\nonumber \\
&=&c_1\intl_{V_{m,k}} du \intl_{\frM_{n,k}}f(yu') |y|_k^{m-n}
dy.\nonumber \eea To prove (\ref{r2}), we represent $\om$ in
(\ref{r1}) in the block form $u=[\eta; \zeta], \; \eta \in
\frM_{k,k}, \; \zeta \in \frM_{k,m-k}$, and change the variable
$\zeta=\eta z$. This gives
\[ (\z_k, f)=c_1\intl_{\frM_{k,k}}
|\eta|^{m-k}d\eta\intl_{\frM_{k,m-k}}dz \intl_{V_{n,k}}f(v[\eta;
\eta z])dv.\] Using Lemma \ref{l2.3} repeatedly, and changing
variables, we obtain \bea  (\z_k, f)&=&2^{-k} \, c_1  \,
\sig_{k,k}\intl_{\P_k} |r|^{(m-k-1)/2}
dr\intl_{\frM_{k,m-k}}dz \intl_{V_{n,k}}f(v[r^{1/2}, r^{1/2}z])dv\nonumber \\
&=&c_1 \, \sig_{k,k}\intl_{\frM_{n,k}}
\frac{dy}{|y|_k^{n-m}}\intl_{\frM_{k,m-k}} f([y; yz]) dz.\nonumber
\eea By  (\ref{con1}), (\ref{2.16}) and (\ref{2.5.2}), this
coincides with (\ref{r2}). \end{proof}

The representation (\ref{r2})
 was obtained in \cite{Sh1} and \cite{Kh} in a different way.

\begin{remark}
Another proof of Theorem \ref{tzk}, which does not utilize
Faraut's Lemma \ref{lFar}, was given  in \cite{Ru10}. The proof
given above mimics the classical argument (cf.  \cite[Chapter 1,
Section 3.9]{GSh1} in the rank-one case) and is very instructive.
\end{remark}
One can also write $ (\z_k,f)$  as \be\label{z-mes}
(\z_k,f)=\intl_{\Ma} f(x)d\nu_k(x),\qquad f\in\s(\Ma), \ee where
$\nu_k$ is a positive locally finite measure $\nu_k$ defined by
\be\label{mes} (\nu_k ,\vp)\equiv c_1\intl_{V_{n,k}}
dv\intl_{\frM_{k,m}} \vp(v \om)d\om, \quad \vp \in C_c(\Ma),
 \ee
$C_c(\Ma)$ being the space of compactly supported continuous
functions on $\Ma$; cf. (\ref{r1}). In order to characterize the
support of $\nu_k$, we recall that $\Mk$ stands  for the
submanifold of matrices $x\in\Ma$ having rank $k$, and  denote \be
\Mkc = \bigcup_{j=0}^{k} \Mj \qquad \text{\rm (the closure of
$\;\Mk$)}. \ee
\begin{lemma}\label{nu}
The following statements hold.\\
 \noindent{\rm(i)}  $\supp\,
\nu_k=\Mkc$.

\noindent{\rm(ii)} The manifold $ \; \Mk$ is an orbit of
 $\;  e_k=\left[\begin{array}{ll}  I_k&0\\
0&0
\end{array}\right]_{n \times m}$
 under the group $G$ of transformations $$x \to g_1xg_2, \qquad  g_1 \in
 GL(n,\bbr), \quad   g_2 \in
 GL(m,\bbr),$$ and \be\label{dimrk} \dim \Ma^{(k)}= k(n+m-k). \ee

\noindent{\rm(iii)} The manifold $ \; \Mkc$ is a collection of all
matrices $x \in \Ma$ of the form \be\label{form1} x=\gam
\left[\begin{array} {c} \om \\ 0
\end{array} \right], \quad \g \in O(n), \quad \om \in
\frM_{k,m}, \ee or \be\label{form2} x=v\om, \qquad v \in V_{n,k},
\quad \om \in \frM_{k,m}. \ee
\end{lemma}
\begin{proof} (i) Let us consider (\ref{mes}).
Since $\rank (v \om)\leq k$, then $(\nu_k ,\vp)=0$ for all
functions $\vp \in C_c(\Ma)$  supported away from $ \Mkc$, i.e.
$\supp\, \nu_k=\Mkc$.

(ii) We have to show that each $x\in\Mk$ is represented as $x=g_1
e_k g_2$ for some $g_1 \in
 GL(n,\bbr)$ and $ g_2 \in
 GL(m,\bbr)$. We write  $x=us$, where
$u\in\vnm$ and $s=(x'x)^{1/2}$ is a positive semi-definite
$m\times m$ matrix of rank $k$ (see Appendix C, {\bf 11}). By
taking into account that
$s=g'_2 \left[\begin{array}{ll}  I_k&0\\
0&0
\end{array}\right]_{m\times m} g_2 \, $ for some $ g_2 \in
 GL(m,\bbr)$, and $u=\gam
\left[\begin{array} {c} I_m
\\ 0
\end{array} \right]$
for some $\g\in O(n)$, we obtain $x=g_1 e_k g_2$ with
$$ g_1=
\g\left[\begin{array}{ll}  g_2'&0\\
0& I_{n-m}
\end{array}\right]\in
 GL(n,\bbr).$$
 Therefore, the manifold $\Mk$ is a homogeneous space of the group
 $G$. This allows us to find the
 dimension of $\Mk$  by the formula
 $$
\dim \Mk=\dim G-\dim G_1,$$ where $G_1$ is the subgroup of $G$
leaving $e_k$ stable. In order to calculate the   dimension of
$G_1$, we write the condition $g_1 e_k g_2=e_k$ in the form
$$
\left[\begin{array}{ll}  a_1&a_2\\
a_3& a_4
\end{array}\right]\left[\begin{array}{ll}  I_k&0\\
0& 0
\end{array}\right]\left[\begin{array}{ll}  b_1&b_2\\
b_3& b_4
\end{array}\right]=\left[\begin{array}{ll}  I_k&0\\
0& 0
\end{array}\right],
$$
where $a_1\in GL(k,\bbr)$ and $b_1\in GL(k,\bbr)$. This gives
$a_1\,b_1=I_k$, $a_3=0$, $b_2=0$. Hence,
$$\dim G_1=n^2-k(n-k)+m(m-k),$$ and   (\ref{dimrk}) follows.

(iii) It is clear that  each matrix of the form (\ref{form1}) or
(\ref{form2}) has rank $\le k$. Conversely, if $\rank (x) \le k$
then, as above, $x=us =\gam \left[\begin{array} {c} s \\ 0
\end{array} \right]$ where $ \g \in O(n)$ and $s=(x'x)^{1/2}$ is a positive semi-definite
$m\times m$ matrix of rank $\le k$. The latter can be written as
\[ s=g\lam g', \qquad  g \in O(m), \qquad \lam =\diag(\lam_1,
\ldots , \lam_k, 0, \ldots , 0),\] and therefore, \[ x=\gam
\left[\begin{array} {c} g\lam \\ 0
\end{array} \right]g'=\g\left[\begin{array}{ll}  g&0\\
0& I_{n-m}
\end{array}\right]\left[\begin{array} {c} \lam \\ 0
\end{array} \right]g'=\g_1\left[\begin{array} {c} \om \\ 0
\end{array} \right]\]
where $\g_1=\g\left[\begin{array}{ll}  g&0\\
0& I_{n-m}
\end{array}\right]$ and $\om \in \frM_{k,m}$.  The representation (\ref{form2})
 follows from (\ref{form1}).
\end{proof}

\begin{corollary} The integral (\ref{z-mes}) can be written as
\be\label{zzm} (\z_k,f)=\intl_{\text{\rm rank} (x) \le \, k} f(x)
d\nu_k (x)=\intl_{\text{\rm rank} (x) = k} f(x) d\nu_k (x).\ee
\end{corollary}
 \begin{proof} The first equality follows from Lemma \ref{nu} (i). The
 second one is clear from the observation that if $\rank (x) \le
 k-1$ then, by (\ref{form2}), $x=v\om, \; v \in V_{n,k-1},
\; \om \in \frM_{k-1,m}$. The set of all such pairs $(v,\om)$ has
measure $0$ in $V_{n,k} \times \frM_{k,m}$.
\end{proof}

\begin{remark}
We cannot obtain simple explicit representation of the conjugate
normalized zeta distributions $(\z_k^-,f), \; k=0,1, \ldots, m-1$.
At the first glance, it would be natural to use the formula \[
(\z_\a^-,f)=c_\a \, (\z^+_{\a+1}, \DP f), \qquad c_\a=(-1)^m
\frac{\Gam (\a+1-m)}{\Gam (\a+1)}, \] see (\ref{kap1}), in which
$(\z^+_{\a+1}, \DP f)$ can be evaluated for $\a=k$ by Theorem
\ref{tzk} or Corollary \ref{nf}. Unfortunately, we cannot do this
because $c_\a=\infty$ for such $\a \, $. On the other hand, Lemma
\ref{lacz} says that $(\z_k^-,f)$ is well defined by (\ref{ffii}),
namely, \be\label{zmi}(\z_k^-,f)=2^{-m}\, \pi^{-m(m-1)/4}\,
(\om_k, \Phi), \ee where \bea\om_k(t_{1,1},\dots ,
t_{m,m})&=&\left. \prod\limits_{i=1}^{m} \frac{|t_{i,i}|^{\a -i}
\, \sgn \, (t_{i,i})}{\Gam ((\a-i)/2 +1)} \, \right |_{\a=k} \, ,
\nonumber
\\ \Phi (t_{1,1},\dots , t_{m,m})&=&\intl_{\bbr^{m(m-1)/2}} \,
dt_{*} \intl_{O(m)} \, f(vt) \,  \sdv
 \, dv. \nonumber \eea
In particular, for $k=0$, \be\label{zmii} \om_0(t_{1,1},\dots ,
t_{m,m})=\left. \prod\limits_{i=1}^{m} \frac{|t_{i,i}|^{\lam} \,
\sgn \, (t_{i,i})}{\Gam (\lam/2 +1)} \, \right |_{\lam=-i} \, ,\ee
where the generalized functions \[ \left .\frac{|s|^{\lam} \, \sgn
\,  (s)}{\Gam (\lam/2 +1)} \right |_{\lam =-i}, \qquad i=1,2,
\ldots, m,\] are defined as follows. For $i$ odd: \bea \hskip-0.5
truecm  &{}&\left ( \left .\frac{|s|^{\lam} \, \sgn \, (s)}{\Gam
(\lam/2 +1)} \, \right |_{\lam =-i} \, , \vp \right
)=\frac{1}{\Gam (1-i/2)} \, (s^{-i}, \vp) \nonumber
\\&=& \! \intl_0^\infty  \! s^{-i}\, \left \{ \vp(s)-\vp(-s)-2\left [s\vp'
(0) + \frac{s^3}{3!}\, \vp''' (0)+ \ldots +
\frac{s^{i-2}}{(i-2)!}\, \vp^{(i-2)} (0)\right ] \right \} ds.
\nonumber \eea For $i$ even: \[\left ( \left .\frac{|s|^{\lam} \,
\sgn \, (s)}{\Gam (\lam/2 +1)}\,  \right |_{\lam =-i} \, , \vp
\right )=\frac{(-1)^{i/2} \, \vp^{(i-1)} (0) \, (i/2-1)!}{(i-1)!}
\,  ;\] see \cite[Chapter 1, Section 3.5]{GSh1} .  Note that the
Fourier transform of the distribution $\z_0^-$ has the form
\be\label{fsi} (\F \z_0^-)(y)=\frac{(-i)^{m}\, \pi^{m^2/2}}{\Gam_m
((m+1)/2)} \, \, \sdy. \ee This follows immediately from
(\ref{az-fm}) and the Parseval formula (\ref{pars}).

In the following the expression $(\z_0^-,f)$ will be understood in
the sense of regularization according to (\ref{zmi}),
(\ref{zmii}).
\end{remark}

\section{Riesz potentials and the generalized Hilbert transform}
  The functional equation (\ref{az-f})
 for the zeta distribution can be written in the form \index{b}{Greek!c@$\gamma_{n,m} (\a)$}
 \be\label{z-f1}  \frac{1}{\g_{n,m}(\a)}\Z(f, \a-n)=(2\pi)^{-nm}
\Z(\F f, -\a)\; , \ee \be\label{gam} \gam_{n,m} (\a)=\frac{2^{\a
m} \, \pi^{nm/2}\, \Gam_m (\a/2)}{\Gam_m ((n-\a)/2)},\qquad \a\neq
n-m+1, \,  n-m+2, \ldots . \ee We recall that $m\geq 2$. The
excluded  values $\a= n-m+1, \, n-m+2, \ldots\; $ are poles of the
gamma function $\Gam_m ((n-\a)/2)$. Note that for $n=m$, all
$\a=1,2,\dots$ are excluded.  The normalization in (\ref{z-f1})
gives rise to the {\it Riesz distribution}\index{a}{Riesz
distribution, $h_\a$} $h_\a$ defined by
 \bea\label{dis} (h_\a, f)&=&2^{-\a m}\, \pi^{-nm/2} \, \Gam_m
 \left (\frac{n-\a}{2}\right ) \, (\z_\a, f) \\
 &=& \mbox{a.c.} \; \frac{1}{\gam_{n,m} (\a)}
  \intl_{\Ma} |x|^{\a-n}_m f(x) dx, \nonumber
  \eea
 where $f\in\s(\Ma)$ and ``a.c." abbreviates analytic continuation in the
 $\a$-variable. For $Re \, \a>m-1$, the distribution $h_\a$ is
 regular and agrees with the usual function $h_\a(x)=|x|^{\a-n}_m
 /\gam_{n,m}(\a)$.
 The  {\it Riesz potential}  of a function $f\in\s(\Ma)$
 is defined by \be\label{rie-z} \index{a}{Riesz potential!on $\Ma$, $I^\a
f$} (I^\a f)(x)=(h_\a, f_x), \qquad f_x(\cdot)= f(x-\cdot).\ee For
$Re \, \a>m-1$, one can represent $I^\a f$ in the classical form
by the absolutely convergent integral
 \be\label{rie} (I^\a f)(x)=\frac{1}{\gam_{n,m} (\a)} \intl_{\Ma}
f(x-y) |y|^{\a-n}_m dy.\ee \index{b}{Latin and Gothic!ia@$I^\a f$}
We also introduce {\it the generalized Hilbert
transform}\index{a}{Generalized Hilbert transform!on $\Ma$, $\H
f$} \index{b}{Latin and Gothic!h@$\H f$}\be\label{hilb} (\H
f)(x)=(\z_0^-, f_x)= (\z_0^- \ast f)(x).\ee By (\ref{fsi}), this
can be regarded as a pseudo-differential operator with the symbol
\be\label{pdo} (\F \z_0^-)(y)=\frac{(-i)^{m}\, \pi^{m^2/2}}{\Gam_m
((m+1)/2)} \, \, \sdy. \ee Clearly, $\H$ extends as a linear
bounded operator on $L^2(\Ma)$. For $m=1$, it coincides (up to a
constant multiple) with the usual Hilbert transform on the real
line.

The following properties of Riesz potentials and Riesz
distributions
 are inherited from those  for the normalized
zeta integrals. We denote\bea\label{c1} \qquad \quad
\g_1&=&2^{-km} \,\pi^{-km/2}
\, \Gam_m \left(\frac{n-k}{2}\right) / \Gam_m \left(\frac{n}{2}\right),  \\
\label{c2}\g_2&=& 2^{-k(m+1)} \, \pi^{-k(m+n)/2} \, \Gam_k
\left(\frac{n-m}{2}\right). \eea
\begin{theorem}\label{ldes} \index{a}{Riesz potential!on $\Ma$, $I^\a
f$!of integral order} Let $f \in \s(\Ma)$, $n\geq m$. Suppose that
$\a=k$ is a positive integer. If $k \neq n-m+1, \, n-m+2, \ldots
$, then \bea(I^k f)(x)\label{des}&=&\g_1\intl_{\frM_{k,m}}d\om
\intl_{O(n)}f \left (x-\gam \left[\begin{array} {c} \om \\ 0
\end{array} \right]  \right ) \, d\gam, \\
\label{des1}&=& \g_2\intl_{V_{n,k}} dv\intl_{\frM_{k,m}} f(x-v
\om)d\om. \eea  Furthermore, \be\label{I0} (I^{0} f)(x)=f(x). \ee
\end{theorem}

This statement is an immediate consequence of  Theorem \ref{tzk}.
\begin{lemma}\label{bp}
Let $f\in\s(\Ma)$, $n\geq m$, $ \; \a\in \bbc$,  $ \; \a\neq
n-m+1, \, n-m+2, \ldots \, \, .$

\noindent {\rm (i)}  The Fourier transform of the Riesz
 distribution $h_\a$ is evaluated by the formula $(\F
 h_\a)(y)=|y|_m^{-\a}$, the precise meaning of which is
\be\label{fou} (h_\a, f)=(2\pi)^{-nm}(|y|_m^{-\a},(\F
f)(y))=(2\pi)^{-nm} \Z(\F f, -\a). \index{a}{Fourier transform,
$\F f$! of the Riesz
 distribution }\ee

\noindent {\rm (ii)} If $k=0,1,\dots$, and $\Del$ is the
 Cayley-Laplace operator, then \be \label{Dkh}
h_{\a}=(-1)^{mk}\Del^k h_{\a+2k}, \quad i.e. \quad ( h_{\a},
f)=(-1)^{mk}( h_{\a+2k}, \Del^k f).\ee

\noindent {\rm (iii)} If  $n=m$, $k=1,2,\dots\; , $ and $\DP$ is
the Cayley operator, then \be \label{Dkhh} h_{\a}=c \, \DP^{2k-1}
\z^-_{\a+2k-1}, \quad i.e. \quad ( h_{\a}, f)=c \,
(-1)^{m}(\z^-_{\a+2k-1},  \DP^{2k-1}f),\ee \[
c=\frac{(-1)^{m(k+1)} \, \Gam_m (1+(m-\a-2k)/2) }{2^{(\a+2k-1)m}
\,\pi^{m^2/2}}.\]
\end{lemma}
\begin{proof}
 (i) follows immediately from the definition (\ref{dis}) and the
functional equation (\ref{z-f1}). To prove (\ref{Dkh}), for
sufficiently large $\a$, according to (\ref{vaz}),  we have \bea
\Del^k h_{\a+2k}(x)&=&\frac{1}{\g_{n,m}(\a+2k)}\Del^k
|x|^{\a+2k-n}_m \nonumber \\ &=&\frac{B_k(\a)}{\g_{n,m}(\a+2k)}
|x|^{\a-n}_m\nonumber \\ &=& ch_{\a}(x),\nonumber \eea where by
(\ref{bka}) and (\ref{Poh}),
$$
c=\frac{B_k(\a)\,\g_{n,m}(\a)}{\g_{n,m}(\a+2k)}=\frac{B_k(\a)\,\Gam_m
(\a/2)\,\Gam_m ((n-\a)/2-k)}{4^{mk}\,\Gam_m (\a/2+k)\,\Gam_m
((n-\a)/2)}=(-1)^{mk}.
$$
For all admissible $\a\in\bbc$, (\ref{Dkh}) follows by analytic
continuation. Let us prove (iii). Owing to (\ref{Dkh}),
$h_{\a}=(-1)^{mk}\, \DP^{2k-1}\, \DP h_{\a+2k}$. Since, by
(\ref{dis} ) and (\ref{kap1}),
\[ h_{\a+2k}=\frac{\Gam_m ((m-\a-2k)/2)}{2^{(\a+2k)m}
\,\pi^{m^2/2}} \, \, \z^+_{\a+2k}, \quad \text{\rm and} \quad \DP
\z^+_{\a+2k}=\frac{\Gam (\a+2k)}{\Gam (\a+2k-m)} \, \,
\z^-_{\a+2k-1}\;,\] then (\ref{Dkhh}) follows after simple
calculation using (\ref{Poh}).
\end{proof}

Lemma \ref{bp} implies the following.
\begin{theorem}\index{a}{Riesz potential!on $\Ma$, $I^\a
f$!of integral order} Let $f\in\s(\Ma)$,  $n \ge m$.

\noindent {\rm (i)} If $k=0,1,2, \dots$,  then \be
\label{Dkf}(I^{-2k} f)(x)=(-1)^{mk}(\Del^k f)(x). \ee

\noindent {\rm (ii)} If $k=1,2, \dots$, and $n>m$, then \be
\label{Dkk}(I^{1-2k} f)(x)=(-1)^{mk}(I^1\Del^k
f)(x)=c_1\intl_{S^{n-1}}dv\intl_{\bbr^m} (\Del^k f)(x-vy')\, dy,
\ee
\[c_1=\frac{(-1)^{mk}\,
 \Gam ((n-m)/2)}{ 2^{m+1} \, \pi^{(m+n)/2}}.\]

\noindent {\rm (iii)} If $k=1,2, \dots$, and $n=m$, then \be
\label{Dkkk}(I^{1-2k} f)(x)=c_2 \, (\H \DP^{2k-1}f)(x), \qquad
c_2=\frac{(-1)^{m(k+1)}\, \Gam _m((m+1)/2)}{\pi^{m^2/2}},\ee $\H$
being the generalized Hilbert transform (\ref{hilb}).
\end{theorem}
\begin{proof}
The equality (\ref{Dkf}) is a consequence of (\ref{rie-z}),
(\ref{fou}), and (\ref{Dkh}). Namely, \bea\nonumber(I^{-2k}
f)(x)&=&(h_{-2k}, f_x)\stackrel{\rm (\ref{Dkh})}{=}(-1)^{mk}(h_0,
\;\Del^k
f_x)\\[14pt]\nonumber &\stackrel{\rm
(\ref{fou})}{=}&(-1)^{mk}(2\pi)^{-nm}\Z(\F(\Del^k f_x),\;
0)\\[14pt]\nonumber&=&(-1)^{mk}(\Del^k f_x)(0)=(-1)^{mk}(\Del^k
f)(x).\eea Similarly, by (\ref{rie-z}), $(I^{1-2k}
f)(x)=(h_{1-2k}, f_x)$. If  $n>m$, then, by (\ref{Dkh}),
\[(h_{1-2k}, f_x)=(-1)^{mk}(h_1,
\;\Del^k f_x)=(-1)^{mk}(h_1, \; (\Del^k f)_x)=(-1)^{mk}(I^1\Del^k
f)(x),\] and it remains to apply (\ref{des1}) (with $k=1$). If
$n=m$, then we apply (\ref{Dkhh} ) with $\a=1-2k$ and get
\[(h_{1-2k}, f_x)=(-1)^{m} c_2 \, (\z_0^-,
\DP^{2k-1}f_x)=c_2 \,(\z_0^-, (\DP^{2k-1}f)_x)=c_2 \, (\H
\DP^{2k-1}f)(x).\]
\end{proof}

Integral representations  (\ref{rie})  and (\ref{des}) (or
(\ref{des1})) can serve as definitions of $I^\a f$ for
 functions $f$ belonging to Lebesgue spaces and $\a$ belonging to the associated
 Wallach-like
 set \index{a}{Wallach-like set, $\V$} \index{b}{Latin and Gothic!wa@$\V$}
\be\label{Wal-r} \V=\{0, 1, 2, \ldots,  k_0 \}\cup \{\a  : Re\, \a
\!
> \! m \! - \! 1; \; \a\neq n-m+1, \, n-m+2, \ldots\},
\ee
$$
k_0=\min(m-1, n-m).
$$

 \begin{theorem}\label{ap} {\rm \cite{Ru10}.}
\index{a}{Riesz potential!on $\Ma$, $I^\a f$!existence}
  Let $f \in L^p(\Ma), \; n \ge m$. If $\a \in\V$, then the Riesz potential $(I^\a
f)(x)$ absolutely converges for almost all $x \in \Ma$ provided
\be\label{sh}  1 \le p <\frac{n}{Re\,\a +m-1}. \ee
\end{theorem}

For $m=1$, the condition (\ref{sh}) is well known \cite{St1} and
best possible. We do not know whether (\ref{sh}) can be improved
if $m>1$.

\subsection{Inversion of Riesz potentials}

Let us discuss the following problem. Given a Riesz potential
$g=I^\a f$, how do we recover its density $f$?  In the rank-one
case, a variety of pointwise inversion formulas for Riesz
potentials is available in a large scale of function spaces
\cite{Ru1}, \cite{Sa}. However, in the higher rank case  we
encounter essential difficulties.  Below we show how the unknown
function $f$ can be recovered in the framework of the theory of
distributions.

First we specify the space of test functions. The Fourier
transform formula $$(h_\a, f)=(2\pi)^{-nm}(|y|_m^{-\a},(\F
f)(y))$$ reveals that the Schwartz class $\s \equiv \s(\Ma)$ is
not good enough
  because it is not invariant under multiplication by
$|y|_m^{-\a}$. To get around this difficulty, we follow
 the  idea of V.I.
Semyanistyi \cite{Se} suggested for  $m=1$. Let $\Psi \equiv
\Psi(\Ma)$  \index{b}{Greek!psi@$\Psi(\Ma)$} be the subspace of
functions $\psi (y) \in \s$ vanishing on the set of all matrices
of rank $<m$, i.e.,  \be\label{sets} \{y: \, y \in \Ma, \;  \rank
(y)<m \}=\{y: \,y \in \Ma, \; |y'y|=0 \},\ee with all derivatives.
The coincidence of both sets in (\ref{sets}) is clear because
$\rank (y)=\rank (y'y)$, see, e.g., \cite[p. 5]{FZ}. The set
$\Psi$ is a closed linear subspace of $\s$. Therefore, it can be
regarded as a linear topological space with the induced topology
of $\s$. Let $ \Phi\equiv\Phi (\Ma)$
\index{b}{Greek!phi@$\Phi(\Ma)$} be the
 Fourier image of $\Psi$. Since the Fourier transform
$\F$ is an automorphism of $\s$ (i.e., a topological isomorphism
of $\s$ onto itself), then $\Phi$ is a closed linear subspace of
$\s$. Having been equipped with the induced topology of $\s$, the
space $\Phi$ becomes a linear topological space isomorphic to
$\Psi$ under the Fourier transform. We denote by  $
\Phi'\equiv\Phi' (\Ma)$ the space of all linear continuous
functionals (generalized functions) on  $\Phi$.  Since for any
complex $\a$, multiplication by $|y|_m^{-\a}$ is an automorphism
of $\Psi$, then, according to the general theory \cite{GSh2},
$I^\a$, as a convolution with $h_\a$, is an automorphism of
$\Phi$, and we have $$\F[I^\a f](y)=|y|_m^{-\a}(\F f)(y)$$ for all
$\Phi'$-distributions $f$. In the  rank-one case, the spaces
$\Phi$, $\Psi$, their duals  and generalizations were  studied by
P.I. Lizorkin, S.G. Samko and others in view of applications to
the theory of function spaces and fractional calculus; see
\cite{Sa}, \cite{SKM}, \cite{Ru1} and references therein.

\begin{theorem} \index{a}{Riesz potential!on $\Ma$, $I^\a
f$!inversion}Let $\a\in\V$ and let $f$ be a locally integrable
function such that
 $g(x)=(I^\a f)(x)$ is well defined as an absolutely convergent
  integral for almost all $x\in\Ma$. Then
 $f$ can be recovered from $g$ in the sense of
$\Phi'$-distributions by the formula \be\label{inv-r1} (f,\phi)=(
g, I^{-\a}\phi), \qquad \phi \in \Phi, \ee where
\be\label{rd}(I^{-\a}\phi)(x)=(\F^{-1}|y|_m^{\a}\F\phi)(x).\ee In
particular, if $\a=2k$ is even, then \be\label{inv-r2}
(f,\phi)=(-1)^{mk}( g, \Del^k\phi), \qquad \phi \in \Phi,\ee
$\Del$ being the Cayley-Laplace operator (\ref{K-L}).

\end{theorem}

\chapter{Radon transforms}\label{s4}

\setcounter{equation}{0}

\section{Matrix planes}\label{s4.1}

\begin{definition}\label{def-plane}
Let $k,n$, and $m$ be positive integers, $0<k<n$, $\vnk$ be the
Stiefel manifold of orthonormal $(n-k)$-frames in $\bbr^n$, $\Mt$
be the space of $(n-k) \times m $ real matrices. For $\;
\xi\in\vnk$ and $t\in\Mt$, the set
\index{b}{Greek!pa@$\tau(\xi,t)$} \be\label{plane} \tau\equiv
\tau(\xi,t)=\{x:x\in\Ma;\; \eq\} \ee
 will be called a {\it matrix $k$-plane} in $\Ma$. \index{a}{Matrix $k$-plane, $\tau(\xi,t)$} For  $k=n-m$,  the plane
 $\tau$ will be called
a {\it matrix hyperplane}. We denote by $\Gr$  the manifold of all
matrix $k$-
 planes.
\end{definition}
 The parameterization $\tau=
\tau(\xi,t)$ by the points $(\xi,t)$ of a "matrix cylinder" $\cd$
is not one-to-one because for any orthogonal transformation
$\theta\in O(n-k)$, the pairs $(\xi,t)$ and $(\xi\theta ', \theta
t)$ define the same plane $\tau$. We identify  functions
$\varphi(\t)$ on $\Gr$ with the corresponding  functions $\fc$ on
$\cd$ satisfying $\varphi(\xi\theta ',\theta t)=\fc$ for all
$\theta\in O(n-k)$, and supply $\Gr$ with the measure $d\tau$ so
that \be\label{grm} \intl_{\Gr} \varphi(\t) \, d\tau=\intl_{\cd}
 \fc \, d\xi dt. \ee

The plane $\t$ is, in fact, a usual $km$-dimensional plane in
$\Rnm$. To see this, we write  $x=(x_{i,j}) \in\Ma$ and
$t=(t_{i,j}) \in \Mt$ as column vectors \be \bar
x=\left(\begin{array} {c} x_{1,1}  \\ x_{1,2} \\ ... \\ x_{n,m}
\end{array} \right) \in \Rnm, \qquad \bar
t=\left(\begin{array} {c} t_{1,1}  \\ t_{1,2} \\ ... \\ t_{n-k,m}
\end{array} \right) \in \Rnkm,\ee
and denote \be \bar\xi=\diag (\xi, \ldots , \xi) \in
V_{nm,(n-k)m}.\ee Then (\ref{plane}) reads \be\label{4.5.1} \tau=
\tau(\bar \xi,\bar t)=\{\bar x:\bar x\in\Rnm; \; \eqv\}. \ee The
$km$-dimensional planes (\ref{4.5.1})  form a subset of measure
zero in the affine Grassmann manifold of  all $km$-dimensional
planes in $\Rnm$.

 The manifold $\Gr$ can be regarded as a fibre
bundle,  the base of which is the
 ordinary Grassmann manifold $G_{n,k}$ of  $k$-dimensional linear subspaces  of
$\rn$, and whose fibres are homeomorphic to $\Mt$. Indeed, let
$$
\pi:\Gr\to G_{n,k},
$$
be the canonical projection which assigns to each matrix plane
$\tau(\xi,t)$ the subspace
\be\label{4.6}\eta=\eta(\xi)=\{y:y\in\rn;\;\xi ' y=0\} \in
G_{n,k}\;.\ee
 The fiber $H_\eta=\pi^{-1}(\eta)$  is the set of
all matrix planes (\ref{plane}), where  $t$ sweeps the space
$\Mt$.
 Regarding $\Gr$  as a fibre bundle, one can utilize a
parameterization which is alternative to (\ref{plane}) and
one-to-one. Namely, let \be x=[x_1 \dots x_m],  \qquad t=[t_1
\dots t_m], \ee where $x_i\in\rn$ and $t_i\in\rnk$ are
column-vectors. For $\tau= \tau(\xi,t)\in \Gr$, we have
$$\tau=\{x: x\in\Ma;\; \xi ' x_i=t_i,\quad i=1,\dots, m \}.$$ Each
ordinary $k$-plane $\tau_i=\{x_i: x_i\in\rn;\;\xi ' x_i=t_i\}$ can
be parameterized by the pair $(\eta, \lam_i)$, where $\eta$ is the
subspace (\ref{4.6}), $\lam_i\in\eta^\bot, \quad i=1, \ldots , m,$
are columns of the matrix $\lam=\xi t \in \Ma$, and $\eta^\perp$
denotes the orthogonal complement of $\eta $ in $\rn$. The
corresponding
 parameterization
\be\label{p2} \tau =\tau (\eta,\lam), \qquad \eta \in G_{n,k},
\quad \lam=[\lam_1 \dots \lam_m], \quad \lam_i \in \eta^\perp, \ee
is one-to-one. Both parameterizations (\ref{plane}) and (\ref{p2})
will be useful in the following.

\section{Definition and elementary properties of the Radon transform}

The Radon transform $\hat f$ of a function $f(x)$ on $\Ma$ assigns
to $f$ a collection of integrals of $f$ over all matrix $k$-planes
(\ref{plane}).
 Namely, \[ \hat f (\tau)=\int_{x \in \tau} f(x), \qquad \tau \in \Gr.\] In order to
 give this integral precise meaning, we note that the matrix plane
$\tau= \tau(\xi,t), \; \xi\in\vnk$, $t\in\Mt$, consists of
``points"
\[ x=g_\xi \left[\begin{array} {c} \omega \\ t \end{array}
\right],  \] where $\om \in \Mkm,$ and $ g_\xi \in SO(n)$ is a
rotation satisfying \index{b}{Latin and Gothic!gk@$g_\xi$}
  \be\label{4.24}
g_\xi\xi_0=\xi, \qquad \xi_0=\left[\begin{array} {c}  0 \\
I_{n-k} \end{array} \right] \in \vnk. \ee
 This observation leads to the following \index{b}{Latin and Gothic!fh@$\hat f(\tau)$}
\begin{definition}\label{def-radon}\index{a}{Radon transform, $\hat f
(\tau)$} The Radon transform of a function $f(x)$ on $\Ma$ is
defined as a function  on the ``matrix cylinder" $\cd$ by the
formula \be\label{4.9} \hat f (\tau) \equiv \rf=\intl_{\Mkm}
f\left(g_\xi \left[\begin{array} {c} \omega \\t
\end{array} \right]\right)d\omega.
\ee
\end{definition}

It is worth noting that the expression (\ref{4.9}) is independent
of the choice of the rotation $g_\xi:\xi_0\rightarrow\xi$. Indeed,
if $g_1$ and $g_2$ are two such rotations, then $g_1=g_2 g$ where
$g$ belongs to the isotropy subgroup of $\xi_0$. Hence $g$ has the
form \[ g=\left[\begin{array} {cc} \theta & 0\\ 0 & I_{n-k}
\end{array} \right], \qquad \theta \in SO(k). \]
Multiplying matrices and changing variable $\theta \om \to \om$,
we get \bea \intl_{\Mkm} f\left(g_1 \left[\begin{array} {c} \omega
\\t
\end{array} \right]\right)d\omega&=&\intl_{\Mkm} f\left(g_2 \left[\begin{array} {cc} \theta & 0\\ 0 & I_{n-k}
\end{array} \right]
\left[\begin{array} {c} \omega \\t
\end{array} \right]\right)d\omega \nonumber \\
&=&\intl_{\Mkm} f\left(g_2 \left[\begin{array} {c} \omega \\t
\end{array} \right]\right)d\omega. \nonumber \eea
Furthermore, since \be g_\xi \left[\begin{array} {c} \omega \\t
\end{array} \right]=g_\xi
\left[\begin{array} {c} \omega \\0
\end{array} \right]+g_\xi
\left[\begin{array} {c} 0 \\t \end{array} \right]=g_\xi
\left[\begin{array} {c} \omega
\\0 \end{array} \right]+\xi t, \ee
one can write \bea \rf\label{4.1} &=&\intl_{\Mkm} f\left(g_\xi
\left[\begin{array} {c} \omega \\0
\end{array} \right]+\xi t\right)d\omega \\ \label{rad1}
&=&\intl_{\{y\in\Ma \; : \; \xi 'y=0\}} f(y+\xi t) \, dy. \eea

If $m=1$, then $\rf$ is the ordinary $k$-plane Radon transform
that assigns to a function $f(x)$ on $ \bbr^n$ the collection of
integrals of $f$ over all $k$-dimensional planes. If  $m=1$ and $
k=n-1$, the definition (\ref{4.1}) gives the classical hyperplane
Radon transform \cite{H1}, \cite{GGV}.

In terms of the one-to-one parameterization (\ref{p2}), where
$\tau =\tau (\eta,\lam)$,  $\eta \in G_{n,k}, \; \lam=[\lam_1
\dots \lam_m] \in \Ma$, and  $ \lam_i \in \eta^\perp$, the Radon
transform (\ref{4.9}) can be written as \be\index{a}{Radon
transform, $\hat f (\tau)$} \hat f(\tau)=\intl_\eta dy_1\dots
\intl_\eta f([y_1+\lam_1 \dots y_m+\lam_m ]) \, dy_m. \ee For
$m=1$, this is the well known form of the $k$-plane transform in
$\rn$; cf. \cite[p. 30, formula (56)]{H1}.

 The following sufficient conditions  of the existence of the Radon transform
 $\hat f$ immediately follow from  Definition \ref{def-radon}. More subtle results will be
 presented in Section \ref{s5.4}.

\begin{lemma}\label{l4.2}\hskip 10truecm
\index{a}{Radon transform, $\hat f (\tau)$! existence}

 \noindent {\rm(i)} If $f\in L^1(\Ma)$, then the Radon transform $\rf$ exists for
all $\xi\in\vnk$ and almost all $t\in\Mt$. Furthermore,
\be\label{4.41}\intl_{\Mt}\rf\; dt=\intl_{\Ma} f(x) \, dx, \qquad
\forall  \, \xi\in\vnk.\ee

\noindent {\rm (ii)} Let $||x||= (\tr(x'x))^{1/2}=(x_{1,1}^2+
\ldots +x_{n,m}^2)^{1/2}$. If $f$ is a continuous function
satisfying \be\label{4.21}
 f(x)=O(||x||^{-a}), \qquad a>km, \ee then $\rf$ exists  for all
 $\xi\in\vnk$ and all $t\in\Mt$.
\end{lemma}
\begin{proof} (i) is a consequence of the Fubini
theorem: \bea \intl_{\Mt}\rf \;dt&=&\intl_{\Mt} dt \intl_{\Mkm}
f\left(g_\xi \left[\begin{array} {c} \omega \\t
\end{array} \right]\right)d\omega \nonumber \\
 &=&\intl_{\Ma} f(g_\xi x) \,
dx=\intl_{\Ma} f(x) \, dx. \nonumber \eea (ii)
 becomes obvious from (\ref{4.5.1}) if we regard $\tau=\tau (\xi, t)$  as a  $km$-dimensional plane in
 $\Rnm$.
\end{proof}

\begin{lemma}\label{l4.3} Suppose that the Radon transform \[ f
(x) \longrightarrow \rf, \qquad x \in \Ma, \quad (\xi, t) \in \cd,
\]
 exists (at least almost everywhere). Then\\
 \noindent {\rm (i)} $\rf$ is   a ``matrix-even'' function, i.e.,
 \be\label{4.22} \hat f(\xi\theta ',\theta t)=\rf,  \qquad \forall \theta\in O(n-k). \ee
 \noindent {\rm (ii)} The Radon transform  commutes with the group $M(n,m)$ of matrix  motions.
 Specifically,
  if $ \; g(x)=\gamma x\beta +b \; $ where $ \;
\gamma\in O(n), \quad  \b\in O(m)$, and $b\in\Ma$, then
\be\label{4.23} (f \circ g)^\wedge (\xi, t)= \hat
f(\gamma\xi,t\beta +\xi'\gamma 'b). \ee

\end{lemma}
\begin{proof}{\rm (i)}
Formula (\ref{4.22}) is a matrix analog of the ``evenness
property" of the classical Radon transform (the case $m=1$,
$k=n-1$): $\hat f(-\xi,-t)=\rf$ \cite[p. 3]{H1}. By (\ref{4.9}),\[
\hat f(\xi\theta ',\theta t)=
\intl_{\Mkm} f\left(g_{\xi \theta'} \left[\begin{array} {c} \omega \\
\theta t
\end{array} \right]\right)d\omega,\] where one can choose $g_{\xi
\theta'}=g_{\xi} \left[\begin{array} {cc} I_k & 0\\ 0 &
\theta'\end{array} \right]$. Hence
\[ g_{\xi \theta'} \left[\begin{array} {c} \omega \\ \theta t
\end{array} \right]=g_{\xi} \left[\begin{array} {cc} I_k & 0\\ 0 & \theta'\end{array}
\right]\left[\begin{array} {c} \omega \\ \theta t
\end{array} \right]=g_{\xi}\left[\begin{array} {c} \omega \\t
\end{array} \right]\] which gives (\ref{4.22}).

{\rm (ii)} By (\ref{4.9}), \bea (f \circ g)^\wedge (\xi,
t)&=&\intl_{\Mkm} f\left( \gam \left (g_{\xi} \left[\begin{array}
{c} \omega \b \\ t\b
\end{array} \right]\right)  +b \right)d\omega \qquad (\omega \b \to \om) \nonumber \\
&=& \intl_{\Mkm} f\left( g_{\gam \xi} \left ( \left[\begin{array}
{c} \omega  \\ t\b
\end{array} \right] +g_{\gam \xi}^{-1}b \right)\right)
d\omega. \nonumber \eea  We set $g_{\gam
\xi}^{-1}b=\left[\begin{array} {c} p \\q
\end{array} \right]$,  $\; p \in \Mkm, \; q \in \frM_{n-k,m}$,
and change the variable $\om +p \to \om$. Then  the last integral
reads \[ \intl_{\Mkm} f\left( g_{\gam \xi} \left[\begin{array} {c}
\omega  \\ t\b +q
\end{array} \right]\right)d\omega.\]
By (\ref{4.24}), $q=\xi'_0 g_{\gam \xi}^{-1}b=\xi'_0
g'_{\xi}\gam'b=(g_{\xi} \xi_0)'\gam 'b =\xi '\gam 'b,$ and we are
done.

\end{proof}

\begin{corollary}
For $x,y\in\Ma$, let $f_x(y)=f(x+y)$. Then \be \label{4.4}
 \hat f_x(\xi ,t)= \hat f(\xi , \xi 'x+t).  \ee
\end{corollary}
\begin{corollary}\label{c3.6}
If the function $f : \Ma \to \bbc$ is $O(n)$ left-invariant, i.e.,
$f(\gamma x)=f(x)$ for all $\gamma\in O(n)$, then \be \hat
f(\gamma\xi\theta ',\theta t)=\rf, \qquad \forall \; \gamma\in
O(n), \; \theta\in O(n-k).\ee
\end{corollary}

\section[Interrelation]{Interrelation between the Radon transform and the  Fourier transform}
We recall that the Fourier transform  \index{a}{Fourier transform,
$\F f$} of a function $f\in L^1(\Ma)$ is defined by \be\label{4.5}
 (\F f)(y)=\intl_{\Ma} \exp(i \, \tr(y'x)) \, f(x) \, dx, \qquad y\in\Ma.
 \ee
The following statement is a matrix generalization of the
so-called {\it projection-slice theorem}.
\index{a}{Projection-Slice Theorem} It links together  the Fourier
transform (\ref{4.5}) and the Radon transform (\ref{4.9}).

For $y=[y_1 \dots y_m]\in\Ma$, let $\L (y)=\lin(y_1, \dots, y_m)$
be the linear hull of the vectors $y_1, \dots, y_m$, i.e., the
smallest linear subspace containing $y_1, \dots, y_m$. Suppose
that $\rank (y)=\ell$. Then $\dim\L (y)=\ell\leq m$.

 \begin{theorem}\label{CST}
 Let $f\in L^1(\Ma)$, $1\le k\le n-m$. If  $y \in\Ma$, and
 $\zeta$ is a $(n-k)$-dimensional plane containing  $\L (y)$, then for
 any orthonormal frame $\xi\in\vnk$ spanning
 $\zeta$, there exists $b\in\Mt$ such that $y=\xi b$. In this
 case,
 \be\label{4.19}
(\F f)(y)=\intl_{\Mt} \exp (i\, \tr(b't)) \, \rf \, dt, \ee
 or
\be\label{4.20} (\F f)(\xi b)=\F [\hat f(\xi,\cdot)](b), \quad
\xi\in\vnk, \quad b\in\Mt. \ee
 \end{theorem}

\begin{proof} Since each vector $y_j \; (j=1, \ldots , m) \; $ lies in $\zeta$, it
decomposes as   $y_j=\xi b_j$ for some $b_j \in \bbr^{n-k}$. Hence
  $y=\xi b$ where $b=[b_1 \ldots b_m]\in\Mt$. Thus it remains to  prove (\ref{4.20}). By (\ref{4.9}),
$$
\F [\hat f(\xi,\cdot)]( b)=\intl_{\Mt} \exp (i\, \tr(b't))  \, dt
\intl_{\Mkm} f\left(g_\xi \left[\begin{array} {c} \omega \\t
\end{array} \right]\right)d\omega.$$
If $x=g_\xi \left[\begin{array} {c} \omega \\t
\end{array} \right]$, then, by (\ref{4.24}), $$\xi' x=\xi'_0 g'_\xi g_\xi \left[\begin{array} {c} \omega \\t
\end{array} \right]=\xi'_0 \left[\begin{array} {c} \omega \\t

\end{array} \right]=t,$$
and the Fubini theorem yields
$$
\F [\hat f(\xi,\cdot)](b)=\intl_{\Ma} \exp (i \, \tr(b'\xi 'x)) \,
f(x) \, dx=(\F f)(\xi b).
$$
\end{proof}
\begin{remark}\label{rCS} It is clear that $\xi$ and $b$ in the basic equality
(\ref{4.19}) are not uniquely defined. If $\rank(y)=m$ one can
choose some $\xi$ and $b$ as follows. By taking into account that
$n-k \ge m$, we set
$$
\xi_0= \! \left[\begin{array} {c}  0 \\
I_{n-k}
\end{array} \right] \! \in \vnk, \quad \om_0= \! \left[\begin{array} {c}  0 \\
I_{m}
\end{array} \right]  \! \in V_{n-k,m}, \quad v_0= \! \left[\begin{array} {c}  0 \\
I_{m}\end{array} \right] \!  \in \vnm,
$$
so that $\xi_0\om_0=v_0$. Consider the polar decomposition $$y=v
r^{1/2},\qquad v\in\vnk, \quad r=y'y\in\p,$$ and let $g_v$ be a
rotation with the property $g_v v_0=v$. Then $$y=v r^{1/2}=g_v v_0
r^{1/2}=g_v \xi_0\om_0 r^{1/2}=\xi b,$$ where \be\label{4.25}
\xi=g_v\xi_0\in\vnk, \qquad  b=\om_0 r^{1/2}\in\Mt.\ee
\end{remark}

\begin{theorem}\label{inj} \index{a}{Inversion formula!for the Radon
transform} {}\hfil

\noindent {\rm (i)}  \ If $1\le k\le n-m$, then the Radon
transform $f\rightarrow\hat f$ is injective on the Schwartz space
$\s(\Ma)$, and $f$ can be recovered by the formula \bea
\label{4.26} f(x) &=& \frac{2^{-m}}{(2\pi)^{nm}}\intl_{\p}
|r|^{\frac{n-m-1}{2}}dr \nonumber \\ && \\ &\times& \intl_{\vnm}
\exp (-i \, \tr(x'vr^{1/2}))(F\hat f(g_v \xi_0,\cdot))(\xi_0 'v_0
r^{1/2})dv. \nonumber \eea

\noindent {\rm (ii)}  \ For $k>n-m$, the Radon transform  is
non-injective on $\s(\Ma)$.
\end{theorem}
\begin{proof}  By Theorem \ref{CST},  given the Radon transform  $\hat f$  of
 $f \in \s(\Ma)$,
the Fourier transform $(\F f)(y)$ can be recovered at each point
$y \in \Ma$ by the formula (\ref{4.19}), so that if $\hat f \equiv
0$ then $\F f \equiv 0$. Since $\F$ is injective, then $f \equiv
0$, and we are done.    Remark \ref{rCS} allows us  to reconstruct
$f$ from $\hat f$,
 because (\ref{4.25})
 expresses $\xi$ and $b$ through $y \in \Ma$ explicitly. This gives (\ref{4.26}).

 To prove {\rm (ii)}, let $\psi\not\equiv 0$ be a
  Schwartz function, the Fourier transform of which is
supported in the set $\Mm$ of matrices $x\in\Ma$ of rank $m$. This
is an open set in $\Ma$.  By (\ref{4.20}), \be\label{4.14} \F
[\hat \psi(\xi,\cdot)](b)=\hat\psi(\xi b)=0 \quad \forall \;
\xi\in\vnk,\quad  \forall \; b\in\Mt, \ee because $\xi b \notin
\Mm$ (since $n-k<m$, then $\rank(\xi b)<m$). By injectivity of the
Fourier transform in (\ref{4.14}),  we obtain $ \hat \psi(\xi,t)=0
\; \forall \xi, t$. Thus for $k>n-m$, the injectivity of the Radon
transform fails.
\end{proof}

\section{The dual Radon transform}

\begin{definition}\label{def-dual}\index{a}{Dual Radon transform, $\df$}   Let $\tau=\tau (\xi,t)$ be  a matrix
plane (\ref{plane}), $(\xi,t)\in\cd$. We say that the plane $\tau
\equiv \tau (\xi,t) $ contains a ``point" $x\in\Ma$ if
  $\eq$. The  dual Radon transform $\check \vp
(x)$ assigns to a function $\vp(\tau)$ on $\Gr$  the integral of
$\vp$ over all matrix $k$-planes containing $x$. Namely,
\index{b}{Greek!ph@$\df$}
$$
\df=\intl_{\tau\ni x}\vp(\tau), \qquad x\in\Ma.
$$
The precise  meaning of this integral is \bea\label{4.2}\qquad
\df&=&\frac{1}{\sigk}\intl_{\vnk}
\varphi(\xi,\xi'x)\,d\xi \\
&=&\intl_{SO(n)}\varphi(\g\xi_0,\xi_0'\g 'x)\,d\g, \qquad \xi_0=\left[\begin{array} {c}  0 \\
I_{n-k} \end{array} \right] \in \vnk. \nonumber \eea
\end{definition}

Clearly, $\df$ exists for all $x\in\Ma$ if
 $\vp$ is  a continuous function. Later we shall prove that  $\df$
 is finite a.e. on $\Ma$ for any locally   integrable  function  $\vp$.

\begin{remark} The definition (\ref{4.2}) is independent of the parameterization $\tau=\tau
(\xi,t)$ in the sense that for any other parameterization
$\tau=\tau (\xi\theta',\theta t), \; \theta\in O(n-k)$  (see
Section \ref{s4.1}), the equality (\ref{4.2}) gives the same
result:
\[\frac{1}{\sigk}\intl_{\vnk}
\varphi(\xi\theta',\theta\xi'x)d\xi=\frac{1}{\sigk}\intl_{\vnk}
\varphi(\xi_1,\xi'_1 x)d\xi_1=\df.\]
 \end{remark}

\begin{lemma}\index{a}{Duality}
The duality relation

 \bea\label{4.3}
\intl_{\Ma} f(x)\df dx&=&\frac{1}{\sigk}
\intl_{\vnk}d\xi\intl_{\Mt}\fc\rf dt \\
\Big(or \intl_{\Ma} f(x)\df dx&=&\frac{1}{\sigk}\intl_{\Gr} \vp
(\tau) \hat f (\tau) \, d\tau \Big)\nonumber
 \eea
is valid  provided that either side of this equality is finite for
$f$ and $\vp$ replaced by $|f|$ and $|\vp|$, respectively.
\end{lemma}
\begin{proof}
By (\ref{4.9}), the right-hand side of (\ref{4.3})  equals
\be\label{rs} \frac{1}{\sigk} \intl_{\vnk}d\xi\intl_{\Mt}\fc \, dt
\intl_{\Mkm} f\left(g_\xi \left[\begin{array} {c} \omega \\t
\end{array} \right]\right)d\omega.\ee
Setting $x=g_\xi \left[\begin{array} {c} \omega \\t
\end{array} \right]$, we have
\[ \xi'x=(g_\xi \xi_0)'g_\xi \left[\begin{array} {c} \omega \\t
\end{array} \right]=\xi_0'\left[\begin{array} {c} \omega \\t
\end{array} \right]=t.\]
Hence,  by the Fubini theorem,  (\ref{rs}) reads \[
\frac{1}{\sigk} \intl_{\vnk}d\xi\intl_{\Ma} \varphi(\xi,\xi'x)
f(x) \, dx=\intl_{\Ma} \df  f(x) \, dx.\]
\end{proof}

\begin{lemma} The dual
Radon transform commutes with the group $M(n,m)$ of matrix
motions. Specifically, if  $ \; gx=\gamma x\beta +b \; $ where $
\; \gamma\in O(n), \quad  \b\in O(m)$, and $b\in\Ma$, then \be
(\vp \circ g)^\vee (x)= \check \vp(gx). \ee More precisely, if
$\tau=\tau(\xi,t)$, then $$(\vp \circ g)(\xi,t)=\vp(\g\xi, t\b+\xi
' \g'b)$$ and \be(\vp \circ g)^\vee (x)= \check\vp(\gamma
x\beta+b).\ee
\end{lemma}
\begin{proof}
By (\ref{4.2}), \bea (\vp \circ g)^\vee
(x)&=&\frac{1}{\sigk}\intl_{\vnk} \varphi(\g\xi,\xi'x\b+\xi '\g
'b)d\xi \qquad (\g\xi\to\xi)
\nonumber \\
&=&\frac{1}{\sigk}\intl_{\vnk} \varphi(\xi,\xi'(\g
x\b+b))d\xi=\check\vp(\gamma x\beta+b). \nonumber \eea
\end{proof}

\begin{corollary}\label{c4.1}
If $\gamma\in O(n)$ and $ \vp(\gamma\xi,t)=\fc$, then
$\check\vp(\gamma x)=\df$.
\end{corollary}

\section{Radon
transforms of  radial functions}

In this section, we show that the Radon transform and the dual
Radon transform of radial functions are represented by
G{\aa}rding-Gindikin fractional integrals studied in Chapter
\ref{s3}. This phenomenon is well known in the rank-one case when
diverse Radon transforms of radial functions are represented by
the usual Riemann-Liouville fractional integrals. In the higher
rank case, an exceptional role of the G{\aa}rding-Gindikin
fractional integrals in the theory of the Radon transform on
Grassmann manifolds   was demonstrated in \cite{GR}.  We recall
(see Section \ref{s2.4}) that a function $f(x)$ on $\Ma$ is radial
if it is $O(n)$ left-invariant. Each such function has the form
$f(x)=\f0 (x'x)$. In the similar way one can define  radial
functions of matrix planes.

\begin{definition} \index{a}{Radial functions! on $\Gr$} \label{d5.3} For $0<k<n$ and $ m \ge 1$,
let $\Gr$ be the manifold of matrix planes $\tau=\tau(\xi,t)$,
$(\xi,t)\in\cd$; see (\ref{plane}). A function
$\vp(\tau)\equiv\fc$ is called radial if it is $O(n)$
left-invariant in the $\xi$ variable, i.e.,
$\varphi(\gamma\xi,t)=\fc$ for all $\gamma\in O(n)$ and all (or
almost all) $\xi$ and $t$.
\end{definition}
Note that if $\vp$  is a radial function, then \be\label{2.12}
\varphi(\gamma\xi\theta ',\theta t)=\fc, \qquad \forall \gamma\in
O(n), \; \theta\in O(n-k).
 \ee
This equality is a result of parameterization which is not
one-to-one; see Section \ref{s4.1}.

\begin{lemma} Every function $\vp$ on $\Gr$ of the form $\fc=\Fc0 (t't)$ is
radial. Conversely, if  $1\le k \le n-m$ and $\fc$ is radial, then
  there is a function
$\Fc0(s)$ on $\cpm$ such that  $\fc=\Fc0 (t't)$ for all
$\xi\in\vnk$ and all matrices $\, t\in\Mt$  (if $\vp$ is a
``rough" function then  ``all" should be replaced by ``almost
all").
 \end{lemma}
 \begin{proof} The first statement is clear.
  By the polar decomposition (see Appendix C, {\bf 11}), for all matrices $t\in\Mt$
we have \[ t=us^{1/2}, \quad u\in V_{n-k,m}, \quad
  s=t't\in\cpm.\] Let
$$
\xi_0= \left[\begin{array} {c}  0 \\  I_{n-k} \end{array} \right]
\in \vnk, \qquad u_0= \left[\begin{array} {c}  0 \\  I_m
\end{array} \right] \in V_{n-k,m}.
$$
We choose $\gamma\in O(n)$ and $\theta\in O(n-k)$ so that
$u=\theta' u_0$, $\gamma\xi\theta' =\xi_0$. Then \bea
\fc&=&\varphi(\xi,\theta' u_0s^{1/2}) \nonumber \\ &\stackrel{\rm
(\ref{2.12})}{=}&\varphi(\gamma\xi\theta', u_0s^{1/2}) \nonumber \\
&=&\varphi(\xi_0, u_0s^{1/2}) \nonumber \\ &\stackrel{\rm
def}{=}&\Fc0(s). \nonumber \eea This is what we need.
\end{proof}
\label{}
\begin{theorem}\label{t5.5} \index{a}{Radon transform, $\hat f
(\tau)$! of radial functions} Let $f$ be a radial function on
$\Ma$ so that $f(x)=\f0 (r)$, $r=x'x$. Let $I_{-}^{k/2}f_0$ be
 the G{\aa}rding-Gindikin fractional integral  (\ref{4.1b}).
Then \be\label{4.10} \rf=\pi^{km/2} (I_{-}^{k/2}f_0)(s),\qquad
s=t't\in\cpm,
 \ee provided that either side of this equality  exists in the Lebesgue
sense.
\end{theorem}

\begin{proof} The statement follows immediately from  (\ref{4.9}):
\bea\label{5.22} \rf&=&\intl_{\Mkm} f\left(g_\xi
\left[\begin{array} {c} \omega
\\t
\end{array} \right]\right)d\omega= \intl_{\Mkm} f_0( \omega
'\omega+t't) \, d\omega\\&=&\pi^{km/2} (I_{-}^{k/2}f_0)(s).
\nonumber\eea

\end{proof}
\begin{remark} The rank of  $s$   in  (\ref{4.10}) does not exceed $m$. If $\rank(s)=m$
then $s \in \p$. If  $\rank (s)<m$ (it always happens if $k>n- m$)
then $s$ is a boundary point of the cone $\p$. The fact, that for
radial $f$ the Radon transform $\hat f$ is also radial, follows
immediately by Corollary \ref{c3.6}, Definition \ref{d5.1}, and
Definition \ref{d5.3}.
\end{remark}

Let us pass to the dual Radon transform.

\begin{theorem}\label{t5.7} \index{a}{Dual Radon transform, $\df$!of radial functions}
 For $(\xi, t) \! \in  \! \cd$, let $ \; \fc \! = \! \vp_0(t't)$. Suppose that $1\le k\le n-m$ and denote
 \[\Phi_0(s)=|s|^\del \vp_0 (s), \quad \del=(n-k)/2-d, \quad
d=(m+1)/2,\quad c=\frac{\pi^{km/2}\sigma_{n-k,m}}{\sigma_{n,m}}\;
.\] Let $I_{+}^{k/2}\Phi_0$ be the G{\aa}rding-Gindikin fractional
integral (\ref{4.1a}) of $\Phi_0$.
 Then for any $x\in\Ma$ of rank $m$,
 \be\label{4.12}
\df=c |r|^{d-n/2}(I_{+}^{k/2} \Phi_0)(r), \quad r=x'x\in\p,\ee
 provided that either side of this equality  exists in the Lebesgue
sense.
\end{theorem}

\begin{proof} By (\ref{4.2}),
$$
\df=\intl_{SO(n)}\vp_0(x'\gamma\xi_0\xi_o '\gamma' x)d\gamma,
\qquad
\xi_0=\left[\begin{array} {c}  0 \\
I_{n-k}
\end{array} \right] \in \vnk.
$$
We write $x$ in the polar coordinates $ x=vr^{1/2}$, $v \in \vnm$,
$r\in\p$, and get  \bea
\df&=&\intl_{SO(n)}\vp_0(r^{1/2}v'\gamma\xi_0\xi_o
'\gamma'vr^{1/2})\,d\gamma \nonumber \\&=&
\frac{1}{\sigma_{n,m}}\intl_{\vnm}\vp_0(r^{1/2}v'\xi_0\xi_o
'vr^{1/2})\,dv. \label{5.14}\eea
 Since $n-k\geq m$, one can transform (\ref{5.14}) by making use of  the
bi-Stiefel decomposition from Lemma \ref{l2.1}. Let us consider
the cases $k<m$ and $k\ge m$ separately.

$1^0$. {\bf The case $k<m$.} We  set
 \be\nonumber
  v= \left[\begin{array} {cc} a \\ u(I_m -a'a)^{1/2}
 \end{array} \right], \qquad a\in \frM_{k, m}, \quad u \in
 V_{n-k,m}.
 \ee
Multiplying matrices and using (\ref{2.10}), we obtain
 \bea
\df&=&c_1\intl_{0<a'a<I_m}|I_m
-a'a|^\del\vp_0(r^{1/2}(I_m-a'a)r^{1/2}) \, da \quad \nonumber
\\
&=& c_1 \intl_{0<bb'<I_m}|I_m
-bb'|^\beta\vp_0(r^{1/2}(I_m-bb')r^{1/2}) db, \nonumber \eea
 where $c_1=\sigma_{n-k,m}/\sigma_{n,m}$ . Then we pass to polar
 coordinates $$b=vq^{1/2}, \quad v \in \vmk, \quad q \in \P_k, $$
 and use the equality $|I_m-bb'|=|I_k-b'b|$. This gives
 \bea
\df&=&2^{-k}c_1 \intl_{\vmk} dv \intl_0^{I_k}|q|^{(m-k-1)/2}|I_k
-q|^{\del}\vp_0(r^{1/2}(I_m-vqv')r^{1/2}) dq \nonumber \\
& =&2^{-k}c_1 |r|^{-\delta}\intl_{\vmk} dv
\intl_0^{I_k}|q|^{(m-k-1)/2} \Phi_0(r^{1/2}(I_m-vqv')r^{1/2}) dq.
\nonumber \eea Hence, by  (\ref{3.12}),  \[ \df=\pi^{km/2} \,
\frac{\sigma_{n-k,m}}{\sigma_{n,m}} \, |r|^{d-n/2}(I_{+}^{k/2}
\Phi_0)(r),\] and (\ref{4.12}) follows.

$2^0$. {\bf The case $k \ge m$.} Let us transform (\ref{5.14}) by
the formula (\ref{2.11}). We have \bea \df&=&c_2\intl_0^{I_m}
|s|^\gam |I_m-s|^\del \vp_0(r^{1/2}(I_m-s)r^{1/2}) \,  ds
\nonumber \\ &=& c_2\intl_0^{I_m} |I_m-s|^\gam |s|^\del
\vp_0(r^{1/2}s r^{1/2}) \, ds, \nonumber \eea
\[ c_2=\frac{2^{-m} \sigma_{k,m} \,
\sigma_{n-k,m}}{\sigma_{n,m}}, \qquad \gam=k/2 -d.\] By changing
variable $ \; r^{1/2}s r^{1/2}=s_0$ (see Lemma \ref{12.2} (ii)) so
that $|s|=|r|^{-1} |s_0|$, $ds_0=|r|^d ds$, and
\[   |I_m-s|=|r^{-1/2} (r- r^{1/2}s r^{1/2})r^{-1/2}|= |r|^{-1}
|r-s_0|, \] we obtain \bea \df&=&c_2 |r|^{-\gam -\del -d}
\intl_0^r |r-s_0|^\gam |s_0|^\del \vp_0 (s_0) \, ds_0\nonumber
\\ &=&c_2 |r|^{d-n/2} \intl_0^r |r-s_0|^{k/2 -d}\Phi_0(s_0) \,
ds_0. \nonumber \eea By (\ref{2.16}) this coincides with
(\ref{4.12}).
\end{proof}
\begin{remark}
In Theorem \ref{t5.7}, we have assumed $1\le k\le n-m$, although,
Definition \ref{def-dual} is also meaningful for  $k>n-m$.  In the
last case the integral  on the right-hand side of (\ref{4.12})
requires regularization.  Since the main object of our study here
is the Radon transform $\rf$, and the condition $k\leq n-m$
constitutes a natural framework of the inversion problem  (see
Theorem \ref{inj}), we do not focus on the case $k>n-m$ and leave
it for subsequent publications.
\end{remark}

\section{Examples}\label{s5.3}

Formulas (\ref{4.10}) and (\ref{4.12}) allow us to compute the
Radon transform and the dual Radon transform of some elementary
functions. The following   examples are useful in different
occurrences. As above, we suppose $$f \equiv f(x), \quad x\in\Ma,
\qquad \vp \equiv \vp (\xi, t), \quad (\xi, t) \in \cd,$$ and
write $f \xrightarrow { \wedge } \vp$ ($\vp \xrightarrow { \vee }
f$) to indicate that $\vp$ is the Radon transform of $f$ ($f$ is
the dual Radon transform of $\vp$). Given a symmetric $m \times m$
matrix $s$ and $\a \in \bbc$, we denote \be s_{+}^\a= \left \{
\begin{array} {cc} |s|^\a
 & \mbox{ if   $\qquad s\in \p$}, \\
{} \\
 0
  & \mbox{if $\qquad s\not\in\p$}.
\end{array}
\right. \ee

\begin{lemma}\label{ex1}\index{a}{Radon transform, $\hat f
(\tau)$! of elementary functions} Let \be\label{5.17}
\lambda_1=\frac{\pi^{km/2}\Gamma_m((\lam-k)/2)}{\Gamma_m(\lam/2)},
\qquad Re\,\lam
>k+m-1.\ee
The following formulas hold.

\be\label{5.1} |x'x|^{-\lam/2}
 \xrightarrow {  \wedge } \lambda_1  \, |t't|^{(k-\lam)/2},
\ee

 \be \label{5.2}
|I_m+x'x|^{-\lam/2}
 \xrightarrow {  \wedge } \lambda_1 \, |I_m+t't|^{(k-\lam)/2},
\ee

\be\label{5.3} (a-x'x)_{+}^{(\lam-k)/2-d}
 \xrightarrow {  \wedge } \lambda_1 \, (a-t't)_{+}^{\lam/2-d}, \qquad
 a \in \p.
\ee
\end{lemma}
\begin{proof} We  denote  \bea
f_1(x)&=&|x'x|^{-\lam/2}, \nonumber \\  f_2(x) &=&|I_m+x'x|^{-\lam/2}, \nonumber \\
 f_3(x) &=&(a-x'x)_{+}^{(\lam-k)/2-d}.\nonumber \eea Then
$\hat f_i  (\xi, t)=\vp_i (t't)$ where by  (\ref{5.22}), \bea
\vp_1(s)&=&\intl_{\Mkm} |\om'\om+s|^{-\lam/2}d\om,
 \nonumber \\
\vp_2(s)&=&\intl_{\Mkm} |I_m+\om'\om+s|^{-\lam/2}d\om, \nonumber
\\ \vp_3(s)&=&\intl_{\Mkm}
(a-\om'\om-s)_{+}^{(\lam-k)/2-d}d\om.\nonumber \eea Owing to
(\ref{2.15}), $\vp_1(s)=\lam_1 |s|^{(k-\lam)/2}$  and
$\vp_2(s)=\lam_1 |I_m+s|^{(k-\lam)/2}$. For the third integral we
have $\vp_3(s)=0$ if $a-s\not\in\p$. Otherwise,
$$
\vp_3(s)=\!\!\intl_{\{\om\in\Mkm\,:\,\om'\om<a-s\}}\!\!\!\!
|a-\om'\om-s|^{(\lam-k)/2-d}d\om,
$$
 and the formula (\ref{2.15.1}) gives $\vp_3(s)=\lam_1 |a-s|^{\lam/2-d}$.
\end{proof}

\begin{lemma}\label{ex2} \index{a}{Dual Radon transform, $\df$!of elementary functions}
Let $1\le k\le n-m$,  \be\label{5.21}
\lambda_2=\frac{\Gamma_m(n/2)\Gamma_m((\lam-k)/2)}{\Gamma_m(\lam/2)\Gamma_m((n-k)/2)},
\qquad Re\,\lam
>k+m-1. \ee
Then

\be\label{5.4} |t't|^{(\lam-n)/2}
 \xrightarrow { \vee } \lambda_2 |x'x|^{(\lam-n)/2},
\ee

 \be\label{5.5}
|t't|^{(\lam-n)/2}|I_m+t't|^{-\lam/2}
 \xrightarrow { \vee } \lambda_2 |x'x|^{(\lam-n)/2}|I_m+x'x|^{(k-\lam)/2},
\ee

\be\label{5.6} |t't|^{(k-n)/2+d}(t't-a)_{+}^{(\lam-k)/2-d}
 \xrightarrow { \vee } \lambda_2
 |x'x|^{d-n/2}(x'x-a)_{+}^{\lam/2-d}.
\ee
\end{lemma}

\begin{proof} Let  \bea\nonumber
\vp_1(\xi, t)&=&|t't|^{(\lam-n)/2},\\\nonumber \vp_2(\xi,
t)&=&|t't|^{(\lam-n)/2}|I_m+t't|^{-\lam/2},\\\nonumber \vp_3(\xi,
t)&=&|t't|^{(k-n)/2+d}(t't-a)_{+}^{(\lam-k)/2-d}.\nonumber \eea
 By (\ref{4.12}),
  \be\label{5.7} \check\vp_i(x)=c|r|^{d-n/2}(\I+k
\Phi_i)(r),\qquad r=x'x ,\quad i=1,2,3,\ee
 where $c=\pi^{km/2}\sigma_{n-k,m}/\sigma_{n,m}$,
\bea \Phi_1(s)&=&|s|^{(\lam-k)/2-d}, \nonumber \\
\Phi_2(s)&=&|s|^{(\lam-k)/2-d}|I_m+s|^{-\lam/2},\nonumber \\
\Phi_3(s)&=&(s-a)_{+}^{(\lam-k)/2-d}.\nonumber \eea For $
\check\vp_1(x)$, owing to (\ref{5.7}) and (\ref{4.1a}), we obtain
\[
\check\vp_1(x)=c|r|^{d-n/2}\pi^{-km/2}\intl_{\{\om\in\Mkm:\;\om'\om<r\}}|r-\om'\om|^{(\lam-k)/2-d}d\om.\]
 Hence (\ref{2.15.1}) and (\ref{2.16}) yield
$$
\check\vp_1(x)=\frac{\Gamma_m(n/2)\;\Gamma_m((\lam-k)/2)}{\Gamma_m(\lam/2)\;\Gamma_m((n-k)/2)}
\; \; |r|^{(\lam-n)/2},\qquad r=x'x.
$$
This coincides with (\ref{5.4}).

For $ \check\vp_3(x)$, according to (\ref{5.7}), we have
$$
\check\vp_3(x)=c|r|^{d-n/2}\pi^{-km/2}\intl_{\{\om\in\Mkm:\;\om'\om<r\}}
(r-\om'\om-a)_{+}^{(\lam-k)/2-d}d\om.
$$
Hence $ \check\vp_3(x)=0$ if $r-a\not\in\p$, and
$$
\check\vp_3(x)=c|r|^{d-n/2}\pi^{-km/2}\intl_{\{\om\in\Mkm:\;\om'\om<r-a\}}
|r-a-\om'\om|^{(\lam-k)/2-d}d\om
$$
 if $r-a\in\p$. Applying (\ref{2.15.1}), we obtain
$$
\check\vp_3(x)=\frac{\Gamma_m(n/2)\;\Gamma_m((\lam-k)/2)}{\Gamma_m(\lam/2)\;\Gamma_m((n-k)/2)}\;|r|^{d-n/2}\;|r-a|^{\lam/2-d}.
$$
This gives (\ref{5.6}).

 In order to prove (\ref{5.5}), we consider the cases $k \ge
 m$ and $k < m$ separately.

  $1^0$. {\bf The case $k \ge m$.}
By (\ref{5.7}),
 \bea\nonumber
 \check\vp_2(x)&=&c_1|r|^{d-n/2}\,\intl_0^r
|s|^{(\lam-k)/2-d}\;|r-s|^{k/2-d}\;
|I_m+s|^{-\lam/2}\;ds\\\nonumber
 &=&c_1
|r|^{d-n/2}\intl_{I_m}^{I_m+r}
|r+I_m-s|^{k/2-d}\;|s-I_m|^{(\lam-k)/2-d}\; |s|^{-\lam/2}\;ds,
\eea
 $c_1=2^{-m}\sigma_{n-k,m}\sigma_{k,m}/\sigma_{n,m}$. Owing to
 (\ref{2.9}), we obtain
\be\nonumber \check\vp_2(r)=c_1 B_m(k/2, (\lam -k)/2) \,
|r|^{(\lam-n))/2}\;|I_m+r|^{(k-\lam)/2},\ee and (\ref{5.5})
follows.

$2^0$. {\bf The case $k<m$.} By (\ref{5.7}) and (\ref{3.7}),
\[ \check\vp_2(x)=c\;|r|^{d-n/2}(\I+k
\Phi_2)(r)=c\;|r|^{(k-n)/2}(I_{-}^{k/2}G_2)(r^{-1})\] where
$$
G_2(s)=|s|^{-k/2-d}\Phi_2(s^{-1})=|I_m+s|^{-\lam/2}.
$$
 Hence, in view of  (\ref{3.7}),
\[  \check\vp_2(x)=c\;\pi^{-km/2}\,|r|^{(k-n)/2} \intl_{\Mkm}|I_m+r^{-1}+\omega
'\omega|^{-\lam/2}d\omega.\] This integral can be evaluated by the
formula (\ref{2.15}), and we obtain  \bea\nonumber
\check\vp_2(x)&=&\frac{c\; \Gamma_m((\lam-k)/2)}{\Gamma_m(\lam/2)}
 \; |r|^{(k-n)/2}\;|I_m+r^{-1}|^{(k-\lam)/2}\\
&=& \lam_2 \, |r|^{(\lam-n)/2}\;|I_m+r|^{(k-\lam)/2}, \qquad
r=x'x. \nonumber \eea
\end{proof}

\section{Integral identities. Existence of the Radon transform} \label{ident}

\subsection{Integral identities}
 It is important to have precise information about the behavior
of the Radon transform at infinity and near certain manifolds. The
same is desirable for the dual Radon transform. Examples in
Section \ref{s5.3} combined with  duality (\ref{4.3}) give rise to
integral identities which provide this information in integral
terms. For $m=1$, similar results were obtained in \cite{Ru7}; see
also \cite{Ru9} for Radon transforms on affine Grassmann
manifolds.

We present these results in two theorems. The constants
$\lambda_1$ and $ \lambda_2$ below are defined by (\ref{5.17}) and
(\ref{5.21}), respectively, $d=(m+1)/2$, and all equalities  hold
provided that either side exists in the Lebesgue sense.

\begin{theorem}
If $1\le k\leq n-m$, $Re\, \lam>k+m-1$, then \bigskip
\be\label{5.8}
 \begin{array}{ll}
{\displaystyle \frac{1}{\sigk} \intl_{\vnk}d\xi\intl_{\Mt}\rf\,
|t't|^{(\lam-n)/2}\,dt} \\[30pt]
={\displaystyle \lambda_2\intl_{\Ma} f(x)\,|x'x|^{(\lam-n)/2}\,
dx},
\end{array}
\ee

\be\label{5.9}
 \begin{array}{ll}
{\displaystyle \frac{1}{\sigk} \intl_{\vnk}d\xi\intl_{\Mt}\rf\,
|t't|^{(\lam-n)/2}\,|I_m+t't|^{-\lam/2}\,dt} \\[30pt]
={\displaystyle \lambda_2\intl_{\Ma}
f(x)\,|x'x|^{(\lam-n)/2}\,|I_m+x'x|^{(k-\lam)/2} \,dx},
\end{array}
\ee

\be\label{5.10}
 \begin{array}{ll}
{\displaystyle \frac{1}{\sigk} \intl_{\vnk}d\xi\intl_{\Mt}\rf\,
|t't|^{-\del}\,(t't-a)_{+}^{(\lam-k)/2-d}\,dt} \\[30pt]
={\displaystyle \lambda_2\intl_{\Ma}
f(x)\,|x'x|^{d-n/2}\,(x'x-a)_{+}^{\lam/2-d} \,dx},
\end{array}
 \ee
 $$
 \del=(n-k)/2-d.$$
 \end{theorem}

\begin{theorem}
If $\;1\le k\le n-m$, $Re\, \lam>k+m-1$, then

 \be\label{5.11}
 \begin{array}{ll}
 {\displaystyle \intl_{\Ma} \df\,|x'x|^{-\lam/2}\, dx}
 \\[30pt]
={\displaystyle \frac{\lambda_1}{\sigk}
\intl_{\vnk}d\xi\intl_{\Mt}\fc \,|t't|^{(k-\lam)/2}\,dt},
\end{array}
\ee \be\label{5.12}
 \begin{array}{ll}
 {\displaystyle \intl_{\Ma} \df\,|I_m+x'x|^{-\lam/2} \,dx}
 \\[30pt]
={\displaystyle \frac{\lambda_1}{\sigk}
\intl_{\vnk}d\xi\intl_{\Mt}\fc \,|I_m+t't|^{((k-\lam)/2}\,dt},
\end{array}
\ee \be\label{5.13}
 \begin{array}{ll}
 {\displaystyle \intl_{\Ma} \df \, (a-x'x)_{+}^{(\lam-k)/2-d} \,dx}
 \\[30pt]
={\displaystyle \frac{\lambda_1}{\sigk}
\intl_{\vnk}d\xi\intl_{\Mt}\fc \,(a-t't)_{+}^{\lam/2-d}\, dt}.
\end{array}
\ee
 \end{theorem}

\subsection{Existence of the Radon transform}\label{s5.4}

In most of the formulas presented above we a priori assumed that
the Radon transform  is well defined. Now we arrive at one of the
central questions: for which functions $f$ does the Radon
transform  $\rf$  exist? In other words, which functions are
integrable over all (or almost all) matrix planes $\xi'x=t$ in
$\Ma$? The
 crux is how to specify  the behavior of $f(x)$ at infinity.
Below we  study this problem if (a) $f$ is a continuous function,
(b) $f$ is a locally integrable function,  and (c) $f \in L^p$.

\begin{theorem}\label{t5.10}\index{a}{Radon transform, $\hat f
(\tau)$! existence} Let $f(x)$ be a continuous function on $\Ma$,
satisfying \be\label{5.15} f(x)=O(|I_m+x'x|^{-\lam/2}). \ee If
$\lam>k+m-1$
 then the Radon transform $\rf$  is finite for all $(\xi, t) \in \cd$.
  If $\lam \le k+m-1$ then there is a function $f_\lam (x)$ which
  obeys (\ref{5.15}) and $\hat f_\lam (\xi, t) \equiv \infty$.
\end{theorem}
\begin{proof} It suffices to consider the function $f_\lam
(x)=|I_m+x'x|^{-\lam/2}$. By (\ref{5.2}), $
 \hat f_\lam(\xi,t)= \lambda_1 \, |I_m+t't|^{(k-\lam)/2}
$ provided  $\lam>k+m-1$. This proves the first statement of the
theorem. Conversely, owing to (\ref{4.10}), the Radon transform  $
\hat f_\lam(\xi,t)$ is finite if and only if
$(I_-^{k/2}f_\lam)(t't) < \infty$ (we utilize the same notation
$f_\lam$  both for $|I_m+x'x|^{-\lam/2}$ and $|I_m+r|^{-\lam/2}$,
$r=x'x$). If $\lam\leq k+m-1$, then $I_-^{k/2}f_\lam \equiv
\infty$; see  Remark  \ref{r3.15}. Thus $\hat f_\lam (\xi, t)
\equiv \infty$ for the same $\lam$.
\end{proof}

\begin{theorem}\label{t5.2}
Let $f$ be a locally integrable function on $\Ma$. If
\be\label{5.18} \intl_{\{y\in\Ma\,:\, y'y>R\}}
|y'y|^{(k-n)/2}|f(y)| \, dy<\infty \qquad \text{for all $ \; R \in
\p$},
 \ee
then $\rf$  is finite for almost all $(\xi, t) \in \cd$. If
(\ref{5.18}) fails for some $ R \in \p$, then $\rf$ may be
identically infinite.
\end{theorem}
\begin{proof} STEP 1.
Suppose  first that $f$ is a radial function so that $f(x)=\f0
(r)$, $r=x'x$. Then (\ref{5.18}) becomes
\be\label{5.19}\intl_{R}^\infty |r|^{k/2-d}\,|f_0(r)|\, dr<\infty
\qquad \text{for all $ \; R \in \p$}.\ee By Theorem \ref{t5.5},
$$
\rf=\pi^{km/2} (I_{-}^{k/2}f_0)(t't),$$ and therefore, $\rf
<\infty$ if and only if $(I_{-}^{k/2}f_0)(t't)<\infty$.
 Now the result follows from Lemma \ref{l3.14} and Remark
\ref{r3.15}.

STEP 2. The general case reduces to the radial one. Indeed, let
$$ F(x)=\intl_{SO(n)}f(\gam x)d\gam $$ be the ``radialization" of
$f$. If (\ref{5.19}) holds, then, by STEP 1, the Radon transform
of $F$ is a radial function of the form $\hat F(\xi, t)=\Phi_0
(s), \; s=t't$, which is finite for almost all $t$.  On the other
hand, owing to  commutation (\ref{4.23}), $\Phi_0 (t't)$ is the
mean value of $\rf$ over all $\xi \in \vnk$, and therefore,
$$ \frac{1}{\sigk}\intl_{\vnk}\rf d\xi= \Phi_0 (t't)<  \infty
$$ for almost all $t$. It
follows that $\rf < \infty$  for almost all $(\xi, t) \in \cd$.
\end{proof}

\begin{theorem}\label{t5.1} \index{a}{Radon transform, $\hat f
(\tau)$! existence} Let $f\in\lp$. The Radon transform $\rf$  is
finite for almost all $(\xi, t) \in \cd$  if and only if
\be\label{lp} 1\leq p<p_0=\frac{n+m-1}{k+m-1}. \ee
\end{theorem}
\begin{proof}
By H\"older's inequality,
$$\intl_{\{y\in\Ma\,:\, y'y>R\}} |y'y|^{(k-n)/2}|f(y)| \, dy\leq A\nf,$$
where \bea\label{5.20} A^{p'}&=& \intl_{\{y\in\Ma\,:\,y'y>R\}}
|y'y|^{p'(k-n)/2}\,dy \nonumber \\
 &=&c \intl_{R}^\infty
|r|^{p'(k-n)/2+n/2-d}dr \qquad (d=(m+1)/2, \; c=c(n,m)) \nonumber
\\
 &=&c \intl_0^{R^{-1}} |r|^{p'(n-k)/2-n/2-d}\,dr<\infty \nonumber \eea
provided $p<p_0$. Thus, by Theorem \ref{t5.2}, the Radon transform
$\rf$ is  finite for almost all $(\xi, t) \in \cd$. Conversely, if
$p>p_0$, then one can choose $\lam$ satisfying $$ \frac{n+m-1}{p}
<\lam \le k+m-1.$$ For such $\lam$, the function $f_\lam (x)=|I_m
+x'x|^{-\lam /2}$  belongs to $\lp$ (use the formula
(\ref{2.15})), and $\hat f_\lam (\xi, t) \equiv \infty$; see the
proof of Theorem \ref{t5.10}. In order to cover the case $p=p_0$,
we need a more subtle counter-example. One can show that for $p
\ge p_0$, the function
\be\label{5.16}F(x)=|2I_m+x'x|^{-(n+m-1)/2p}(\log|2I_m+x'x|)^{-1}\ee
 belongs to $\lp$, and $\hat
F(\xi,t)\equiv\infty$. The proof of this statement is technical,
and presented in Appendix B.
\end{proof}
\begin{remark} {\rm(i)} Theorems \ref{t5.10} and \ref{t5.1} agree with each
other in the sense that
  all  functions satisfying (\ref{5.15}) belong to
 $\lp$ where $$\frac{n+m-1}{\lam} <p<p_0.$$ The last statement  holds by the formula (\ref{2.15}),
 according to which  the function $|I_m +x'x|^{-\lam p/2}$ is
integrable for such $p$.
\end{remark}

{\rm(ii)} The equality (\ref{5.9}) implies  part (ii) of Theorem
\ref{t0.5} and provides alternative proof of the ``if part" of
Theorem \ref{t5.1}. Indeed, choosing $\lam>\max\{n,\; k+m-1\}$, we
have $|x'x|^{(\lam-n)/2}\leq |I_m+x'x|^{(\lam-n)/2}$, and
\be\label{4.27}
 \begin{array}{ll}
{\displaystyle \frac{1}{\sigk} \intl_{\vnk}d\xi\intl_{\Mt}|\rf|\,
|t't|^{(\lam-n)/2}\,|I_m+t't|^{-\lam/2}\,dt} \\[30pt]
\leq {\displaystyle \lambda_2\intl_{\Ma}
|f(x)|\,|I_m+x'x|^{(k-n)/2} \,dx}.
\end{array}
\ee Hence, if $|I_m+x'x|^{\a/2}f(x)\in L^1(\Ma)$ for some $\a\geq
k-n$,  then $\rf$ is finite for almost all $(\xi,t) \in \cd$.
Moreover, by H\"older's inequality, the right-hand side of
(\ref{4.27}) does not exceed $A\nf$, where \be A^{p'}=\intl_{\Ma}
|I_m+x'x|^{p'(k-n)/2}dx.  \ee  By (\ref{2.15}), the last integral
is finite when $p'(k-n)>k+m-1$, i.e.,  for $p<p_0$. Thus the
left-hand side of (\ref{4.27}) is finite too, and therefore, the
Radon transform $\rf$ is well defined for almost all $(\xi,t) \in
\cd$.

\subsection{Existence of the dual Radon transform}
The special case  $\lam=k+2d$ in
 (\ref{5.13}) deserves particular mentioning. In this case
\be\label{du}
  \intl_{x'x<a} \df  dx
= \frac{\lambda_1}{\sigk} \intl_{\vnk}d\xi\intl_{\Mt}\fc\,
(a-t't)_{+}^{k/2}\,dt, \ee  and therefore,  \be
  \intl_{x'x<a} |\df|\,  dx
\le c \intl_{\vnk}d\xi\intl_{t't<a}|\fc| \,dt, \quad c=
\frac{\lambda_1 \, |a|^{k/2} }{\sigk}. \ee Since $a \in \p$ is
arbitrary, this implies  the following.
\begin{theorem}  \index{a}{Dual Radon transform, $\df$!existence}If $\fc$ is a  locally integrable function  on
the set $\cd$, $1\leq k\le n-m$, then the dual Radon transform
$\df$ is finite
 for almost all $x \in \Ma$.
\end{theorem}

\chapter{Analytic families associated to the Radon
transform}\label{s5}

\section{Matrix distances and shifted  Radon transforms}

In this section, we  introduce important mean value operators
which can be called  {\it the shifted  Radon transforms. } Here we
adopt the terminology of F. Rouvi\`ere \cite{Rou} that  admirably
suits our case. In contrast to the rank-one case $m=1$, where
averaging parameters are positive numbers, in the higher rank case
$m>1$ these parameters are matrix-valued.

Some geometric preliminaries are in order. We shall need a natural
analog of the euclidean distance for the space $\Ma$ of
rectangular $n\times m$ matrices. Unlike  $m=1$, when the distance
is represented by a positive number, in the  higher rank case our
``distance'' will be matrix-valued and represented by a positive
semi-definite $m\times m$ matrix.

\begin{definition}\label{def-dist}
A {\it matrix distance} between two points $x$ and $y$ in $\Ma$ is
defined by \be\label{md}\index{a}{Matrix distance} \index{b}{Latin
and Gothic!dx@$d(x,y)$} d(x,y)=[(x-y)'(x-y)]^{1/2}.\ee Given a
point $x\in\Ma$ and a matrix $k$-plane $\tau=\tau(\xi, t) \in\Gr$
(see (\ref{plane})), a matrix distance between  $x$ and $\tau$ is
defined by \index{b}{Latin and Gothic!dxt@$d(x,\tau)$}
 \be\label{dist}  d( x, \tau)
=[(\xi ' x-t)'(\xi ' x-t)]^{1/2}.\ee Abusing  notation,  we  write
(cf. (\ref{xm})) \index{b}{Latin and
Gothic!xy@\texttt{"|}$x-y$\texttt{"|}$_m$}  \index{b}{Latin and
Gothic!xt@\texttt{"|}$x-\tau$\texttt{"|}$_m$} \be\label{detdist}
|x-y|_m=\det(d(x,y)),\qquad |x-\tau|_m=\det(d(x,\tau)).\ee
\end{definition}

Let us comment on this definition. We first note that for $m=1$,
(\ref{md}) is the usual euclidean distance between points in
$\rn$. If $y=0$, then $d(x,0)=(x'x)^{1/2}$. This agrees with the
polar decomposition $x=vr^{1/2}$, $r=x'x$, $v\in\vnm$. If $\rank
(x)=m$, then $d(x,0)=r^{1/2}$ is a positive definite matrix. If
$\rank (x)<m-1$ then  $d(x,0)\in \d\p$ and $\det(d(x,0))=0$.

Let us explain (\ref{dist}). For $\xi\in\vnk$, let $\{\xi\}$
denote the $(n-k)$-dimensional subspace spanned by $\xi$, and
$\Pr_{\{\xi\}}$ the orthogonal projection onto $\{\xi\}$ which is
a linear map with  $n\times n$  matrix $\xi \xi '$. Then, as in
the euclidean case, it is natural to define the distance between
the point $x\in\Ma$ and the plane $\tau=\tau(\xi, t) \in\Gr$ as
that between two points, namely, $\Pr_{\{\xi\}} x=\xi \xi 'x$ and
$\xi t$. By (\ref{md}), we have
$$d( x, \tau)=d(\xi
\xi 'x, \xi t)=[(\xi \xi 'x - \xi t)'(\xi \xi 'x- \xi
t)]^{1/2}=[(\xi ' x-t)'(\xi ' x-t)]^{1/2}.$$ This can be regarded
as a matrix distance between $x$ and the projection $\Pr_\tau x$
of $x$ onto the plane $\tau$.

\begin{lemma}\label{l5.2} {}\hfil

\noindent{\rm(i)} \ The group  $M(n,m)$   of matrix motions,
 $$\Ma\ni x\longrightarrow \g x\b+b,\qquad\g\in O(n),\qquad \b\in
 O(m),\qquad b\in\Ma,$$ acts on the manifolds $\Ma$ and $\Gr$ transitively.

\noindent{\rm(ii)} \ The determinants $|x-y|_m$ and $|x-\tau|_m$
are invariant under the action of $M(n,m)$. Namely,
$$
|gx-gy|_m=|x-y|_m,\qquad |gx-g\tau|_m=|x-\tau|_m,\qquad g\in
M(n,m).
$$

\noindent{\rm(iii)} \ The distances $d( x, y)$ and $d( x, \tau)$
are invariant under the subgroup $M(n)$ of $M(n,m)$, acting by the
rule $x \to\g x+b$, $\g\in O(n)$, $b\in\Ma$ (i.e. $\b=0$).
\end{lemma}

\begin{proof}
The statements follow immediately from  Definition \ref{def-dist},
by taking into account that \be\label{gtau} g:\;
\tau(\xi,t)\longrightarrow \tau(\g\xi, t\b+\xi '\g'b). \ee

\end{proof}

\begin{definition} Let $r\in\cpm$.
 The {\it shifted  Radon transform} \index{a}{Shifted Radon transform, $\hat f_r(\t)$} of a function $f(x)$ on $\Ma$ is
 defined by
 \bea\label{shr} \index{b}{Latin
and Gothic!fr@$\hat f_r(\tau)$} \hat f_r(\xi,t)&=&
\frac{1}{\sigma_{n-k,m}}\intl_{V_{n-k,m}}dv \intl_{\frM_{k,m}}
f\left(g_\xi \left[\begin{array} {c} \om
\\t+vr^{1/2}
\end{array} \right]\right)d\om\\[14pt]
&=& \frac{1}{\sigma_{n-k,m}}\intl_{V_{n-k,m}} \hat f(\xi,
t+vr^{1/2})dv\label{shr1},
  \eea
  where  $(\xi,t)\in\cd$, $ g_\xi \in SO(n)$ is a
rotation satisfying (\ref{4.24}).
  \end{definition}

If $r=0$ then $\hat f_r(\xi,t)$ coincides with the Radon transform
$\rf$; cf. (\ref{4.9}).
 The integral (\ref{shr}) averages $f(x)$ over all $x$ satisfying
 $\xi' x=t+vr^{1/2}$. This means that $\hat f_r(\xi,t)$ is the integral of $f(x)$ over all $x$ at matrix distance
 $r^{1/2}$ from the plane $\tau=\tau(\xi,t)$ (see Definition
 \ref{def-plane}).
 One can write
\be\label{dualsh} \hat f_r (\tau)=\hat f_r(\xi,t)= \intl_{d(x
,\tau)=r^{1/2}} f(x).\ee

 \begin{definition}\label{shift}
 Let $1\le k\leq n-m$, $r\in\cpm$,
 $\Gr$ be the manifold of all matrix $k$-planes
$$ \tau\equiv
\tau(\xi,t)=\{x\in\Ma:\eq\},\qquad (\xi,t)\in\cd. $$ Given a point
$z\in\Mt$ at distance $r^{1/2}$ from the  origin (i.e. $z'z=r$),
the {\it shifted dual Radon transform}\index{a}{Shifted dual Radon
transform, $\check \vp_r (x)$}  of a function $\vp(\tau)\equiv\fc$
on $\Gr$
 is defined by \index{b}{Greek!pha@$\check \vp_r (x)$}

 \be\label{6.1}
\check \vp_r (x)= \frac{1}{\sigk}\intl_{\vnk}
\varphi(\xi,\xi'x+z)d\xi.
  \ee

\end{definition}

Let us comment on this definition. We first note that (\ref{6.1})
is independent of the choice of $z$ with $z'z=r$. Indeed, by
passing to polar coordinates
 \[ z=\theta u_0 r^{1/2}, \quad \theta \in O(n-k), \quad u_0=\! \left[\begin{array} {c}  0 \\
I_{m}
\end{array} \right]  \! \in V_{n-k,m}, \]
 we have
\bea \check \vp_r (x)&=& \frac{1}{\sigk}\intl_{\vnk}
\varphi(\xi,\xi'x+\theta u_0 r^{1/2})\,d\xi\qquad \mbox{ (set
$\xi=\eta \theta'$)} \nonumber
\\\nonumber &=& \frac{1}{\sigk}\intl_{\vnk} \varphi(\eta \theta',
\theta
(\eta'x+ u_0 r^{1/2})\,d\eta \quad (\mbox{use} \quad\varphi(\xi\theta ',\theta t)=\fc)\\
 &=& \frac{1}{\sigk}\intl_{\vnk} \varphi(\eta ,
\eta'x+ u_0 r^{1/2})\,d\eta. \label{7.1} \eea

If $r=0$ then $\check \vp_r (x)$   is the usual  dual Radon
transform  (\ref{4.2}). Owing to (\ref{dist}), the matrix distance
$d( x, \tau)$ between $x$ and $\tau=\tau(\xi, \xi'x+z)$ is
$r^{1/2}$.  Hence,  $\check \vp_r (x)$ may be regarded as the mean
value of $\fc$ over all matrix planes  at distance $r^{1/2}$ from
$x$: \be \check \vp_r (x)=\intl_{d(x ,\tau)=r^{1/2}} \vp(\tau)\;
.\ee

\section{Intertwining operators}

\index{a}{Intertwining operators} Given a sufficiently good
function $w(r)$ on $\bbr_+$,  we define the following intertwining
operator \index{b}{Latin and Gothic!wf@$W f$} \bea\label{wxt} (W
f)(\tau)\equiv (W f)(\xi, t)&=& \intl_{\Ma} f(x)\,w( |\xi '
x-t|_m)\,dx\\\nonumber &=&\intl_{\Ma} f(x)\,w( |x-\tau|_m)\,dx,
\eea which assigns to a sufficiently good function $f(x)$ on $\Ma$
a function $(W f)(\tau)$ on the manifold $\Gr$ of matrix
$k$-planes $\tau=\tau(\xi,t)$, $(\xi,t)\in\cd$. The dual operator,
which maps a function $\vp(\tau)$ on $\Gr$ to the corresponding
function $(W^{\ast}\vp)(x)$ on $\Ma$, is defined by
\index{b}{Latin and Gothic!wfs@$W^{\ast}\vp$} \bea
(W^{\ast}\vp)(x)&=&\frac{1}{\sigk}\intl_{\vnk}\!\!d\xi\!\!\intl_{\Mt}\!\!\!\!\!\!\vp(\xi,t)\,w(
|\xi ' x-t|_m)\, dt\label{wdxt}\\\nonumber &=&\frac{1}{\sigk}
\intl_{\Gr} \vp(\tau)\, w( |x-\tau|_m)\,d\tau .\eea We recall that
$|x-\tau|_m$ denotes for the determinant of the matrix distance
between $x$ and $\tau$; see (\ref{detdist}). For $m=1$, operators
(\ref{wxt}) and (\ref{wdxt}) were introduced in \cite{Ru7}. The
following statement follows immediately from Lemma \ref{l5.2}.

\begin{lemma}
Operators $W$ and $W^{\ast}$ commute with the group $M(n,m)$ of
matrix motions.
\end{lemma}

\begin{lemma}
If $\tau=\tau(\xi,t)$, $ (\xi,t)\in\cd $, then \be\label{w-r} (W
f)(\tau)\equiv(W f)(\xi, t)=\intl_{\Mt}\hat f(\xi, z) \,w
(|t-z|_m)\, dz, \ee provided that either side of (\ref{w-r}) is
finite for $f$
 replaced by $|f|$.
\end{lemma}
\begin{proof}
Let  $ g_\xi \in SO(n)$ be a rotation satisfying $$
g_\xi\xi_0=\xi, \qquad \xi_0=\left[\begin{array} {c}  0 \\
I_{n-k} \end{array} \right] \in \vnk. $$ The change of variable
$x=g_\xi y$ in (\ref{wxt}) gives
$$
(W f)(\xi, t)= \intl_{\Ma} f(g_\xi y)\,w( |\xi_0 ' y-t|_m)\,dy,\qquad \xi_0= \left[\begin{array} {c}  0 \\
I_{n-k}
\end{array} \right]  \in V_{n, n-k}.
$$
By setting
 \[ y=\left[\begin{array} {c} \om \\ z
\end{array} \right], \qquad \om \in \frM_{k,m}, \qquad z \in
\frM_{n-k,m},\] we have \bea \nonumber(W f)(\xi, t)&=&
\intl_{\frM_{n-k,m}}w( |t-z|_m)\; dz\intl_{\frM_{k,m}}
f\left(g_\xi \left[\begin{array} {c} \om \\ z
\end{array} \right]\right)d\om\\\nonumber
&=&\intl_{\frM_{n-k,m}}\hat f(\xi, z)\,w( |t-z|_m)\; dz. \eea

\end{proof}

 \begin{lemma}
Let $\;1\le k\le n-m$, $\del=(n-k)/2-d$, $d=(m+1)/2$,

$$\tau=\tau(\xi,t),\qquad (\xi,t)\in\cd,\qquad u_0=\left[\begin{array} {c}  0 \\
I_{m}
\end{array} \right]  \! \in V_{n-k,m}.$$ Then \bea\label{wpm} (W f)(\xi,
t)&=&2^{-m}\sigma_{n-k,m}\intl_{\p}|r|^{\del}\,w(|r|^{1/2})\,\hat
f_r(\xi,t)\,dr,\\[14pt]
(W^{\ast}\vp)(x)&=&2^{-m}
\sigma_{n-k,m}\intl_{\p}|r|^{\del}\,w(|r|^{1/2})\,\check \vp_{
r}(x)\,dr .\label{dwpm}\eea
\end{lemma}
\begin{proof}

From (\ref{w-r}), by passing to the polar coordinates (see Lemma
\ref{l2.3}), we obtain

 \bea \nonumber(W f)(\xi, t)&=&\intl_{\Mt}\hat f(\xi, t+z)\, w
(|z|_m)\, dz\\\nonumber
&=&2^{-m}\intl_{\p}|r|^{\del}\,w(|r|^{1/2})\,dr \intl_{V_{n-k,m}}
\hat f(\xi, t+vr^{1/2})\,dv. \eea
 By (\ref{shr1}), this
coincides with (\ref{wpm}). Furthermore, (\ref{wdxt}) yields
\bea\nonumber
(W^{\ast}\vp)(x)&=&\frac{1}{\sigk}\intl_{\vnk}\!\!d\xi\!\!\intl_{\Mt}\!\!\!\!\!\!\vp(\xi,t+\xi
' x)\,w( |t|_m) \,dt\\\nonumber
&=&\frac{2^{-m}}{\sigk}\intl_{\vnk}\!\!d\xi\intl_{\p}|r|^{\del}\,w(|r|^{1/2})\,dr\!
\!\!\intl_{V_{n-k,m}}\!\!\!\vp(\xi, u r^{1/2}+\xi ' x) \,du \\[14pt]
&=&\frac{2^{-m}\sigma_{n-k,m}}{\sigk}\!\!
\intl_{\vnk}\!\!d\xi\intl_{\p}|r|^{\del}\,w(|r|^{1/2})\,dr\!
\!\intl_{O(n-k)}\!\!\vp(\xi, \theta u_0 r^{1/2}+\xi ' x)\,
d\theta. \nonumber\eea Now we change the order of integration and
replace $\xi$ by $\xi\theta '$. Since $\vp(\xi\theta',\theta
t)=\fc$, then
$$
(W^{\ast}\vp)(x)=\frac{2^{-m}\sigma_{n-k,m}}{\sigk}
\intl_{\p}|r|^{\del}\,w(|r|^{1/2})\,dr\!
\!\intl_{\vnk}\!\!\vp(\xi, u_0 r^{1/2}+\xi ' x)\, d\xi.
$$
By (\ref{7.1}), this gives (\ref{dwpm}).
\end{proof}

\section{The generalized Semyanistyi fractional integrals}

To avoid possible confusion, we shall discriminate between
operators acting on $\Ma$ and the similar operators on $\Mt$. As
before, the notation $I^\a$, $\Del$, $\DP$, and $\H$ will be used
for the Riesz potential, the Cayley-Laplace operator, the Cayley
operator, and the generalized Hilbert transform on $\Ma$. We write
$\tilde I^\a$,\index{a}{Riesz potential! on $\Mt$, $\tilde I^\a
f$} \index{b}{Latin and Gothic!iat@$\tilde I^\a\vp$} $\tilde
\Del$, \index{a}{Cayley-Laplace operator!on $\Mt$, $\tilde\Del$}
$\tilde \DP$, \index{a}{Cayley operator!on $\Mt$, $\tilde\DP$} and
$\tilde \H$ \index{a}{Generalized Hilbert transform!on
$\Mt$,$\tilde\H f$}\index{b}{Latin and Gothic!h@$\tilde\H f$}
\index{b}{Latin and Gothic!da@$\tilde\DP$}
\index{b}{Greek!dela@$\tilde\Del$}for the similar operators on
$\Mt$. These will be applied to functions
 $\rf$ and $\fc$ in the $t$-variable. We assume $1\le k\leq n-m$, and
denote by $\s(\frT)$ \index{b}{Latin and Gothic!st@$\s(\frT)$} the
space of functions $\fc$ which are infinitely differentiable in
the $\xi$-variable and belong to the Schwartz space $\Mt$ in the
$t$-variable uniformly in $\xi \in \vnk$.

 Consider the following operators
 \index{b}{Latin and Gothic!pf@$P^{\a} f$} \index{b}{Latin and Gothic!pfd@$\pd f$}
\be\label{p-rr} \index{a}{Generalized Semyanistyi fractional
integrals,  $P^{\a}$, $\pd$} P^{\a}f=\tilde I^\a\hat f, \qquad
\pd\vp=(\Tilde I^\a\vp)^\vee , \ee where $f \in \s(\Ma)$,  $\vp
\in \s(\frT)$, and \be\label{a-p}\a\in\bbc, \quad \a\neq n-k-m+1,
n-k-m+2, \ldots\; .\ee The right-hand sides of equalities in
(\ref{p-rr}) absolutely converge for $Re\,\a >m-1$ and are
understood in the sense of analytic continuation for other values
of $\a$. Below we prove a series of lemmas giving explicit
representation of operators (\ref{p-rr}) for different values of
$\a$.
\begin{lemma}
Let $f \in \s(\Ma)$,  $\vp \in \s(\frT)$. If $$Re\,\a
>m-1,\qquad \a\neq n-k-m+1,\; n-k-m+2, \ldots\; ,$$ then $ P^{\a}$ and
$\pd$ are  intertwining operators  of the form  (\ref{wxt}) and
(\ref{wdxt}), respectively. Namely, \bea\label{p}
 \qquad \qquad (P^\a f)(\xi, t)&=&
\frac{1}{\g_{n-k,m}(\a)}\intl_{\Ma} f(x)\,|\xi ' x-t|_m^{\a+k-n}\,dx,\\
(\pd\vp)(x)\!\!&=&\!\!\frac{1}{
\g_{n-k,m}(\a)}\intl_{\vnk}\!\!d_\ast
\xi\!\!\intl_{\Mt}\!\!\!\!\!\!\vp(\xi,t)\,|\xi' x-t|_m^{\a+k-n}\,
dt\label{pd}, \eea where $d_\ast \xi=\sigk^{-1} \, d\xi$ is the
normalized measure on $\vnk$ and $\g_{n-k,m}(\a)$ is the
normalized constant for the Riesz potential on $\Mt$; cf.
(\ref{gam}). \index{b}{Latin and Gothic!dst@$d_\ast \xi$}
\end{lemma}

\begin{proof}
The formula (\ref{pd})  follows immediately from  the definitions
(\ref{p-rr}), (\ref{4.2}), and (\ref{rie}). Furthermore, by
(\ref{p-rr}), (\ref{4.9}), and (\ref{rie}), we have \bea
\nonumber(P^\a f)(\xi, t)&=& \frac{1}{
\g_{n-k,m}(\a)}\!\!\!\!\!\intl_{\Mt}\!\!\!\!\! |y|_m^{\a+k-n} dy
\!\!\!\!\!\intl_{\Mkm}\!\!\!\!\! f\left(g_\xi \left[\begin{array}
{c} \omega
\\t-y
\end{array} \right]\right)d\omega \\\nonumber
&=& \frac{1}{\g_{n-k,m}(\a)}\intl_{\Ma} f(x)\,|\xi '
x-t|_m^{\a+k-n}\, dx. \eea This proves (\ref{p}).
\end{proof}

\begin{lemma} Let $f \in \s(\Ma)$,  $\vp \in \s(\frT)$, $1\le k\le n-m$. If
 $\ell$ is a positive integer  so that $\ell \le n-k-m$ (cf.
(\ref{a-p})), then    \bea\label{pl} \index{a}{Generalized
Semyanistyi fractional integrals,  $P^{\a}$, $\pd$}
\qquad(P^{\ell} f)(\xi,t)&=&c_\ell \intl_{\frM_{k+\ell,m}}d z
\intl_{O(n-k)}f\left(\xi t - g_\xi\left[\begin{array} {cc} I_k &
0\\ 0 & \g
\end{array}
\right]\left[\begin{array} {c} z \\
0
\end{array} \right]\right)d\g,\\
\label{pdl} (\pdl\vp)(x)&=&c_\ell
\intl_{\vnk}d_\ast\xi\intl_{\frM_{\ell,m}}dz
\intl_{O(n-k)}\vp\left(\xi,\;\xi 'x-\g\left[\begin{array} {c} z \\
0
\end{array} \right]\right)d\g,
\eea where \be\label{c-p} c_\ell=2^{-\ell m} \,\pi^{-\ell m/2} \,
\Gam_m\Big ( \frac{n-k-\ell}{2}\Big ) / \Gam_m\Big (
\frac{n-k}{2}\Big ). \ee Moreover, \be\label{p0}
(P^{0}f)(\xi,t)=\rf, \qquad (\pdz\vp)(x)=\df.  \ee
\end{lemma}
\begin{proof}
By  (\ref{p-rr}),  analytic continuation of $P^{\a} f$ and
$\pd\vp$ reduces to that of the Riesz potential on $\Mt$. One can
readily see that $\rf$, defined by (\ref{4.9}),  is a Schwartz
function in the $t$-variable. Thus, \be P^{\ell}f=\tilde I^\ell
\hat f, \qquad \pdl\vp=(\Tilde I^\ell\vp)^\vee . \ee  Now
(\ref{pdl}) follows  from (\ref{des}), and (\ref{p0}) is a
consequence of (\ref{I0}). Furthermore, by (\ref{des}) and
(\ref{4.9}),

 \bea \nonumber(P^{\ell} f)(\xi,t)&=&c_\ell \intl_{\frM_{\ell,m}}d
z
\intl_{O(n-k)}\hat f\left(\xi, t+\g\left[\begin{array} {c} z \\
0
\end{array} \right]\right)d\g \\[14pt]
&=&c_\ell \intl_{\frM_{\ell,m}}dz
\intl_{O(n-k)}d\g\intl_{\frM_{k,m}} f\left(g_\xi
\left[\begin{array} {c} \om
\\t+\g\left[\begin{array} {c} z \\
0
\end{array} \right]\end{array}\right]\right) d\om .\nonumber
\eea Since $$ g_\xi \left[\begin{array} {c} \om
\\t+\g\left[\begin{array} {c} z \\
0
\end{array} \right]\end{array}\right]=\xi t + g_\xi\left[\begin{array} {cc} I_k &
0\\ 0 & \g
\end{array}
\right]\left[\begin{array} {c}\om \\z \\
0
\end{array} \right],
$$
then (\ref{pl}) follows if we change the notation
$\left[\begin{array} {c}\om
\\z
\end{array} \right]\to z$.
\end{proof}

\begin{remark}
We call $P^\a f$ and $\pd\vp$ {\it the generalized Semyanistyi
fractional integrals}; see Remark \ref{r5.14} for comments.
  Formulas (\ref{p})--(\ref{p0}) can serve as definitions of
$P^\a f$ and $\pd \vp$ if $f$ and $\vp$ are arbitrary locally
integrable functions so that the corresponding integrals converge.
\end{remark}

\begin{lemma} \label{lpm1} Let $f \in \s(\Ma)$,  $ \; \vp \in \s(\frT)$.
If $\ell=1,2,\dots\;$, then \be\label{pm1} \index{a}{Generalized
Semyanistyi fractional integrals,  $P^{\a}$, $\pd$} \qquad
(P^{-2\ell}f)(\xi,t)=(-1)^{m\ell} \tilde \Delta^\ell\rf,\ee and
\be\label{pdm1} (\stackrel{*}{P}\!{}^{-2\ell}\vp) (x)
=(-1)^{m\ell}[\tilde\Delta^{\ell}\vp]^\vee (x)  =
(-1)^{m\ell}\intl_{\vnk} \left .\tilde\Delta^{\ell}\fc \right
|_{t=\xi'x} \,
 \, d_\ast\xi.  \ee
\end{lemma}
\begin{proof} The statement follows from (\ref{p-rr}) and
(\ref{Dkf}).
\end{proof}

\begin{lemma}\label{lpm2} Let $f \in \s(\Ma)$,  $\vp \in \s(\frT)$,  $ \;
\ell=1,2,\dots\;$. We denote \be\label{c1c2}
c_1=\frac{(-1)^{m\ell}\,
 \Gam ((n-k-m)/2)}{ 2^{m+1} \, \pi^{(m+n-k)/2}}, \qquad
 c_2=\frac{(-1)^{m(\ell+1)}\, \Gam _m((m+1)/2)}{\pi^{m^2/2}}. \ee

\noindent {\rm(i)} If   $1\le k\le n-m$ and  $F_\xi
(t)=\tilde\Delta^\ell\hat f(\xi,t)$  then \bea\label{pm1}
(P^{1-2\ell}f)(\xi,t)&=& (-1)^{m\ell} (\tilde
I^1 F_\xi)(t)\\[14pt]\nonumber
&=&c_1 \, \intl_{S^{n-k-1}}dv\intl_{\bbr^m} F_\xi(t-vy')\, dy
\nonumber\eea and \bea \label{pm11}\qquad\qquad\qquad\quad
(\stackrel{*}{P}\!{}^{1-2\ell}\vp) (x)
&=&(-1)^{m\ell}[\tilde I^1\tilde\Delta^{\ell}\vp]^\vee (x)\\[14pt]\nonumber
&=& c_1 \! \intl_{\vnk}\!\!\! d_\ast\xi\!\!\!
\intl_{S^{n-k-1}}\!\!\!dv\intl_{\bbr^m} \left .
\tilde\Delta^{\ell}\vp(\xi,t )\right |_{t=\xi'x-vy'}\, \, dy .
\nonumber \eea

\noindent {\rm(ii)} If   $k=n-m$, then
 \be\label{ppm1}(P^{1-2\ell}f)(\xi,t)=c_2 \, (\tilde \H
\tilde\DP^{2\ell-1}
 \hat f(\xi, \cdot))(t)\ee
 and
 \bea
 \label{ppm11} (\stackrel{*}{P}\!{}^{1-2\ell}\vp) (x)
&=&c_2 \, (\tilde \H \tilde\DP^{2\ell-1}\vp(\xi, \cdot ))^\vee (x)
\\&=&c_2 \,
\intl_{\vnk} (\tilde \H \tilde\DP^{2\ell-1}\vp(\xi, \cdot
))(\xi'x) \,  d_\ast\xi, \nonumber \eea $\tilde \H$ being the
generalized Hilbert transform; cf. (\ref{hilb}).

\end{lemma}
\begin{proof} (i) follows from (\ref{p-rr}) and
(\ref{Dkk}); (ii) is a consequence of  (\ref{p-rr}) and
(\ref{Dkkk}).
\end{proof}

The following statement is  the main results of this chapter.

\begin{theorem}\label{bas}
Let $1\le k\le n-m$, $\a\in\bbc$; $\a\neq n-k-m+1, n-k-m+2,
\ldots\;$.

 \noindent{\rm(i)} \ If $f\in\s(\S_m)$ then
\be\label{gfu} \index{a}{Fuglede formula!generalized} (\pd\hat
f)(x)= c_{n,k,m} (I^{\a+k} f)(x)\ee (the generalized Fuglede
formula), where
 \be\label{cnkm}
 c_{n,k,m}=2^{km}\pi^{km/2}\gm\left(\frac{n}{2}\right)/\gm\left(\frac{n-k}{2}\right).
 \ee

 \noindent{\rm(ii)}\ Let $\ell_0=\min\{m-1, n-k-m \}$, $$
\A=\{0, 1, 2, \ldots, \ell_0\} \cup \{\a  :  Re\,\a \!
> \! m \! - \! 1 ; \; \a  \! \neq \!  n-k \! - \! m +\! 1,  n-k \! - \! m +\! 2,\dots
\;\}.
$$
Suppose that  $f\in L_{loc}^1(\Ma)$ and   $\a\in\A$. Then equality
(\ref{gfu}) holds provided that the  Riesz potential $(I^{\a+k}
f)(x)$ is finite for $f$  replaced by $|f|$ {\rm(}e.g., for $f\in
L^p$, $1\leq p< n/{\rm(}Re\,\a+k+m-1${\rm)}{\rm)}.
\end{theorem}
\begin{proof}
{\rm(i)} Let $f\in\s(\S_m)$. We make use of the equality
(\ref{5.8}) with $\lam=\a+k$ and $Re\,\a
>m-1$. This gives \be\label{5.8.1}
 \begin{array}{ll}
{\displaystyle \frac{1}{\sigk} \intl_{\vnk}d\xi\intl_{\Mt}\rf\,
|t|_m ^{\a+k-n}\,dt} \\[30pt]
={\displaystyle
\frac{\Gamma_m(n/2)\,\Gamma_m(\a/2)}{\Gamma_m((\a+k)/2)\,\Gamma_m((n-k)/2)}\intl_{\Ma}
f(y)\,|y|_m^{\a+k-n} \,dy}.
\end{array}
\ee Replacing $f(y)$ by the shifted function $f_x(y)=f(x+y)$ and
taking into account (\ref{4.4}), we get \be\label{5.8.2}
 \begin{array}{ll}
{\displaystyle \frac{1}{\sigk} \intl_{\vnk}d\xi\intl_{\Mt}\hat
f(\xi, \xi 'x+t)\,
|t|_m^{\a+k-n}\,dt} \\[30pt]
={\displaystyle
\frac{\Gamma_m(n/2)\,\Gamma_m(\a/2)}{\Gamma_m((\a+k)/2)\,\Gamma_m((n-k)/2)}\intl_{\Ma}
f(x+y)\,|y|_m^{\a+k-n}\, dy},
\end{array}
\ee cf. (\ref{pd})  and (\ref{rie}).  Hence, (\ref{gfu}) follows
when $Re\,\a
>m-1$  with the constant
\bea\nonumber
c_{n,k,m}&=&\frac{\Gamma_m(n/2)\,\Gamma_m(\a/2)\,\g_{n,m}(\a+k)}{\Gamma_m((\a+k)/2)\,
\Gamma_m((n-k)/2)\,\g_{n-k,m}(\a)}\\[14pt]
&\stackrel{\rm
(\ref{gam})}{=}&2^{km}\pi^{km/2}\gm\left(\frac{n}{2}\right)/\gm\left(\frac{n-k}{2}\right).\nonumber
\eea By analytic continuation, it is true for all $\a\in\bbc$,
$\a\neq n-k-m+1, n-k-m+2, \ldots\;$.

{\rm(ii)}  Suppose $f\in L_{loc}^1(\Ma)$. For $Re\,\a
>m-1$, (\ref{gfu}) follows from (\ref{5.8.2}) by taking into account that (\ref{5.8.2}) was derived from
(\ref{5.8}), and the latter is also true for locally integrable
functions.
 For $\a=\ell$, $\ell=1,2,\dots
\ell_0$, we have \bea\nonumber (\pdl \hat
f)(x)\!\!\!&=&\!\!\!\frac{c_\ell}{\sigk }
\intl_{\vnk}\!\!\!d\xi\intl_{\frM_{\ell,m}}\!\!\!dz
\intl_{O(n-k)}\!\!\!d\g\intl_{\frM_{k,m}} \!\!\!f\left(x-g_\xi
\left[\begin{array} {c} \om \\ \g\left[\begin{array} {c} z
\\
0
\end{array} \right]
\end{array} \right]\right)d\om \\[14pt]
&=&c_\ell\intl_{O(n)}d\b\intl_{\frM_{\ell,m}}dz
\intl_{O(n-k)}d\g\intl_{\frM_{k,m}} f\left(x-\b
\left[\begin{array} {c} \om \\ \g\left[\begin{array} {c} z \\
0
\end{array} \right]
\end{array} \right]\right)d\om,\nonumber
\eea where $c_\ell$ is the constant (\ref{c-p}). We  write
$$
\left[\begin{array} {c} \om \\ \g\left[\begin{array} {c} z \\
0
\end{array} \right]\end{array} \right]=\left[\begin{array} {cc} I_k &
0\\ 0 & \g
\end{array}
\right]\left[\begin{array} {c}\om \\z \\
0
\end{array} \right].\\[14pt]
$$
Then the change of variables  $\b\left[\begin{array} {cc} I_k & 0\\
0 & \g
\end{array}
\right]\to\b$ gives

\bea\nonumber (\pdl \hat f)(x)&=&c_\ell\intl_{\frM_{\ell,m}}dz
\intl_{\frM_{k,m}}d\om\intl_{O(n)}f\left(x-\b
\left[\begin{array} {c}\om \\z \\
0
\end{array} \right]\right)d\b\\[14pt]\nonumber
&=&c_\ell\intl_{\frM_{\ell+k,m}}dy \intl_{O(n)}f\left(x-\b
\left[\begin{array} {c} y\\
0
\end{array} \right]\right)d\b\\ [14pt]\nonumber
&=& c_{n,k,m} (I^{\ell+k} f)(x),\eea where \bea\nonumber c_{n,k,m}
&=& c_\ell \, 2^{(\ell+k) m} \,\pi^{(\ell+k)m/2} \,\Gam_m\Big (
\frac{n}{2}\Big )  /\Gam_m\Big ( \frac{n-\ell-k}{2}\Big )
\\[14pt]
&=&2^{km}\pi^{km/2}\gm\left(\frac{n}{2}\right)/\gm\left(\frac{n-k}{2}\right).\nonumber
\eea
\end{proof}

\begin{remark}\label{r5.14} Some comments are in order. The idea to study
Radon transforms as members of the corresponding analytic families
goes back to Semyanistyi  \cite{Se} who considered the hyperplane
Radon transform on $\rn$ (the case $m=0, \; k=n-1$). This approach
was extended by
 Rubin to Radon transforms of different kinds, see, e.g., \cite{Ru11}. In the rank-one case $m=1$
the formula (\ref{gfu}) is due to Fuglede \cite{Fu} for $\a=0$ and
to Rubin \cite{Ru11}, \cite{Ru7} for any $\a$. In the higher  rank
case $m>1$, it was established for $\a=0$ in \cite{OR} (for
sufficiently good functions) and justified for $f\in L^p$ in
\cite{Ru10}.

The formula (\ref{gfu}) has the same nature as the classical
 decomposition of  distributions  in plane waves. To see
this, one should formally set $x=0$ in (\ref{gfu}) and regard this
equality  in the framework of the theory of distributions. The
idea of  decomposition of a function in integral of plane waves
amounts to  Radon \cite{R} and John \cite{Jo}, and has proved to
be very fruitful in PDE. For the distribution $|x|^\lam$ on $\rn$
an account of this theory is presented in  \cite[Section
3.10]{GSh1}. This method  was generalized by  Petrov \cite{P2} for
the matrix case for $k=n-m$, and briefly outlined by  Shibasov
\cite{Sh2} for $k<n-m$. Our way of thinking differs from that in
\cite{Sh2}, \cite{P2}, and \cite{GSh1}, and the result is more
general because we allow $f$ to be a ``rough'' function.
\end{remark}

For the sake of completeness, we also present the following
statement which follows from Theorem \ref{bas} and the semigroup
property of Riesz potentials. Concerning this property, see
\cite{Kh} and \cite{Ru10}.

\begin{theorem}
Let $$Re\, \a >m-1, \qquad Re\, \b>m-1, \qquad Re\,(\a+
\b)<n-k-m+1.$$ If the integral $ I^{\a+\b+k}f$  absolutely
converges then \be\label{pab} \pd P^\b f= c_{n,k,m}I^{\a+\b+k}f,
\ee $c_{n,k,m}$ being the constant (\ref{cnkm}).
\end{theorem}
\begin{proof} By (\ref{gfu}) and (\ref{p-rr}), we have
\[ c_{n,k,m}I^{\a+\b+k}f=\stackrel{*}{P}{}^{\a+\b} \hat f=
(\tilde I^{\a+\b} \hat f)^\vee= (\tilde I^{\a}\tilde I^{\b}\hat
f)^\vee=\stackrel{*}{P}{}^\a P^\b f. \]
\end{proof}

\chapter{Inversion of the Radon transform}\label{s6}

\section{The radial case} Theorems \ref{t5.5}, \ref{t5.1},
  and the inversion formula (\ref{inm})  for
the G{\aa}rding-Gindikin fractional integrals imply the following
result for the
 Radon transform  of radial functions.

\begin{theorem}\index{a}{Inversion formula!for the Radon
transform!of radial functions} Suppose that $f(x)\equiv\f0 (r)$, $
\; x\in\Ma$, $ \; r=x'x$, and let \be\label{6.p} f\in\lp, \qquad
1\leq p<\frac{n+m-1}{k+m-1}.\ee Then the Radon transform $\rf$ is
well defined by (\ref{4.9}) for almost all $(\xi, t) \in \vnk
\times \Mt$ and represents a radial function, namely,
\be\label{rad}\rf=\pi^{km/2} (I_{-}^{k/2}f_0)(s)=\vp_0(s),\qquad
s=t't\in\cpm .\ee If $1\le k\le n-m$, then   $f_0$ can be
recovered from $\vp_0$ by the formula \be\label{6.4}
f_0(r)=\pi^{-km/2}(D _-^{k/2} \vp_0)(r),\ee
 where  $D_-^{k/2}$ is defined in the sense of $
\D'(\p)$-distributions by (\ref{dmin}).
\end{theorem}

\begin{remark} The condition (\ref{6.p}) can be replaced by
 the weaker one. Indeed, owing to Theorem \ref{l3.14}, it suffices to  assume that
 \be\label{6.R} \intl_{R}^\infty
|r|^{(k-m-1)/2}\,|f_0(r)| \, dr<\infty \qquad \text{for all $ \; R
\in \p$}. \ee  Furthermore, we know that by Theorem \ref{inj},
 the assumption $ k\leq n-m$ is necessary for injectivity of the Radon transform.
 One might expect that once we restrict to   radial functions, then this assumption can be reduced.
  However, it is not so, because the
function $\psi$ in the proof of Theorem \ref{inj} can be chosen to
be radial. Note that if $ k>n- m$, then the exterior variable $s$
in (\ref{rad}) ranges on the boundary of the cone, and we ``lose
the dimension".
\end{remark}

In the same manner, Theorem \ref{t5.7} and Lemma \ref {l-inv1}
allow us to obtain  an inversion formula for the dual Radon
transform.

\begin{theorem}
Let $\fc \equiv\vp_0 (s)$, $(\xi, t) \in \vnk \times \Mt$,
$s=t't$. We assume $1\le k\le n-m$ and denote
\[ \Phi_0(s)=|s|^\del \vp_0 (s), \qquad \del=(n-k)/2-d,  \qquad d=(m+1)/2.\]
If $\Phi_0(s) \in L^1_{loc}(\p)$ then the dual Radon transform
$\check\vp (x)$ is well defined  by (\ref{4.2}) for almost all $x
\in \Ma$ and represents the radial function \be\label{radd} \df=c
|r|^{d-n/2}(I_{+}^{k/2} \Phi_0)(r)=f_0(r), \quad r=x'x\in\cpm,\ee
  $c=\pi^{km/2}\sigma_{n-k,m}/\sigma_{n,m}$.
 The function $\vp_0$ can
be recovered from $f_0$ by the formula \be\label{6.4.1}
\vp_0(s)=c^{-1}|s|^{-\del}(D_{+}^j I_+^{j-k/2}F_0)(s), \qquad
F_0(r)=|r|^{n/2-d}f_0(r),\ee for any integer $j\ge k/2$.
 If $k$
is even, then \be\label{6.4.2}
\vp_0(s)=c^{-1}|s|^{-\del}(D_{+}^{k/2}F_0)(s). \ee
 The differential operator $D_{+}$ in (\ref{6.4.1}) and (\ref{6.4.2}) is
  understood in the sense of $\D'(\p)$-distributions (cf. (\ref{inp1}) and  (\ref{inp2})).
\end{theorem}

\section{The  method of mean value
operators } As it was  mentioned in Introduction, after the 1917
paper by Radon \cite{R}, the following two inversion methods have
been  customarily  used for reconstruction of functions from their
integrals over affine planes. These are  {\it the method of mean
value operators}, and  {\it the method of Riesz potentials}. The
third classical method presented  in \cite{GSh1} is  {\it the
method of plane waves}. Below we focus on the first method and
extend it to functions of matrix argument.

\begin{definition}  Given a function $f(x)$ on $\Ma$,  we define

 \be\label{6.2}
(M_r f)(x)= \frac{1}{\sigma_{n,m}}\intl_{\vnm}
f(vr^{1/2}+x)dv,\qquad r\in \p.
  \ee
  This   is a matrix generalization of the usual
  spherical mean \index{a}{Spherical mean (matrix generalization)}\index{b}{Latin and Gothic!mr@$\M_r f$} on $\bbr^n$.
  \end{definition}
We say that
  $r\in \p$ tends to zero when $\tr(r)\to 0$.

  \begin{lemma}\label{l7.6} If $f \in L^p(\Ma), \; 1 \le p < \infty$, then
 \be\label{6.9} \lim\limits_{r\to 0}(M_r f)(x)=f(x)\ee
in the $L^p$-norm. If $f \in C_0 (\Ma)$, i.e., $f(x) \to 0$ as
$||x||=(\tr (x'x))^{1/2} \to \infty$, this limit   is uniform on
 $\Ma$.
 \end{lemma}
 \begin{proof}
By the generalized Minkowski inequality, $$ \| M_r f -f\|_p\leq
\frac{1}{\sigma_{n,m}}\intl_{\vnm}\|
f(vr^{1/2}+\cdot)-f(\cdot)\|_p\; dv. $$ Since the integrand does
not exceed $2\|f\|_p$, by the Lebesgue theorem on dominated
convergence, one can pass to the limit under the sign of
integration.  Let $y=vr^{1/2}\in\Ma$. We represent the
corresponding  $nm$-vector $\bar y=(y_{1,1},\dots, y_{n,m})$  in
polar coordinates
$$\bar y=\theta \rho,\qquad \rho=\|y\|=(\tr(y'y))^{1/2},\qquad
\theta\in S^{nm-1},$$
 so that
$$
\lim\limits_{r\to 0}\|
f(vr^{1/2}+\cdot)-f(\cdot)\|_p^p=\lim\limits_{\rho\to
0}\intl_{\bbr^{nm}} | f(\theta \rho +\bar x)-f(\bar x)|^p d\bar
x=0.
$$
This gives (\ref{6.9}) in the $L^p$-norm. The proof of the second
part of the statement is the same.
\end{proof}

The following lemma which combines  the Radon transform and the
shifted dual Radon transform (\ref{dualsh}) is the core of the
method. It reduces the inversion problem for the Radon transform
to the case of radial functions. Thus ``dimension of the problem"
becomes essentially smaller.

  \begin{lemma}\label{l6.6}
  For  fixed  $x \in \Ma$, let $F_x(r)=(M_r f)(x), \; r\in \p.$ Then
  \be\label{6.3}
(\hat f)_s^{\vee}(x)=\pi^{km/2} (I_{-}^{k/2}F_x)(s), \qquad
s\in\p, \ee provided that either side of this equality exists in
the Lebesgue sense.
\end{lemma}

\begin{proof}
We denote $\fc=\rf$, $f_x(y)=f(x+y)$. Let $z\in\Mt$ be a matrix at
distance $s^{1/2}$ from the origin. By (\ref{6.1}) and
(\ref{4.4}), \bea
 (\hat f)_s^{\vee}(x)&=&
\frac{1}{\sigk}\intl_{\vnk} \hat f_x (\xi, z)d\xi \nonumber \\
&=&\intl_{SO(n)} \hat f_x (\gam\xi, z)d\gam. \nonumber \eea
Interchanging the order of the Radon transform and integration
over $SO(n)$, and using (\ref{4.23}), we get

$$
\check \vp_s (x)=\intl_{\Mkm}d\omega \intl_{SO(n)} f_x\left(\g
\left[\begin{array} {c} \omega \\z
\end{array} \right]\right)d\g.
$$
Hence,  $z \to \check \vp_s (x)$ is  the Radon transform
  of the radial function
$$\tilde f_x(y)=\intl_{SO(n)} f_x(\g y)d\g= \intl_{SO(n)} f(x+\g y)d\g,
\qquad y\in \Ma.$$ In other words, \be\label{con} \check \vp_s
(x)=(\tilde f_x)^\wedge (\xi,z), \ee where
$$\tilde f_x(y)=\frac{1}{\sigma_{n,m}}\intl_{\vnm}
f(x+vr^{1/2})dv=F_x(r), \qquad  r=y'y.$$ Owing to (\ref{con}), the
desired equality (\ref{6.3}) follows from the representation
(\ref{4.10}) for the Radon transform of a radial function.
\end{proof}

\begin{corollary}
Let $f\in\lp$, $1\leq p<(n+m-1)/(k+m-1)$. Then (\ref{6.3}) holds
 for almost all $x\in\Ma$. In particular, (\ref{6.3}) holds for
 any
continuous function $f$  satisfying \[
f(x)=O(|I_m+x'x|^{-\lam/2}), \qquad \lam>k+m-1. \]
\end{corollary}
\begin{proof} It suffices to show that  the right side of (\ref{6.3}) is finite for
almost all $x\in\Ma$ whenever  $f\in\lp$, $1\leq
p<(n+m-1)/(k+m-1)$. By Theorem \ref{l3.14}, the integral
$(I_{-}^{k/2}F_x)(s)$ is well defined provided
\be\label{fin}I=\intl_{R}^\infty |r|^{k/2-d}|F_x(r)| \, dr<\infty
\qquad \text{for all $ \; R \in \p$}. \ee The proof of (\ref{fin})
is simple. Indeed, taking into account (\ref{6.2}) and using Lemma
\ref{l2.3},
 we obtain \bea\nonumber
I&\leq&\frac{1}{\sigma_{n,m}}\intl_{R}^\infty
|r|^{k/2-d}\,dr \intl_{\vnm} |f(vr^{1/2}+x)|\,dv\\
&=& \frac{2^m}{\sig_{n,m}}\intl_{\{y\in\Ma\,:\,y'y>R\}}
|y'y|^{(k-n)/2} |f(x+y)|\,dy.\nonumber
 \eea
By H\"older's inequality, $I\leq A \nf$, where the constant $A$ is
the same as in the proof of Theorem \ref{t5.1}. If
$p<(n+m-1)/(k+m-1)$ then $A<\infty$, and we are done.
\end{proof}
Now we can prove the  main result.
\begin{theorem}
Let $1\le k\leq n-m$, \be f\in\lp, \qquad 1\leq
p<\frac{n+m-1}{k+m-1}.\ee Then the Radon transform
$\vp(\xi,t)=\rf$ is well defined by (\ref{4.9}) for almost all
$(\xi, t) \in \vnk \times \Mt$ and can be inverted by the formula
\be\label{6.7}\index{a}{Inversion formula!for the Radon
transform}f(x)=\pi^{-km/2}\lim\limits_{r\to 0}^{(L^p)}(D_-^{k/2}
\Phi_x )(r), \qquad \Phi_x (s)=\check \vp_
 s (x),\ee where  $D_-^{k/2}$ is defined in the sense of $
\D'(\p)$-distributions by (\ref{dmin}).  If $f$ is a continuous
function satisfying \be \label{dec} f(x)=O(|I_m+x'x|^{-\lam/2}),
\qquad \lam>k+m-1, \ee  then the limit in (\ref{6.7}) can be
treated in the sup-norm.
\end{theorem}

\begin{proof}
By Lemmas \ref{l6.6}  and  \ref{ligm}, one can recover
$F_x(r)=(M_r f)(x)$ and get
$$ (M_r
f)(x)=\pi^{-km/2}(D_-^{k/2}\Phi_x)(r), \qquad r \in \p. $$ Now the
result follows by   Lemma \ref{l7.6}.
\end{proof}

\section{The  method of Riesz potentials}\label{mrp}
The second traditional inversion method for  the Radon transform
reduces the problem to inversion of  the  Riesz potentials. This
method relies on the formula (\ref{gfu}) where $\a$ is in our
disposal. The case $\a=0$ corresponds to the Fuglede formula
\be\label{fu}\index{a}{Fuglede formula} (\hat f)^{\vee} (x) \! =
\! c_{n,k,m} (I^k f)(x), \ee which together with   Theorem
\ref{ap} implies the following.
\begin{theorem}\label{t8.10} Let $f \in L^p (\Ma), \; 1\leq
p<n/(k+m-1)$. Then the Radon transform $\vp=\hat f$ is well
defined, and $f$ can be recovered from $\vp$ in the sense of
$\Phi'$-distributions by the formula \be\label{inv-rr}
\index{a}{Inversion formula!for the Radon transform}
c_{n,k,m}(f,\phi)=(\check \vp, I^{-k}\phi), \ee $$ \phi \in
\Phi,\qquad
c_{n,k,m}=2^{km}\pi^{km/2}\gm\left(\frac{n}{2}\right)/\gm\left(\frac{n-k}{2}\right),
$$
where the operator $I^{-k}$ is defined by
$$(I^{-k}\phi)(x)=(\F^{-1}|y|_m^{k}\F\phi)(x).$$ In particular, for
$k$ even,  \be  c_{n,k,m}(f,\phi)=(-1)^{mk/2}(\check \vp,
\Del^{k/2}\phi), \ee
 $\Del$ being the Cayley-Laplace operator (\ref{K-L}).

\end{theorem}

\begin{remark} For $k$ odd,  the Radon transform  can be
inverted under more restrictive assumptions as follows. Let
$k<n-m$. We choose $\a=1$ in (\ref{gfu}) and get \be
 \tilde I^1\hat f = c_{n,k,m}I^{k+1}f. \ee If $f \in
L^p (\Ma), \; 1\leq p<n/(k+m)$, then $f$ can be recovered from
$\vp=\hat f$  by the formula \be
c_{n,k,m}(f,\phi)=(-1)^{m(k+1)/2}(\tilde I^1 \vp,
\Del^{(k+1)/2}\phi), \qquad \phi \in \Phi. \ee
\end{remark}

\section{Decomposition in plane waves}

This method  was developed in \cite{P2} for the case $k=n-m$, and
outlined in \cite{Sh2} for $1\le k\le n-m$. It is based on
decomposition of distributions  in (matrix) plane waves. As we
have already noted in Remark \ref{r5.14},  our formula (\ref{gfu})
has the same nature. Together  with  (\ref{pdm1}), (\ref{pm11}),
and (\ref{ppm11}), it enables us to invert the Radon transform of
functions $f\in\s(\Ma)$.

\begin{theorem}\label{t6.13}
Let $1\le k\le n-m$, $f\in\s(\S_m)$.  The Radon transform
$\fc=\rf$ can be inverted by the following  formulas.
\index{a}{Inversion formula!for the Radon transform}

 \noindent{\rm(i)} \ For $k$ even,
\be \label{pw1}\index{a}{Inversion formula!for the Radon
transform} f(x)=(-1)^{mk/2}c_{n,k,m}^{-1}\intl_{\vnk} \left
.\tilde\Delta^{k/2}\fc \right |_{t=\xi'x} \,
 \, d_\ast\xi \;,\ee
$$
 c_{n,k,m}=2^{km}\pi^{km/2}\gm\left(\frac{n}{2}\right)/\gm\left(\frac{n-k}{2}\right)
 $$
 (we recall that $ d_\ast\xi$ stands for the invariant measure on
 $\vnk$ of total mass 1).

 \noindent{\rm(ii)} \ For $k$ odd and $k<n-m$, \be \label{pw2} f(x)=
c_1 \! \intl_{\vnk}\!\!\! d_\ast\xi\!\!\!
\intl_{S^{n-k-1}}\!\!\!dv\intl_{\bbr^m} \left .
\tilde\Delta^{(k+1)/2}\vp(\xi, t )\right |_{t=\xi'x-vy'}\, \, dy
,\ee
$$
c_1=(-1)^{m(k+1)/2} 2^{-km-m-1}\pi^{(k-n)/2-m(k/2+1)}
\Gam_{m+1}\left(\frac{n-k}{2}\right)/\gm\left(\frac{n}{2}\right).
$$

 \noindent{\rm(iii)} \  \index{a}{Inversion formula!for the Radon transform} For $k$
odd and $k=n-m$, \be \label{pw3} f(x)= c_2\intl_{\vnk} (\tilde \H
\tilde\DP^{k}\vp(\xi, \cdot ))(\xi'x) \,  d_\ast\xi,\ee

$$
c_2=(-1)^{m(k+3)/2} 2^{-km}\pi^{-m(k+m)/2}
\Gam_{m}\left(\frac{m}{2}\right)\gm\left(\frac{m+1}{2}\right)/\gm\left(\frac{n}{2}\right).
$$

\end{theorem}
\begin{proof}
We write (\ref{gfu}) with $\a=-k$ so that
  \be\label{gfu2} f(x)=c_{n,k,m}^{-1}(\stackrel{*}{P}\!{}^{-k}\vp)(x),\ee
where $\stackrel{*}{P}\!{}^{-k}$ is the operator (\ref{p-rr}). Now
it remains to apply formulas (\ref{pdm1}), (\ref{pm11}), and
(\ref{ppm11}).

\end{proof}

Formulas (\ref{pw1})--(\ref{pw3}) differ from those in \cite{P2}
and \cite{Sh2}.

\begin{appendix}
\chapter {Table of integrals}

Let $ k,m \in \bbn; \;\a,\b,\g \in \bbc; \; d=(m+1)/2.$ We recall
that
 $\S_m, \p$, and $\cpm$ denote the set of all  $m \times m$ real symmetric
 matrices, the cone of positive definite
matrices in  $\S_m$, and the closed cone of positive semi-definite
matrices in $\S_m$, respectively. The following formulas hold.

 \begin{equation}\label{2.7}
 \intl_a^b |r-a|^{\a -d} |b-r|^{\beta -d} dr= B_m (\a
,\b) |b-a|^{\a+\beta -d},
 \end{equation}
 $$ a \in \S_m,\quad  b>a, \quad  Re \, \a >d-1,  \quad Re \, \b
>d-1;$$

    \begin{equation}\label{2.8}
 \intl_c^b |b-r|^{\a -d} |r-c|^{\beta -d}
\frac{dr}{|r|^{\a+\beta}}= \frac{B_m (\a ,\b)}{|a|^\a |b|^\beta}
|b-a|^{\a+\beta -d},
 \end{equation}
$$ a \in \S_m,\quad  b>a,\quad  c=b^{1/2}a^{1/2}b^{-1}a^{1/2}b^{1/2},\quad  Re \, \a >d-1,  \quad Re \, \b
>d-1;$$

  \begin{equation}\label{2.9}
 \intl_{I_m}^b |b-r|^{\a -d} |r-I_m|^{\beta -d}
\frac{dr}{|r|^{\a+\beta}}= \frac{B_m (\a ,\b)}{ |b|^\beta}
|b-I_m|^{\a+\beta -d},
 \end{equation}
$$  b>I_m,\quad  Re \, \a >d-1,  \quad Re \, \b
>d-1;$$
\begin{equation}\label{2.14} \intl_s^\infty
|r|^{-\gam}|r-s|^{\a-d}dr=|s|^{\a-\gam}B_m (
\a,\gam-\a),\end{equation}
$$ s\in\p, \quad Re\,\a>d-1, \quad Re\, (\gam -\a)> d-1;$$

 \begin{equation}\label{2.13}\intl_s^\infty
|I_m+r|^{-\gam}|r-s|^{\a-d}dr=|I_m+s|^{\a-\gam}B_m (
\a,\gam-\a),\end{equation}
$$ s\in\cpm, \quad Re\,\a>d-1, \quad Re\, (\gam -\a)> d-1;$$

  \begin{equation}\label{2.15}
\intl_{\Mkm}|b+y'y|^{-\lam/2}dy=\frac{\pi^{km/2}\gm((\lam-k)/2)}{\gm(\lam/2)}|b|^{(k-\lam)/2},
 \end{equation}
$$ b\in\p, \quad Re\,\lam>k+m-1;$$

  \begin{equation}\label{2.15.1}
 \intl_{\{y\in\Mkm:\;y'y<b\}}\!\!\!\!\!\!\!
|b-y'y|^{(\lam-k)/2-d
}dy=\frac{\pi^{km/2}\gm((\lam-k)/2)}{\gm(\lam/2)}|b|^{\lam/2-d},
\end{equation}
 $$ b\in\p, \quad Re\,\lam>k+m-1.$$

\[ \text {\rm PROOF} \]

{\bf (\ref{2.7}).} We have \[ I \equiv \intl_a^b |r-a|^{\a -d}
|b-r|^{\beta -d} dr= \intl_0^c |s|^{\a -d} |c-s|^{\beta -d} \, ds,
\]
 $c= b-a$.  Let  $s=c^{1/2}tc^{1/2}$. By  Lemma \ref{12.2} (ii), \[
 I=|c|^{\a +\b-d}\intl_0^{I_m} |t|^{\a -d} |I_m
-t|^{\beta -d} \, dt=|c|^{\a +\b-d} B_m (\a ,\b). \]

{\bf (\ref{2.8}), (\ref{2.9}).}
 Let us write  (\ref{2.7}) from the right to the left, and set $$ r=a^{1/2}b^{1/2}\tau
b^{1/2}a^{1/2}, \qquad  dr=(|a||b|)^dd\tau.$$ This gives  \[ B_m
(\a ,\b)|b-a|^{\a+\beta -d} \! =\!
(|a||b|)^{\a+\beta-d}\intl_{b^{-1}}^{c^{-1}} |\tau-b^{-1}|^{\a -d}
|c^{-1}-\tau|^{\beta -d} d\tau. \] Then we set $\tau=s^{-1}, \;
d\tau=|s|^{-2d}ds, $ and get \[  B_m (\a ,\b)|b-a|^{\a+\beta -d}
\! =\!|a|^\a |b|^\beta \intl_{c}^{b} |b-s|^{\a -d} |s-c|^{\beta
-d} \frac{ds}{|s|^{\a+\beta}}\] which was required. The equality
(\ref{2.9})  follows from (\ref{2.8}).

 {\bf (\ref{2.14}), (\ref{2.13}).}
 By setting $r=q^{-1}$, $dr=|q|^{-m-1}dq$, one can write the left-hand side of
(\ref{2.14}) as \be |s|^{\a-d}\intl_0^{s^{-1}}
|q|^{\gam-\a-d}|s^{-1}-q|^{\a-d}dq\\\nonumber \stackrel{\rm
(\ref{2.7})}=|s|^{\a-\gam}B_m ( \a,\gam-\a), \ee and we are done.
The equality  (\ref{2.13}) follows from (\ref{2.14}) if we replace
$s $ and $r$ by $I_m +s$ and $I_m +r$, respectively.

 {\bf  (\ref{2.15}), (\ref{2.15.1}).}  By changing variable $y\to y b^{1/2}$, we obtain
\[
\intl_{\Mkm}|b+y'y|^{-\lam/2}dy=|b|^{(k-\lam)/2}J_1,\] \[
\intl_{\{y\in\Mkm:\;y'y<b\}}
|b-y'y|^{(\lam-k)/2-d}dy=|b|^{\lam/2-d}J_2,\] where \bea\nonumber
J_1&=&\intl_{\Mkm}|I_m+y'y|^{-\lam/2}dy,\\
J_2&=&\intl_{\{y\in\Mkm:\;y'y<I_m\}}
|I_m-y'y|^{(\lam-k)/2-d}dy.\nonumber \eea Let us show that \[
J_1=J_2=\frac{\pi^{km/2}\gm((\lam-k)/2)}{\gm(\lam/2)}.\]

 {\bf The case $k \ge m$.} We write both integrals in the polar coordinates according to Lemma
 \ref{l2.3}. For $J_1$ we have
\bea J_1&=&2^{-m}\sigma_{k,m}\intl_{\p} |r|^{k/2
-d}|I_m+r|^{-\lam/2}dr \nonumber \\
&=&2^{-m}\sigma_{k,m}B_m\left(\frac{k}{2},\frac{\lam-k}{2}\right)
\nonumber \eea (the second equality holds by (\ref{2.13}) with
$s=0, \;  \a=k/2, \; \g=\lam/2 $). Similarly,
 \bea
J_2&=&2^{-m}\sigma_{k,m}\intl_{0}^{I_m}|r|^{k/2
-d}|I_m-r|^{(\lam-k)/2-d}dr \nonumber \\
&=&2^{-m}\sigma_{k,m}B_m\left(\frac{k}{2},\frac{\lam-k}{2}\right).\nonumber
\eea Now the result follows by (\ref{2.6}) and (\ref{2.16}).

 {\bf The case $k < m$.} We replace $y$ by $y'$ and pass to the
polar coordinates. This yields
 \bea J_1&=&\intl_{\Mmk}|I_m+yy'|^{-\lam/2}dy \nonumber \\
&=&2^{-k}\intl_{\vmk}dv\intl_{\pk}
|I_m+vqv '|^{-\lam/2}|q|^{(m-k-1)/2}dq, \nonumber \\
&&(|I_m+vqv '|=|I_k+q|)\nonumber \\
&=& 2^{-k}\sigma_{m,k}
\intl_{\pk}|q|^{(m-k-1)/2}|I_k+q|^{-\lam/2}dq. \nonumber \eea By
(\ref{2.13}) (with $s=0, \;  m=k, \; \g=\lam/2 $) and
(\ref{2.5.2}), \bea
J_1&=&2^{-k}\sigma_{m,k}B_k\left(\frac{m}{2},\frac{\lam-m}{2}\right)\nonumber
\\&=& \frac{\pi^{km/2}\gk((\lam-m)/2)}{\gk(\lam/2)} \nonumber
\\&=&
\frac{\pi^{km/2}\gm((\lam-k)/2)}{\gm(\lam/2)}.\nonumber \eea
Similarly, \bea
J_2&=&\intl_{\{y\in\Mmk:\;yy'<I_m\}}|I_m-yy'|^{(\lam-k)/2-d}dy
\nonumber \\
&=&2^{-k}\intl_{\vmk}dv\intl_{\{q\in\pk:\;vqv'<I_m\}} |I_m-vqv
'|^{(\lam-k)/2-d}|q|^{(m-k-1)/2}dq \nonumber \\
&=&2^{-k}\sig_{m,k} \intl_0^{I_k}
|I_k-q|^{(\lam-m-k-1)/2}|q|^{(m-k-1)/2}dq  \nonumber \\
&=&2^{-k}\sigma_{m,k}B_k\left(\frac{m}{2},\frac{\lam-m}{2}\right),
\nonumber \eea and we get the same.

\chapter{A counterexample to Theorem  \ref{t5.1}}

Let us prove that the function
\begin{equation}\label{b1}
F(x)=F_0(x'x)=|2I_m+x'x|^{-(n+m-1)/2p}(\log|2I_m+x'x|)^{-1}\end{equation}
 belongs to $\lp$, but  \begin{equation}\label{b2} \hat
F(\xi,t)\equiv\infty \quad \text{\rm for} \quad p\geq
p_0=(n+m-1)/(k+m-1). \end{equation}  We shall put ``${}\simeq
{}$", ``${}\lesssim{}$" , and ``${} \gtrsim {}$" instead of
``$=$", ``$\le$", and ``$\ge$", respectively, if the corresponding
relation holds up to a constant multiple. To prove (\ref{b2}),
owing to (\ref{5.22}), we have
$$
\hat F(\xi,t)=\intl_{\Mkm} |2I_m+\om '\om+s|^{-(n+m-1)/2p}
(\log|2I_m+\om '\om+s|)^{-1}d\om,
$$
where $s=t't$.  This gives  ( set $\om=y(2I_m+s)^{1/2}, \quad
d\om=|2I_m+s|^{k/2}dy$)

\bea \hat F(\xi,t)&=&|2I_m+s|^{k/2-(n+m-1)/2p}I_1(z),\qquad z=\log|2I_m+s|,\nonumber
\\[14pt]
I_1(z)&=& \intl_{\Mkm}|I_m+y 'y|^{-(n+m-1)/2p} (z+\log|I_m+y
'y|)^{-1}dy.\nonumber \eea For $k \ge m$, by Lemma
 \ref{l2.3}, we obtain
$$
I_1(z) \simeq \intl_{\p}|r|^{k/2}|I_m+r|^{-(n+m-1)/2p}
(z+\log|I_m+r|)^{-1} d_\ast r\; ,$$ $d_\ast r=|r|^{-(m+1)/2}dr$.
Let us pass to polar coordinates on $\p$:
\begin{equation}\label{b3} r=\g 'a\g, \quad \g\in O(m), \quad
a=\diag(a_1,\dots, a_m), \quad a_j>0.
 \end{equation} Then
\begin{equation}\label{b4} d_\ast r=c_m \,v_m(a)\Big(\prod\limits_{j=1}^m
a_j^{-(m+1)/2}\,da_j \Big) d\g, \end{equation} where
$$
v_k(a)=\prod\limits_{1\leq i< j\leq k}|a_i-a_j|,\quad k=2,\dots,
m,\quad c_m^{-1}=\pi^{-(m^2+m)/4}\prod\limits_{j=1}^m
j\;\Gam(j/2),
$$
\cite[ p. 23, 43]{T}. By (\ref{b3}) and  (\ref{b4}),

\bea I_1(z) &\simeq & \intl_{0}^\infty\dots\intl_{0}^\infty v_m(a)
\left[z+\sum\limits_{j=1}^m \log(1+a_j)\right
]^{-1} \nonumber \\
&\times&  \!\! \left[\prod\limits_{j=1}^m a_j^{(k-m-1)/2}
(1+a_j)^{-(n+m-1)/2p}da_j\right]. \nonumber \eea The
transformation $a_j+1= b_j$ yields \be \label{b5.1}I_1(z) \gtrsim
\!\intl_{2}^3\!  b_1^ \lam db_1\!\intl_{4}^5\! b_2^ \lam db_2\!
\dots\!\intl_{2m}^{\infty}v_m(b)\left[z \! + \!
\sum\limits_{j=1}^m \log b_j\right ]^{-1} \!\!  b_m^ \lam db_m,
\ee where $\lam=(k-m-1)/2-(n+m-1)/2p$. Note that in (\ref{b5.1}),
$$
v_m(b)=\prod\limits_{1\leq i< j\leq
m}|b_j-b_i|=\prod\limits_{i=1}^ {m-1}|b_m-b_i|\prod\limits_{1\leq
i< j\leq m-1}|b_j-b_i|\gtrsim (b_m-2m)^{m-1}
$$
and $\log b_1< \ldots <\log b_m$. Hence, for $p \ge p_0$, \[
I_1(z) \gtrsim \intl_{2m}^{\infty} b_m^{\lam}
(b_m-2m)^{m-1}\left[z+m\log
 b_m\right ]^{-1} db_m =\infty.
\]
If  $k<m$ then \bea &&I_1(z)=\intl_{\Mmk} |I_m+\om
\om'|^{-(n+m-1)/2p} (z+\log|I_m+\om \om'|)^{-1}d\om \nonumber \\
&&\simeq \!\!\intl_{\vmk}\!\!dv\!\!\intl_{\pk}\!|q|^{(m-k-1)/2}
|I_m+vqv'|^{-(n+m-1)/2p} (z+\log|I_m+vqv'|)^{-1} dq.\nonumber \eea
Using the equality $|I_m+vqv'|=|I_k+q| \; $ \cite[p. 575]{Mu}, and
setting $q=\b 'b\b$, $ \; \b\in O(k)$, $ \; b=\diag(b_1,\dots,
b_k)$, $ \; b_j>0$, we have \bea I_1(z)& \simeq&  \intl_{\pk}
|q|^{(m-k-1)/2}|I_k+q|^{-(n+m-1)/2p}(z+\log|I_k+q|)^{-1} dq
\nonumber \\
& \simeq&
\intl_{0}^\infty\dots\intl_{0}^\infty\Big[\prod\limits_{j=1}^k
b_j^{(m-k-1)/2} (1+b_j)^{-(n+m-1)/2p}\Big] \nonumber
\\ &\times& \left[z+\sum\limits_{j=1}^k \log(1+b_j)\right
]^{-1}v_k(b)\,db_1\dots db_k. \nonumber \eea Proceeding as above,
for $\nu=(m-k-1)/2-(n+m-1)/2p, \; p \ge p_0,$ we obtain
$$
I_1(z)\gtrsim\intl_{2k}^{\infty} b_k^{\nu}
(b_k-2k)^{k-1}\left[z+k\log
 b_k\right ]^{-1} db_k=\infty.
$$

Let $I_2=\|F\|_p^p$. To complete the proof it remains to show that
$I_2<\infty$. As above, we have
 \bea
I_2&\simeq&\intl_{\p}|r|^{n/2}|F_0(r)|^p d_\ast r \nonumber
\\ &=&\intl_{\p}|r|^{n/2}|2I_m+r|^{-(n+m-1)/2}(\log|2I_m+r|)^{-p}
d_\ast r \nonumber \\ &\simeq&
\intl_{\bbr_+^m}\Big[\prod\limits_{j=1}^m
a_j^{n/2-(m+1)/2}(2+a_j)^{-(n+m-1)/2}\Big] \nonumber \\
&\times&\left[\sum\limits_{j=1}^m \log(2+a_j)\right
]^{-p}v_m(a)\,da_1\dots da_m, \nonumber \eea where $\bbr_+^m$ is
the set of points $a=(a_1,\dots, a_m)$ with positive coordinates.
Let us split $\bbr_+^m$ into $m+1$ pieces $\Om_0,\dots, \Om_m$,
where \bea\nonumber \Om_0&=&\{a:\; 0<a_j<1\quad \mbox{for
all}\quad
j=1,\dots,m \},\\[10pt]
\nonumber \Om_m&=&\{a:\; a_j>1\quad \mbox{for all}\quad
j=1,\dots,m   \}, \eea and $\Om_\ell$ ($\ell =1,2,\dots, m-1$) is
the set of points $a\in\bbr_+^m$ having $\ell$ coordinates $>1$
and $m-\ell$ coordinates $\le 1$.
 Then $I_2 \simeq \sum\limits_{\ell=0}^m
A_\ell$, where $A_\ell=\int_{\Om_\ell} (...) $.

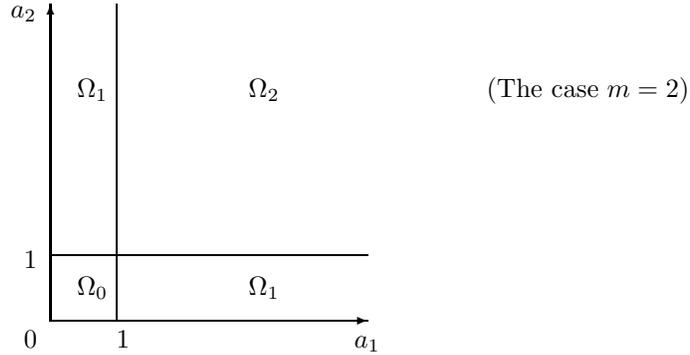
\begin{picture}(400, 150)
\put(85, 15){\vector(1,0){120}} \put(85, 15){\vector(0,1){120}}
\put(200,5){$a_1$} \put(70,130){$a_2$} \put(110,
15){\line(0,1){120}} \put(85,40){\line(1,0){120}}
\put(75,5){$0$} \put(75,35){$1$} \put(110,5){$1$}\put(160,
100){$\Om_2$}
 \put(160, 25){$\Om_1$}
 \put(95, 100){$\Om_1$}
 \put(95, 25){$\Om_0$}
 \put(250, 100){(The case $m=2$)}

\end{picture}
 For $p>1$, we have \bea  A_0 &=&
\intl_{0}^1\dots\intl_{0}^1\Big[\prod\limits_{j=1}^m
a_j^{n/2-(m+1)/2}(2+a_j)^{-(n+m-1)/2}\Big] \nonumber \\
&\times&\left[\sum\limits_{j=1}^m \log(2+a_j)\right
]^{-p}v_m(a)\,da_1\dots da_m \nonumber\\\nonumber &\lesssim&
\prod\limits_{j=1}^m \;\intl_{0}^1 a_j^{n/2-(m+1)/2}da_j<\infty
.\eea In order to estimate $A_m$, we first note that the integrand
is a symmetric function of $a_1,\dots, a_m$. Hence, $A_m=(m!)B_m$,
where $B_m$ is the integral over the set $\{a:\;a\in\Om_m; \;
a_1>a_2>\ldots >a_m\}$. Using obvious inequalities \be\label{b6.1}
a_j^{n/2-(m+1)/2}(2+a_j)^{-(n+m-1)/2}<a_j^{-m},\ee we have

\bea  B_m &\lesssim& \intl_{1}^\infty \frac{da_1} {a_1^{m}\log^p
(2+a_1)} \intl_{1}^{a_1}a_2^{-m} da_2 \dots\intl_{1}^{a_{m-2}}\;
a_{m-1}^{-m}\prod\limits_{1\leq i< j\leq m-1}|a_i-a_j|\, da_{m-1} \nonumber \\
&\times&\intl_{1}^{a_{m-1}}a_m^{-m}\prod\limits_{i=1}^{
m-1}|a_i-a_m|\,da_m. \nonumber \eea Since

\bea\nonumber \intl_{1}^{a_{m-1}}a_m^{-m}\prod\limits_{i=1}^{
m-1}|a_i-a_m|\,da_m&\leq& \Big(\prod\limits_{i=1}^{ m-1} a_i\Big)
\intl_{1}^{a_{m-1}}\frac{da_m}{a_m^{m}}\\\nonumber &=&
\frac{1-a_{m-1}^{1-m}}{m-1}\;\prod\limits_{i=1}^{ m-1} a_i
\lesssim \prod\limits_{i=1}^{ m-1} a_i,\eea it follows that

$$
B_m \lesssim \intl_{1}^\infty \frac{da_1} {a_1^{m-1}\log^p
(2+a_1)} \intl_{1}^{a_1}a_2^{1-m} da_2 \dots\intl_{1}^{a_{m-2}}\;
a_{m-1}^{1-m}\prod\limits_{1\leq i< j\leq m-1}|a_i-a_j|
\,da_{m-1}.
$$
Repeating this process, we get $$B_m \lesssim \intl_{1}^\infty
\frac{d\,a_1 }{ a_1 \log^p (2+a_1)}<\infty.$$ Let us estimate
$A_\ell$, $\ell=1,\dots,m-1$. Owing to the symmetry, it suffices
to consider the integral

\bea \nonumber B_\ell &=& \intl_{1}^\infty da_1
\intl_{1}^{a_1}da_2\dots\intl_{1}^{a_{\ell-1}}\Big[\prod\limits_{j=1}^\ell
a_j^{n/2-(m+1)/2}(2+a_j)^{-(n+m-1)/2}\Big]\\\nonumber &\times&
J(a_1,\dots, a_\ell)\; da_\ell,\eea where \bea\nonumber
J(a_1,\dots, a_\ell)&=& \intl_{0}^{1}\dots
\intl_{0}^{1}\Big[\prod\limits_{j=\ell+1}^m
a_j^{n/2-(m+1)/2}(2+a_j)^{-(n+m-1)/2}\Big]\\\nonumber
&\times&\left[\sum\limits_{j=1}^m \log(2+a_j)\right
]^{-p}v_m(a)\,da_{\ell+1}\dots da_m. \eea We  use the inequality
$$
v_m(a)=\prod\limits_{1\leq i< j\leq m}|a_i-a_j|\leq
\prod\limits_{i=1}^\ell \prod\limits_{j=i+1}^m |a_i-a_j|\leq
\prod\limits_{1\leq i< j\leq
\ell}|a_i-a_j|\prod\limits_{i=1}^\ell(a_i+1)^{m-\ell}.
$$
  This  gives  $$J(a_1,\dots, a_\ell)\le  [\log(2+a_1)]^{-p}\prod\limits_{1\leq i< j\leq
\ell}|a_i-a_j|\prod\limits_{i=1}^\ell (a_i+1)^{m-\ell},$$ and,
owing to (\ref{b6.1}),  \bea B_\ell &\lesssim& \intl_{1}^\infty
\frac{da_1} {a_1^{\ell}\log^p (2+a_1)} \intl_{1}^{a_1}a_2^{-\ell}
da_2 \dots\intl_{1}^{a_{\ell-2}}\;
a_{\ell-1}^{-\ell}\prod\limits_{1\leq i< j\leq \ell-1}|a_i-a_j| \,da_{\ell-1} \nonumber \\
&\times&\intl_{1}^{a_{\ell-1}}a_\ell^{-\ell}\prod\limits_{i=1}^{
\ell-1}|a_i-a_\ell|\,da_\ell. \nonumber \eea  Proceeding as above
(see the argument for $B_m$), we obtain
$$B_\ell \lesssim
\intl_{1}^\infty \frac{d\,a_1 }{ a_1 \log^p (2+a_1)}<\infty.$$
Thus, $I_2<\infty$, and we are done.

\chapter{Some facts from algebra}

 We recall some well-known facts    repeatedly used throughout the monograph.
We do this for convenience of the reader by taking into account
that the literature on this subject is rather sparse. More results
from matrix algebra
 can be found, i.e.,  in \cite{Mu} and \cite{FZ}.

\vskip 0.5truecm

 Notation:

 \vskip 0.5truecm

$\bullet$ $\frM_{n,m}$ is the space of real matrices $x=(x_{i,j})$
having $n$ rows and $m$
 columns.

$\bullet$ $x'$ denotes the transpose of  $x$.

$\bullet$ $\rank (x)=$ rank of $x$.

$\bullet$  $I_m$ is the identity $m \times m$
  matrix.




$\bullet$  $\tr (a)$ is the trace of the square matrix $a$.

$\bullet$ $|a|$ is  the absolute value of the determinant $\det
(a)$.


$\bullet$ $O(m)$ is the group of orthogonal
  $m \times m$ matrices.

$\bullet$ $\vnm= \{v \in \frM_{n,m}: v'v=I_m \}$
 is the Stiefel manifold. \index{a}{Stiefel manifold, $\vnm$}

\vskip 0.5truecm

Some useful facts:

\vskip 0.5truecm

{\bf 1.} If $x$ is $n \times m$ and $y$ is $m \times n$, then
$|I_n + xy|=|I_m + yx|$.

{\bf 2.} If $a$ is a square matrix and $|a|\neq 0$ then
$(a^{-1})'=(a')^{-1}$.

{\bf 3.} $\rank (x)=\rank (x')= \rank (xx')=\rank (x'x)$.

{\bf 4.}  $\rank (xy) \le \min (\rank (x), \rank (y))$.

{\bf 5.} $\rank (x+y) \le \rank (x) + \rank (y)$.

{\bf 6.} $\rank (axb)=\rank (x)$ if $a$ and $b$ are nonsingular
square matrices.

{\bf 7.} $\tr (a)=\tr (a')$.

{\bf 8.}  $\tr (a+b)= \tr (a) + \tr (b)$.

{\bf 9.}   If $x$ is $n \times m$ and $y$ is $m \times n$ then
$\tr (xy)= \tr (yx)$.

{\bf 10.} Let $\S_m$  be the space of $m \times m$ real symmetric
matrices $r=(r_{i,j}),
 \, r_{i,j}=r_{j,i}$. A matrix
 $r \in \S_m$ is called positive definite  (positive semi-definite) if $u'ru>0$ ($u'ru \ge 0$)
 for all vectors $u \neq 0$  in
  $\bbr^m$; this is commonly expressed as $r>0$ ($r \ge 0$). Given
 $r_1$ and  $r_2$ in  $\S_m$, the inequality $r_1 > r_2$ means $r_1 - r_2>0$.

{\bf (i)} If $r>0$ then $r^{-1} >0$.

 {\bf (ii)} $x'x \ge 0$ for any matrix $x$.

 {\bf (iii)} If $r \ge 0$ then $r$ is nonsingular if and only if $r>0$.

{\bf (iv)} If $r>0, \; s>0,$ and $ r-s >0$ then $s^{-1}- r^{-1}>0$
and $|r| > |s|$.

{\bf (v)} A symmetric matrix is  positive definite  (positive
semi-definite) if and only if all its eigenvalues are positive
(non-negative).

{\bf (vi)} If $r \in \S_m$ then there exists an orthogonal matrix
$\gam \in O(m)$ such that $\gam' r \gam = \Lam$ where
$\Lam=$diag$(\lam_1, \ldots , \lam_m)$. Each $\lam_j$ is real and
equal to the $j$th eigenvalue of $r$.

{\bf (vii)} If $r$ is a positive semi-definite $m \times m$ matrix
then there exists a positive semi-definite $m \times m$ matrix,
written as $r^{1/2}$, such that $r=r^{1/2}r^{1/2}$. If $ r = \gam
\Lam \gam'$, $\gam \in O(m)$, then $r^{1/2}= \gam \Lam^{1/2}
\gam'$ where $\Lam^{1/2}=$diag$(\lam_1^{1/2}, \ldots ,
\lam_m^{1/2})$.

{\bf 11.} If $x$ is $n \times k$ and $y$ is $m \times k$, $ n \ge
m$, then $x'x=y'y$ if and only if there exists $v \in V_{n,m}$
such that $x=vy$. In particular, if $x$ is $n \times m, \; n \ge
m$, then there exists $v \in V_{n,m}$ such that $x=vy$, where
$y=(x'x)^{1/2}$.


\end{appendix}

\bibliographystyle{amsalpha}

\Printindex{a}{Subject  index} \Printindex{b}{Index of symbols}

\end{document}